\documentclass{article}

\usepackage{arxiv}
\usepackage{graphicx} 
\usepackage{datetime}
\usepackage{fancyhdr}
\usepackage{CJKutf8}
\usepackage{amsmath}
\usepackage{amsfonts}
\usepackage{amsthm}
\usepackage{babel}
\usepackage{algorithm}
\usepackage{float}
\usepackage{subfig}
\usepackage{algpseudocode}
\usepackage{listings}
\usepackage{titling}
\usepackage{multirow}
\usepackage{hyperref}
\usepackage{cleveref}
\usepackage{authblk}
\allowdisplaybreaks
\newtheorem{theorem}{Theorem}[section]
\newtheorem{lemma}[theorem]{Lemma}
\newtheorem{proposition}[theorem]{Proposition}
\newtheorem{corollary}[theorem]{Corollary}

\newtheorem{remark}{Remark}[section]
\newtheorem{example}{Example}[section]

\usepackage[font=small,labelfont=bf]{caption}
\usepackage[section]{placeins}

\title{A Trace--Logarithmic Variational Functional for Equidistribution and Alignment in Moving Mesh Adaptation}

\author[1]{Wenbin Wang}
\author[2]{Yunqing Huang}
\author[3]{Huayi Wei\thanks{Corresponding author.}}
\affil[1]{Hunan Research Center of the Basic Discipline Fundamental Algorithmic Theory and Novel Computational Methods, Xiangtan University, Xiangtan 411105, Hunan, China}
\affil[2]{Hunan International Scientific and Technological Innovation Cooperation Base of Computational Science, Key Laboratory of Intelligent Computing and Information Processing of Ministry of Education, Xiangtan University, Xiangtan 411105, Hunan, China }
\affil[3]{National Center of  Applied Mathematics in Hunan, Hunan Key Laboratory for Computation and Simulation in Science and Engineering, Xiangtan 411105, Hunan, China}
\date{}

\hypersetup{
pdftitle={ Trace--Logarithmic Variational Functional for Equidistribution and Alignment in Moving Mesh Adaptation},
pdfsubject={65M50, 65N50},
pdfauthor={W. Wang, Y. Huang, H. Wei},
pdfkeywords={variational moving mesh, adaptive mesh generation, equidistribution, alignment, geometric discretization, moving mesh PDEs},
}
\begin{document}
\maketitle

\begin{abstract}
    Existing variational mesh functionals often require an empirical equidistribution--alignment weight or use a strongly nonlinear density.  
    We propose a calibrated trace--logarithmic functional in the inverse pullback tensor $\boldsymbol A=\boldsymbol J^{-1}\boldsymbol M^{-1}\boldsymbol J^{-T}$, 
    with the logarithmic coefficient fixed by the local target $\boldsymbol A=\theta\boldsymbol I$.  We prove strict convexity in the 
    symmetric positive-definite (SPD) tensor variable, inverse-Jacobian polyconvexity and coercivity, and weak minimizer existence under explicit 
    admissibility assumptions.  The density also satisfies the standard coercivity condition for semi-discrete physical-coordinate mesh 
    nonsingularity; the geometric discretization yields a compact computational-coordinate residual for the reported direct-secant implementation.  Numerical tests for metric-induced mesh adaptation, 
    a Burgers benchmark, a Rayleigh--Taylor instability simulation, and a matched high-anisotropy stress test show robust mesh redistribution; in the stress test, 
    the trace--logarithmic runs complete over a broader fixed-protocol range of the compression parameter than the scanned Huang benchmark.
\end{abstract}
\begin{keywords}{variational moving mesh ,adaptive mesh generation ,equidistribution,
    alignment,geometric discretization ,MMPDE}
\end{keywords}

\textbf{MSC codes:} 65M50, 65N50

\section{Introduction}

Moving mesh methods have been extensively developed and applied to fluid dynamics, heat transfer, combustion, and related problems \cite{tang2003adaptive,
tang2005moving,mackenzie2006moving,shen2009efficient}.  They are particularly effective when the numerical solution contains localized features, 
such as shocks, boundary layers, or evolving interfaces, for which a fixed mesh would require a costly global refinement.  By relocating nodes 
according to the solution while preserving mesh topology, moving mesh methods concentrate degrees of freedom where they are needed and retain a 
coarser resolution elsewhere.  This is especially useful for anisotropic phenomena, where the relevant length scale and preferred direction of 
resolution can differ substantially.  Recent discontinuous Galerkin and finite-element moving-mesh developments illustrate the continuing role of this strategy in 
balance-law and free-boundary computations \cite{huang2023ripa,shen2024bernoulli}.

Among the available strategies, variational moving mesh methods provide a systematic route from mesh-quality requirements to mesh equations.  
A mesh functional is constructed and then minimized, or followed through a descent flow, to determine the coordinate transformation.  This 
formulation links mesh generation with elliptic mesh equations, prescribed boundary correspondence, and finite-element implementations 
\cite{huang2010adaptive,carey1997computational,knupp2000framework,knupp1993fundamentals,winslow1981adaptive,dvinsky1991adaptive}.  Variational 
mesh adaptation is one branch of anisotropic adaptation; metric-based approaches include Hessian-driven methods, anisotropic finite elements, 
anisotropic centroidal Voronoi tessellations, and natural metrics \cite{remacle2005anisotropic,vassilevski2005hessian,apel1999anisotropic,huang2013anisotropic}.

The central geometric target in multidimensional adaptation is an \(\boldsymbol M\)-uniform mesh, whose metric tensor \(\boldsymbol M(\boldsymbol x)\) 
prescribes local scale and preferred directions.  Equidistribution controls the metric volume of the elements, whereas alignment controls their 
shape and orientation in the metric.  Their joint satisfaction defines the standard equidistribution--alignment (EA) target.  Huang's variational framework converts these conditions 
into mesh functionals and moving mesh partial differential equations (MMPDEs) \cite{huang2001variational,huang2010adaptive}; the geometric discretization of Huang and Kamenski then 
gives an elementwise differentiation framework for the resulting discrete energies \cite{huang2015geometric}.  Subsequent EA constructions extend 
these principles to surface meshes and to a common treatment of curves, surfaces, and bulk domains \cite{kolasinski2020surface,zhang2025unifying}.

Despite this structure, the choice of an EA density remains nontrivial.  Huang's classical EA functional combines alignment and equidistribution 
terms through an empirical mixing coefficient \(\mu\).  The Frobenius-norm functional of Kolasinski and Huang removes this separately prescribed 
weight by penalizing the combined EA residual directly \cite{kolasinski2018new}.  For bulk meshes in dimensions \(d=2,3\), however, its leading 
algebraic dependence on the inverse-Jacobian variables is stronger than that of the trace term used below.  This motivates a calibrated trace--logarithmic 
alternative that removes an independently prescribed EA balancing weight while retaining a comparatively simple leading trace structure.

In this work, we construct the density using the inverse pullback tensor \(\boldsymbol A\), where \(\boldsymbol A=\boldsymbol J^{-1}\boldsymbol M^{-1}\boldsymbol J^{-T}\). 
This standard EA/MMPDE variable also appears in recent curve/surface/domain formulations \cite{huang2001variational,
huang2015geometric,zhang2025unifying}; the tensor itself is not the source of novelty.  The trace--logarithmic density is calibrated by 
requiring the local EA target \(\boldsymbol A=\theta\boldsymbol I\) to be stationary, where \(\theta\) is determined by the metric volume.  
We prove strict convexity in the independent SPD tensor variable, establish a polyconvex representation in the inverse-Jacobian variable with an 
independent determinant variable, and obtain coercivity and weak minimizer existence under explicit admissibility assumptions.  The same coercive 
growth satisfies the standard hypothesis for semi-discrete physical-coordinate mesh nonsingularity.  At the discrete level, we retain the 
physical-coordinate formulation as a geometric reference and derive the computational-coordinate residual used in the reported direct-secant implementation.

The numerical experiments examine both standard adaptation behavior and a regime in which the two densities separate more clearly.  
Function-induced meshes report equidistribution, alignment, geometric quality, and interpolation behavior.  A Burgers calculation compares functional--metric combinations in a moving PDE setting, and a Rayleigh--Taylor calculation demonstrates robustness for an evolving multiscale interface.  Finally, a localized high-anisotropy stress test fixes the mesh, metric normalization, and artificial-time protocol while increasing a compression parameter, thereby probing the density response beyond standard mesh-quality diagnostics.

The remainder of this paper is organized as follows.  Section~2 reviews equidistribution and alignment principles and constructs the calibrated trace--logarithmic functional.  Section~3 establishes its analytical properties and the connection with semi-discrete mesh nonsingularity.  Section~4 presents the geometric discretization and moving mesh algorithm.  Section~5 reports the numerical experiments.

\section{Equidistribution and Alignment Conditions and Functional Construction}\label{sec:2}

An \(\boldsymbol M\)-uniform mesh is characterized by equidistribution and alignment in a prescribed symmetric positive-definite metric \(\boldsymbol M\) \cite{huang2010adaptive}.  Existing variational formulations encode these conditions through different density choices; the density determines both the analytical structure of the functional and the complexity of its geometric derivatives.  This section recalls the EA target, reviews two representative functionals, and derives the calibrated trace--logarithmic density used below.

\subsection{Equidistribution and Alignment Conditions}

Let \(\Omega_p\subset\mathbb R^d\), \(d>1\), be a bounded polygonal or polyhedral physical domain.  Let \(\boldsymbol M=\boldsymbol M(\boldsymbol x)\) be a symmetric uniformly positive-definite metric tensor on \(\Omega_p\) satisfying
\begin{equation}\label{Mbound}
m_0\boldsymbol I\le \boldsymbol M(\boldsymbol x)\le m_1\boldsymbol I,
\qquad \text{for a.e. }\boldsymbol x\in\Omega_p,
\end{equation}
where \(m_0,m_1>0\).  Let \(\Omega_c\subset\mathbb R^d\) be the computational domain, with coordinates \(\boldsymbol\xi\), and consider an orientation-preserving coordinate transformation
\[
\boldsymbol x=\boldsymbol x(\boldsymbol\xi):\Omega_c\to\Omega_p,
\qquad
\boldsymbol\xi=\boldsymbol x^{-1}(\boldsymbol x),
\qquad
\boldsymbol J=\frac{\partial\boldsymbol x}{\partial\boldsymbol\xi},
\quad \det\boldsymbol J>0.
\]

Define the metric volume
\[
\sigma:=\int_{\Omega_p}\sqrt{\det\boldsymbol M(\boldsymbol x)}\,\mathrm d\boldsymbol x
=\int_{\Omega_c}\sqrt{\det\boldsymbol M(\boldsymbol x(\boldsymbol\xi))}\det\boldsymbol J(\boldsymbol\xi)\,\mathrm d\boldsymbol\xi,
\qquad
\rho(\boldsymbol x):=\sqrt{\det\boldsymbol M(\boldsymbol x)}.
\]
Equidistribution requires
\[
\rho\det\boldsymbol J=\frac{\sigma}{|\Omega_c|}
\qquad\text{in }\Omega_c,
\]
whereas alignment requires the pullback metric to be isotropic:
\[
\boldsymbol J^T\boldsymbol M\boldsymbol J=\eta(\boldsymbol\xi)\boldsymbol I
\quad\Longleftrightarrow\quad
\boldsymbol J^{-1}\boldsymbol M^{-1}\boldsymbol J^{-T}=\theta(\boldsymbol\xi)\boldsymbol I.
\]
Together these conditions give
\begin{equation}\label{EAC}
\boldsymbol J^{-1}\boldsymbol M^{-1}\boldsymbol J^{-T}
=\theta\boldsymbol I,
\qquad
\theta:=\left(\frac{\sigma}{|\Omega_c|}\right)^{-2/d}.
\end{equation}

Let \(\mathcal T_p\) and \(\mathcal T_c\) be corresponding simplicial meshes of \(\Omega_p\) and \(\Omega_c\).  On a paired physical/computational element \(K\), the map is approximated by
\[
\boldsymbol x(\boldsymbol\xi)|_{K_c}\approx\boldsymbol J_K\boldsymbol\xi+\boldsymbol b_K.
\]
With the element-average metric \(\boldsymbol M_K\), set
\[
\rho_K:=\sqrt{\det\boldsymbol M_K},
\qquad
\sigma_h:=\sum_{K\in\mathcal T_h}|K|\rho_K.
\]
The elementwise EA target is
\[
\rho_K\det\boldsymbol J_K=\frac{\sigma_h}{|\Omega_c|},
\qquad
\boldsymbol J_K^{-1}\boldsymbol M_K^{-1}\boldsymbol J_K^{-T}=\theta_h\boldsymbol I,
\qquad
\theta_h:=\left(\frac{\sigma_h}{|\Omega_c|}\right)^{-2/d}.
\]
For a fixed discrete metric, \(\sigma_h\), \(\theta_h\), and
\[
\kappa_{\theta_h}:=p d^{p-1}\theta_h^p
\]
are global scalars rather than nonlinear unknowns.  In the epochwise algorithm of Section~4, they are recomputed at the beginning of an outer epoch and then held fixed during the associated artificial-time solve.  The density construction is independent of the particular metric-generation procedure; the examples below use the metrics specified in Section~5.

\subsection{Basic Functional Structure and Existing Functionals}

We consider functionals of the form
\[
I[\boldsymbol\xi]
=\int_{\Omega_p}G\bigl(\boldsymbol J^{-1},\det(\boldsymbol J^{-1}),\boldsymbol M,\boldsymbol x\bigr)\,\mathrm d\boldsymbol x,
\]
where \(G\) is a smooth density on the orientation-preserving class.

\begin{example}[Huang's functional]\label{ex:huang-functional}
The classical EA functional is \cite{huang2001variational}
\begin{equation}\label{HF}
I[\boldsymbol\xi]
=\mu\int_{\Omega_p}\rho\,[\operatorname{tr}(\boldsymbol J^{-1}\boldsymbol M^{-1}\boldsymbol J^{-T})]^{d\gamma/2}\,\mathrm d\boldsymbol x
+(1-2\mu)d^{d\gamma/2}\int_{\Omega_p}\frac{\rho}{[\det(\boldsymbol J)\rho]^\gamma}\,\mathrm d\boldsymbol x,
\end{equation}
where \(0<\mu\le\frac12\) and \(\gamma>0\).  The first and second terms encode alignment and equidistribution, respectively.
\end{example}

This functional is coercive and polyconvex and admits a minimizer under the standard hypotheses; see \cite{huang2010adaptive}.  Its leading trace term is comparatively moderate, while \(\mu\) remains an externally selected EA balance parameter.

\begin{example}[Kolasinski--Huang functional]\label{ex:kh-functional}
A parameter-free alternative is \cite{kolasinski2018new}
\begin{equation}\label{KHF}
I[\boldsymbol\xi]
=\int_{\Omega_p}\rho\left\|\boldsymbol J^{-1}\boldsymbol M^{-1}\boldsymbol J^{-T}-\theta\boldsymbol I\right\|_F^{2\gamma}\,\mathrm d\boldsymbol x,
\qquad
\theta=\left(\frac{\sigma}{|\Omega_c|}\right)^{-2/d}.
\end{equation}
\end{example}

The functional in \cref{KHF} follows directly from the combined EA condition.  It is coercive but nonconvex, and its algebraic nonlinearity increases with \(\gamma\) \cite{kolasinski2018new}.

Both examples are conveniently expressed through the inverse pullback tensor
\begin{equation}\label{Aalpha}
\boldsymbol A:=\boldsymbol J^{-1}\boldsymbol M^{-1}\boldsymbol J^{-T},
\qquad
 a:=\det\boldsymbol A.
\end{equation}
The determinant notation \(a\) is kept distinct from \(j:=\det\boldsymbol J\) and \(\delta:=\det\boldsymbol J^{-1}=j^{-1}\), which appear in the discrete differentiation of Section~4.  In these variables,
\begin{equation}\label{AHF}
I[\boldsymbol\xi]
=\int_{\Omega_p}\rho\left[
\mu[\operatorname{tr}(\boldsymbol A)]^{d\gamma/2}
+(1-2\mu)d^{d\gamma/2}a^{\gamma/2}
\right] \mathrm d\boldsymbol x,
\end{equation}
while \cref{KHF} becomes
\begin{equation}\label{AKHF}
I[\boldsymbol\xi]
=\int_{\Omega_p}\rho\left\|\boldsymbol A-\theta\boldsymbol I\right\|_F^{2\gamma}\,\mathrm d\boldsymbol x.
\end{equation}
The tensor \(\boldsymbol A\) is the standard inverse-pullback variable in variational EA formulations \cite{huang2001variational,huang2015geometric,zhang2025unifying}; here it provides the natural variable for calibrating and differentiating the proposed density.

\subsection{Derivation of a New Functional}

Let \(\gamma>1\) and set
\[
p:=\frac{d\gamma}{2}>\frac d2,
\qquad
\kappa_\theta:=p d^{p-1}\theta^p=d^p\frac\gamma2\theta^p.
\]
For \(\boldsymbol A\) in the cone \(\operatorname{SPD}(d)\) of \(d\times d\) symmetric positive-definite matrices, consider
\[
T_{p,\kappa}(\boldsymbol A)
:=[\operatorname{tr}(\boldsymbol A)]^p-\kappa\log\det(\boldsymbol A).
\]
Its symmetric-Frobenius derivative is
\[
\nabla_{\boldsymbol A}T_{p,\kappa}(\boldsymbol A)
=p[\operatorname{tr}(\boldsymbol A)]^{p-1}\boldsymbol I-\kappa\boldsymbol A^{-1}.
\]
Requiring the EA target \(\boldsymbol A=\theta\boldsymbol I\) to be stationary yields
\[
\left.
\nabla_{\boldsymbol A}T_{p,\kappa}(\boldsymbol A)
\right|_{\boldsymbol A=\theta\boldsymbol I}
=
\left[p(d\theta)^{p-1}-\kappa\theta^{-1}\right]\boldsymbol I
=
\boldsymbol 0,
\]
and hence
\[
\kappa=p d^{p-1}\theta^p=\kappa_\theta.
\]
Thus the logarithmic coefficient is fixed by the target scale.  With \(a=\det\boldsymbol A\), define
\begin{equation}\label{GTF}
G(\boldsymbol A,a,\boldsymbol M)=\rho\,T(\boldsymbol A,a),
\qquad
T(\boldsymbol A,a):=[\operatorname{tr}(\boldsymbol A)]^p-\kappa_\theta\log a.
\end{equation}
The proposed mesh functional is therefore
\begin{equation}\label{PF}
I[\boldsymbol\xi]
=\int_{\Omega_p}\rho\left([
\operatorname{tr}(\boldsymbol A)]^p-\kappa_\theta\log\det\boldsymbol A\right)\,\mathrm d\boldsymbol x,
\qquad
\theta=\left(\frac{\sigma}{|\Omega_c|}\right)^{-2/d}.
\end{equation}
For \(d=2,3\), the leading growth of the trace term in \(\boldsymbol P=\boldsymbol J^{-1}\) is \(\|\boldsymbol P\|_F^{d\gamma}\), whereas the leading growth of \cref{KHF} is \(\|\boldsymbol P\|_F^{4\gamma}\).  This difference motivates the trace--logarithmic alternative.

\section{Properties of the Proposed Functional}\label{sec:3}

This section identifies the calibrated tensor target and then establishes polyconvexity, coercivity, and weak existence in the inverse-Jacobian representation.  The resulting coercive growth is finally connected with the standard physical-coordinate semi-discrete mesh-nonsingularity result.

\subsection{Scale Invariance}

\begin{proposition}[Scale invariance of the minimizer]\label{thm:scaleinvariance}
Let \(\widetilde{\boldsymbol M}=c\boldsymbol M\) with \(c>0\). Then
\[
I[\boldsymbol\xi;c\boldsymbol M]
=c^{\frac d2-p}\left(
I[\boldsymbol\xi;\boldsymbol M]
+d\kappa_\theta\sigma\log c
\right).
\]
Consequently, over any fixed admissible class independent of \(c\),
\[
\operatorname*{argmin}_{\boldsymbol\xi}I[\boldsymbol\xi;c\boldsymbol M]
=
\operatorname*{argmin}_{\boldsymbol\xi}I[\boldsymbol\xi;\boldsymbol M].
\]
\end{proposition}

\begin{proof}
Under \(\boldsymbol M\mapsto c\boldsymbol M\),
\[
\widetilde\rho=c^{d/2}\rho,\qquad
\widetilde{\boldsymbol A}=c^{-1}\boldsymbol A,\qquad
\widetilde\theta=c^{-1}\theta,\qquad
\kappa_{\widetilde\theta}=c^{-p}\kappa_\theta.
\]
Hence
\[
\widetilde\rho\,T_{\widetilde\theta}(\widetilde{\boldsymbol A})
=c^{\frac d2-p}\rho
\left[T_\theta(\boldsymbol A)+d\kappa_\theta\log c\right].
\]
Integration over \(\Omega_p\), together with \(\sigma=\int_{\Omega_p}\rho\,\mathrm d\boldsymbol x\), proves the identity.  The second term is independent of \(\boldsymbol\xi\), and the prefactor is positive.
\end{proof}

\subsection{Calibrated Tensor Characterization and Polyconvexity}

\begin{theorem}[Calibrated tensor target]\label{tensorconvexity}
Let \(p=d\gamma/2\), \(\gamma>1\), and \(\kappa_\theta=pd^{p-1}\theta^p\).  For \(\boldsymbol A\in\operatorname{SPD}(d)\), define
\[
T_\theta(\boldsymbol A)
:=[\operatorname{tr}(\boldsymbol A)]^p-\kappa_\theta\log\det\boldsymbol A.
\]
Then \(T_\theta\) is strictly convex on \(\operatorname{SPD}(d)\), and its unique minimizer is \(\theta\boldsymbol I\).
\end{theorem}

\begin{proof}
For \(\boldsymbol H\in\operatorname{Sym}(d)\),
\[
D^2T_\theta(\boldsymbol A)[\boldsymbol H,\boldsymbol H]
=p(p-1)[\operatorname{tr}(\boldsymbol A)]^{p-2}[\operatorname{tr}(\boldsymbol H)]^2
+\kappa_\theta\left\|\boldsymbol A^{-1/2}\boldsymbol H\boldsymbol A^{-1/2}\right\|_F^2.
\]
The second term is positive for every nonzero \(\boldsymbol H\), so \(T_\theta\) is strictly convex.  Moreover,
\[
DT_\theta(\theta\boldsymbol I)[\boldsymbol H]
=\left[p(d\theta)^{p-1}-\kappa_\theta\theta^{-1}\right]\operatorname{tr}(\boldsymbol H)=0.
\]
Thus \(\theta\boldsymbol I\) is the unique minimizer.
\end{proof}

\begin{remark}[Target gap]
A direct rearrangement gives
\[
T_\theta(\boldsymbol A)-T_\theta(\theta\boldsymbol I)
=d^p\theta^p\left[
\left(\frac{\operatorname{tr}\boldsymbol A}{d\theta}\right)^p-1
-p\log\left(\frac{\operatorname{tr}\boldsymbol A}{d\theta}\right)
\right]
+\kappa_\theta\log\left(
\frac{(\operatorname{tr}\boldsymbol A/d)^d}{\det\boldsymbol A}
\right).
\]
The first term measures departure from the target trace scale, while the second is nonnegative by the arithmetic--geometric mean inequality and vanishes precisely for isotropic \(\boldsymbol A\).
\end{remark}

\begin{proposition}[Polyconvexity]\label{Polyc}
For almost every \(\boldsymbol x\in\Omega_p\), define the extended integrand
\[
W_{\boldsymbol x}(\boldsymbol P):=
\begin{cases}
G\!\left(\boldsymbol P\boldsymbol M^{-1}\boldsymbol P^T,
\det(\boldsymbol P\boldsymbol M^{-1}\boldsymbol P^T),\boldsymbol M\right),
&\det\boldsymbol P>0,\\
+\infty,&\det\boldsymbol P\le0.
\end{cases}
\]
Then \(W_{\boldsymbol x}\) is polyconvex with respect to
\(\boldsymbol P=\boldsymbol J^{-1}\).
\end{proposition}
\begin{proof}
Let \(q:=2p=d\gamma>1\).  For \(\det\boldsymbol P>0\),
\[
\left[\operatorname{tr}(\boldsymbol P\boldsymbol M^{-1}\boldsymbol P^T)\right]^p
=\|\boldsymbol P\boldsymbol M^{-1/2}\|_F^q,
\qquad
\det(\boldsymbol P\boldsymbol M^{-1}\boldsymbol P^T)
=\det(\boldsymbol M^{-1})(\det\boldsymbol P)^2.
\]
Consequently,
\[
W_{\boldsymbol x}(\boldsymbol P)
=\overline G_{\boldsymbol x}(\boldsymbol P,\det\boldsymbol P),
\]
where
\[
\overline G_{\boldsymbol x}(\boldsymbol Y,r):=
\begin{cases}
\rho(\boldsymbol x)\!\left[
\|\boldsymbol Y\boldsymbol M(\boldsymbol x)^{-1/2}\|_F^q
-2\kappa_\theta\log r
+\kappa_\theta\log\det\boldsymbol M(\boldsymbol x)
\right],&r>0,\\
+\infty,&r\le0.
\end{cases}
\]
Since \(q>1\), the map
\(\boldsymbol Y\mapsto\|\boldsymbol Y\boldsymbol M^{-1/2}\|_F^q\)
is convex as the composition of a convex norm power with a linear map.  Moreover,
\[
r\longmapsto
\begin{cases}
-2\kappa_\theta\log r,&r>0,\\
+\infty,&r\le0
\end{cases}
\]
is convex and lower semicontinuous on \(\mathbb R\).  Thus
\(\overline G_{\boldsymbol x}\) is convex in the independent variables
\((\boldsymbol Y,r)\).

Let \(\mathcal M_k(\boldsymbol P)\) denote the vector of all \(k\)-th order minors and write
\[
\mathcal M(\boldsymbol P)
=\bigl(\boldsymbol P,\mathcal M_2(\boldsymbol P),\ldots,
\mathcal M_{d-1}(\boldsymbol P),\det\boldsymbol P\bigr).
\]
Define on the full minor space
\[
H_{\boldsymbol x}(\boldsymbol Y,\boldsymbol Z_2,\ldots,
\boldsymbol Z_{d-1},r):=\overline G_{\boldsymbol x}(\boldsymbol Y,r).
\]
The function \(H_{\boldsymbol x}\) is convex and is simply independent of the intermediate-minor variables.  By the extended-value definition,
\[
W_{\boldsymbol x}(\boldsymbol P)
=H_{\boldsymbol x}(\mathcal M(\boldsymbol P))
\qquad\text{for every }\boldsymbol P\in\mathbb R^{d\times d},
\]
which is the defining polyconvex representation.
\end{proof}

\subsection{Coercivity in the Inverse-Jacobian Variable and Weak Existence}

Set \(q:=2p=d\gamma\) and
\[
\mathcal A_g^0:=\{\boldsymbol\xi\in W^{1,q}(\Omega_p;\mathbb R^d):
\operatorname{Tr}\boldsymbol\xi=\boldsymbol g\}.
\]
Here \(\operatorname{Tr}\) denotes the Sobolev trace operator.  The following scalar estimate isolates the logarithmic contribution.

\begin{lemma}[Scalar logarithmic lower envelope]\label{scalarenvelope}
For \(b>0\), \(c\in(0,1)\), and \(x>0\),
\begin{equation}\label{scalarenvelopeineq}
x-b\log x\ge cx-b\left[\log\!\left(\frac{b}{1-c}\right)-1\right]_+,
\qquad [s]_+:=\max\{s,0\}.
\end{equation}
\end{lemma}
\begin{proof}
Set \(\phi_{b,c}(x):=(1-c)x-b\log x\).  Then
\(\phi_{b,c}''(x)=b/x^2>0\), and its unique minimizer is
\(x_*=b/(1-c)\).  Hence
\[
\phi_{b,c}(x)\ge
b\left[1-\log\!\left(\frac{b}{1-c}\right)\right],
\]
and replacing a positive lower constant by zero when necessary gives
\eqref{scalarenvelopeineq}.
\end{proof}

\begin{proposition}[Coercivity]\label{coercive}
Assume \eqref{Mbound} and
\(\boldsymbol g\in W^{1-1/q,q}(\partial\Omega_p;\mathbb R^d)\).  Set
\[
\rho_0:=m_0^{d/2},\qquad \rho_1:=m_1^{d/2},\qquad
C_{c,\theta}:=d^p\theta^p
\left[\log\!\left(\frac{\theta^p}{1-c}\right)-1\right]_+,
\quad c\in(0,1).
\]
For almost every \(\boldsymbol x\in\Omega_p\), every
\(\boldsymbol A\in\operatorname{SPD}(d)\), and every
\(\boldsymbol P\in\mathbb R^{d\times d}\) with \(\det\boldsymbol P>0\),
\begin{align}
G(\boldsymbol A,\det\boldsymbol A,\boldsymbol M)
&\ge c\rho_0[\operatorname{tr}(\boldsymbol A)]^p-\rho_1C_{c,\theta},\label{coerciveA}\\
G(\boldsymbol P\boldsymbol M^{-1}\boldsymbol P^T,
\det(\boldsymbol P\boldsymbol M^{-1}\boldsymbol P^T),\boldsymbol M)
&\ge c\rho_0m_1^{-p}\|\boldsymbol P\|_F^q-\rho_1C_{c,\theta}.\label{coerciveP}
\end{align}
Consequently, for each fixed \(c\in(0,1)\), there are \(c_I>0\) and
\(C_I\ge0\), depending only on \(c,q,m_0,m_1,\theta,\Omega_p\) and
\(\|\boldsymbol g\|_{W^{1-1/q,q}(\partial\Omega_p)}\), such that
\begin{equation}\label{coerciveI}
I[\boldsymbol\xi]\ge c_I\|\boldsymbol\xi\|_{W^{1,q}(\Omega_p)}^q-C_I,
\qquad \boldsymbol\xi\in\mathcal A_g^0.
\end{equation}
\end{proposition}
\begin{proof}
Let
\[
s:=\frac{\operatorname{tr}\boldsymbol A}{d},
\qquad x:=s^p,
\qquad b_\theta:=\theta^p.
\]
The arithmetic--geometric mean inequality gives
\(\det\boldsymbol A\le s^d\), and therefore
\(\log\det\boldsymbol A\le d\log s\).  Since
\(\kappa_\theta=pd^{p-1}\theta^p\),
\begin{align*}
T_\theta(\boldsymbol A)
&=[\operatorname{tr}(\boldsymbol A)]^p
-\kappa_\theta\log\det\boldsymbol A\\
&\ge d^ps^p-pd^p\theta^p\log s
=d^p\bigl(x-b_\theta\log x\bigr).
\end{align*}
Applying Lemma~\ref{scalarenvelope} with \(b=b_\theta\) yields
\[
T_\theta(\boldsymbol A)
\ge c[\operatorname{tr}(\boldsymbol A)]^p-C_{c,\theta}.
\]
The metric bounds imply \(\rho_0\le\rho\le\rho_1\).  Multiplying the preceding estimate by \(\rho\), using the lower bound for the positive term and the upper bound for the negative constant, proves \eqref{coerciveA}.

For \(\boldsymbol A=\boldsymbol P\boldsymbol M^{-1}\boldsymbol P^T\),
\[
\operatorname{tr}\boldsymbol A
=\|\boldsymbol P\boldsymbol M^{-1/2}\|_F^2
\ge m_1^{-1}\|\boldsymbol P\|_F^2.
\]
Raising this inequality to the power \(p\), and using \(q=2p\) in
\eqref{coerciveA}, gives \eqref{coerciveP}.

Let \(\boldsymbol\xi\in\mathcal A_g^0\).  If
\(I[\boldsymbol\xi]=+\infty\), then \eqref{coerciveI} is immediate.
Otherwise the extended-value definition gives
\(\det\nabla\boldsymbol\xi>0\) almost everywhere, so integration of
\eqref{coerciveP} with \(\boldsymbol P=\nabla\boldsymbol\xi\) gives
\[
I[\boldsymbol\xi]
\ge c\rho_0m_1^{-p}\|\nabla\boldsymbol\xi\|_{L^q(\Omega_p)}^q
-\rho_1C_{c,\theta}|\Omega_p|.
\]
By the trace theorem, choose an extension
\(\boldsymbol\xi_g\in W^{1,q}(\Omega_p;\mathbb R^d)\) satisfying
\[
\operatorname{Tr}\boldsymbol\xi_g=\boldsymbol g,
\qquad
\|\boldsymbol\xi_g\|_{W^{1,q}(\Omega_p)}
\le C_E\|\boldsymbol g\|_{W^{1-1/q,q}(\partial\Omega_p)}.
\]
Since \(\boldsymbol\xi-\boldsymbol\xi_g\in W_0^{1,q}\), Poincar\'e's inequality and the elementary power inequality imply
\[
\|\boldsymbol\xi\|_{W^{1,q}(\Omega_p)}^q
\le C\|\nabla\boldsymbol\xi\|_{L^q(\Omega_p)}^q+C_g,
\]
where \(C_g\) depends only on \(\Omega_p\) and the trace norm of
\(\boldsymbol g\).  Combining the last two estimates proves
\eqref{coerciveI}.
\end{proof}

\begin{theorem}[Existence of a weak minimizer with positive Jacobian a.e.]
\label{existence}
Let \(q=d\gamma>d\), let
\(\boldsymbol M\in L^\infty(\Omega_p;\mathrm{SPD}(d))\)
satisfy \cref{Mbound}, and let
\(
\boldsymbol g\in
W^{1-1/q,q}(\partial\Omega_p;\mathbb R^d).
\)
Interpret \(I\) on \(\mathcal A_g^0\) through the
orientation-preserving extended integrand in \cref{Polyc}. If
\(\mathcal A_g^0\) contains a finite-energy competitor, then \(I\)
admits a minimizer
\(\boldsymbol\xi_*\in\mathcal A_g^0\) satisfying
\[
\det\nabla\boldsymbol\xi_*>0
\qquad\text{a.e. in }\Omega_p.
\]
\end{theorem}

\begin{proof}
Let \((\boldsymbol\xi_n)\subset\mathcal A_g^0\) be a finite-energy
minimizing sequence. By \cref{coercive}, it is bounded in
\(W^{1,q}(\Omega_p;\mathbb R^d)\). Hence, after extraction,
\[
\boldsymbol\xi_n\rightharpoonup\boldsymbol\xi_*
\quad\text{in }W^{1,q}(\Omega_p;\mathbb R^d).
\]
Weak continuity of the trace operator gives
\(\boldsymbol\xi_*\in\mathcal A_g^0\). Since \(q>d\), weak
continuity of minors
\cite{ball1976convexity,dacorogna2008direct} yields
\[
\det\nabla\boldsymbol\xi_n
\rightharpoonup
\det\nabla\boldsymbol\xi_*
\quad\text{in }L^{q/d}(\Omega_p).
\]

Define
\[
\widehat G_{\boldsymbol x}(\boldsymbol P,r)
:=
\overline G_{\boldsymbol x}(\boldsymbol P,r)
+2\rho(\boldsymbol x)\kappa_\theta r.
\]
For \(r>0\),
\[
\widehat G_{\boldsymbol x}(\boldsymbol P,r)
=
\rho(\boldsymbol x)
\left[
\|\boldsymbol P\boldsymbol M^{-1/2}\|_F^q
+2\kappa_\theta(r-\log r)
+\kappa_\theta\log\det\boldsymbol M
\right],
\]
and \(\widehat G_{\boldsymbol x}=+\infty\) for \(r\le0\).
Since \(r-\log r\ge1\) and
\(\log\det\boldsymbol M\in L^\infty(\Omega_p)\),
\(\widehat G_{\boldsymbol x}\) is a convex lower-semicontinuous
normal integrand with an integrable lower bound. Therefore, weak
lower semicontinuity of convex integral functionals
\cite{dacorogna2008direct} gives
\[
\int_{\Omega_p}
\widehat G_{\boldsymbol x}
(\nabla\boldsymbol\xi_*,
 \det\nabla\boldsymbol\xi_*)\,d\boldsymbol x
\le
\liminf_{n\to\infty}
\int_{\Omega_p}
\widehat G_{\boldsymbol x}
(\nabla\boldsymbol\xi_n,
 \det\nabla\boldsymbol\xi_n)\,d\boldsymbol x.
\]
Moreover,
\(\rho\kappa_\theta\in L^\infty(\Omega_p)
\subset L^{(q/d)'}(\Omega_p)\), so the added determinant term is
weakly continuous. Subtracting it yields
\[
I[\boldsymbol\xi_*]
\le
\liminf_{n\to\infty}I[\boldsymbol\xi_n]
=
\inf_{\boldsymbol\xi\in\mathcal A_g^0}I[\boldsymbol\xi].
\]
Thus \(\boldsymbol\xi_*\) is a minimizer. Its energy is finite, and
the extended-value definition therefore implies
\[
\det\nabla\boldsymbol\xi_*>0
\qquad\text{a.e. in }\Omega_p.
\]
\end{proof}

\subsection{Mesh Nonsingularity and the Inverse-Coordinate Representation}

Fix an inner epoch so that \(I_h\), the metric and scale data, the computational reference mesh, the boundary correspondence, and the positive weights are independent of artificial time.  For the free-node set \(\mathcal N_{\rm free}\), consider
\begin{equation}\label{SEMPDE}
\frac{\mathrm d\boldsymbol x_i}{\mathrm dt}
=-\frac{\omega_i}{\tau}\frac{\partial I_h}{\partial\boldsymbol x_i},
\qquad i\in\mathcal N_{\rm free},
\qquad 0<\omega_{\min}\le\omega_i\le\omega_{\max}.
\end{equation}
Here \(\partial I_h/\partial\boldsymbol x_i\) is the full derivative, including the explicit spatial dependence of a fixed metric field sampled at moving nodes.

\begin{corollary}[Semi-discrete mesh nonsingularity]\label{nonsingularity}
Assume \eqref{Mbound}, a nonsingular initial physical mesh, a fixed nonsingular computational mesh, 
an admissible boundary correspondence, and the smoothness assumptions on the fixed metric/density required in \cite{huang2016nonsingularity}.  
Then \eqref{SEMPDE} remains nonsingular for all \(t>0\).  
There are \(c_h,c_v>0\), depending only on the fixed problem and computational-mesh data, the number of elements, 
the coercivity constants, and \(I_h(\mathcal T_h(0))\), such that
\[
a_{K,\boldsymbol M}(t)\ge c_h,
\qquad |K(t)|\ge c_v,
\qquad K\in\mathcal T_h(t),\quad t>0.
\]
Here \(a_{K,\boldsymbol M}\) is the minimum metric altitude.  If the epoch functional is time independent, every accumulation mesh is nonsingular and free-node stationary.
\end{corollary}
\begin{proof}
Along \eqref{SEMPDE},
\[
\frac{\mathrm d I_h}{\mathrm dt}
=-\frac1\tau\sum_{i\in\mathcal N_{\rm free}}\omega_i
\left\|\frac{\partial I_h}{\partial\boldsymbol x_i}\right\|_2^2\le0.
\]
By \eqref{coerciveA}, with exponent \(p=d\gamma/2>d/2\), \(\alpha:=c\rho_0\), and \(\beta:=\rho_1C_{c,\theta}\),
\[
G(\boldsymbol P\boldsymbol M^{-1}\boldsymbol P^T,
\det(\boldsymbol P\boldsymbol M^{-1}\boldsymbol P^T),\boldsymbol M)
\ge\alpha[\operatorname{tr}(\boldsymbol P\boldsymbol M^{-1}\boldsymbol P^T)]^{p}-\beta.
\]
For fixed epoch data, the trace--logarithmic density is smooth on the nonsingular element-configuration set, and the displayed bound is exactly the coercive inverse-Jacobian condition in \cite[Theorems~4.1 and~4.3]{huang2016nonsingularity}.  Restricting the gradient sums to free nodes leaves the energy identity and admissible variations unchanged, so those arguments give the altitude and volume bounds and, for a time-independent epoch, free-node stationarity of every accumulation mesh.
\end{proof}

\begin{remark}
For every admissible element,
\[
\boldsymbol P_K=\boldsymbol J_K^{-1},
\qquad \det\boldsymbol P_K=(\det\boldsymbol J_K)^{-1}.
\]
Thus both views encode the same orientation, while \(\det\boldsymbol J_K\to0^+\) forces \(\|\boldsymbol P_K\|_F\to\infty\); coercivity therefore penalizes vanishing physical scales, and polyconvexity supplies weak lower semicontinuity.  These continuum properties neither identify the \(\boldsymbol x\)- and \(\boldsymbol\xi\)-view trajectories nor imply fully discrete nonsingularity.  Section~4 uses the same density in computational coordinates, and each metric update defines a new fixed inner-epoch functional.
\end{remark}

\section{Geometric Discretization for Functionals and Moving Mesh Algorithms}\label{sec:4}

This section gives the geometric discretization and the associated moving-mesh algorithm.  The physical-coordinate differentiation follows the established geometric MMPDE construction \cite{huang2001variational,huang2015geometric,zhang2025unifying}; it provides the semi-discrete reference formulation used in Section~3.  Related variational MMPDE constructions have also been developed for parallel tetrahedral mesh-quality improvement \cite{lopez2022parallel}.  The computational-coordinate ($\boldsymbol\xi$-view) formulas use the same elementwise inverse geometry $\boldsymbol P_K=\boldsymbol J_K^{-1}$ and assemble the residual used in the reported solver.

\subsection{Discrete Energy Functional and Gradient Flow}

After mesh discretization, the functional \cref{PF} can be written as
\begin{equation}
I_h = \sum_{K\in \mathcal{T}_h}|K| G( \boldsymbol A_K , a_K , \boldsymbol M_K) = \sum_{K\in \mathcal{T}_h}|K| \rho_K T_K,
\end{equation}
where
\[
\begin{aligned}
\boldsymbol P_K&:=\boldsymbol J_K^{-1},
& j_K&:=\det\boldsymbol J_K,
& \delta_K&:=\det\boldsymbol P_K=j_K^{-1},\\
\boldsymbol A_K&:=\boldsymbol P_K\boldsymbol M_K^{-1}\boldsymbol P_K^T,
& a_K&:=\det\boldsymbol A_K
=\det(\boldsymbol M_K^{-1})\delta_K^2.
\end{aligned}
\]
Here \(\delta_K\) is the determinant variable in the geometric differentiation, while \(a_K\) is the determinant argument of the tensor density \(T\).  We write \(T_K:=T(\boldsymbol A_K,a_K)\) and \(G_K:=\rho_KT_K\) for the cellwise quadrature density.  The element contribution to \(I_h\) is always \(|K|G_K\).  In a \(\boldsymbol\xi\)-variation the factor \(|K|\) is fixed; in a \(\boldsymbol x\)-variation the complete summand \(|K|G_K\) must be differentiated.
Here, $\boldsymbol J_K$ is the element Jacobian.  For the $\boldsymbol x$-view reference formula, we use the standard piecewise-affine metric interpolation \cite[Appendix~A]{huang2015geometric},
\[
\boldsymbol M_h|_K=\sum_{j=0}^{d}\phi_{j,K}\boldsymbol M_{j,K},
\qquad
\boldsymbol M_K=\frac1{d+1}\sum_{j=0}^{d}\boldsymbol M_{j,K},
\qquad
\boldsymbol M_{j,K}=\boldsymbol M_h(\boldsymbol x_j^K).
\]
Thus \(\boldsymbol M_{j,K}\) is the nodal metric sample and \(\boldsymbol M_K\) its cell average.  At outer epoch $n$, let $\boldsymbol M_i^{(n)}$ denote the nodal metric.  We write \(I_h^{(n)}\) for the discrete functional evaluated with the resulting metric and the associated scale \(\theta_h^{(n)}\); the epoch superscript is suppressed in the local formulas below unless it is needed to distinguish stored data.  Define
\begin{equation}\label{eq:metric-balance}
[\boldsymbol M]_i^{(n)}
:=\left[(\theta_h^{(n)})^{-1/2}
\det(\boldsymbol M_i^{(n)})^{1/d}\right]^{1/2},
\qquad
\omega_i^{\xi,(n)}
:=\bigl([\boldsymbol M]_i^{(n)}\bigr)^{-d/2}.
\end{equation}
This balances metric scale in the computational-coordinate velocity; a common positive multiplier is absorbed into $\tau$.  The values $\omega_i^{\xi,(n)}$ are frozen during the associated inner solve and therefore do not enter the stage tangent.

Let \(N_{\rm v}\) be the number of mesh vertices.  The evolution of the mesh in the $\xi$-view is governed by
\begin{equation}\label{eq:xi-view-flow}
\frac{\partial \boldsymbol\xi_i}{\partial t}
=-\frac{\omega_i^{\xi,(n)}}{\tau}
\frac{\partial I_h}{\partial\boldsymbol\xi_i},
\qquad i=0,\ldots,N_{\rm v}-1.
\end{equation}
For frozen epoch data, this is a descent system for \(I_h^{(n)}\).  The physical-coordinate reference system below is the form used for the fixed-metric semi-discrete result in Section~3.  Let \(\ell(i,K)\in\{0,\ldots,d\}\) be the local index for which \(\boldsymbol x_i=\boldsymbol x_{\ell(i,K)}^K\). 
 Since \(|K|\) is fixed in a \(\boldsymbol\xi\)-variation, the global gradient is assembled as
\begin{equation}\label{xiview_equ}
\frac{\partial I_h}{\partial\boldsymbol\xi_i}
=\sum_{K\in S_i}|K|\frac{\partial G_K}{\partial\boldsymbol\xi_{\ell(i,K)}^K},
\qquad i=0,\ldots,N_{\rm v}-1,
\end{equation}
where \(S_i\) denotes the element star associated with node \(i\).

Alternatively, one may directly evolve the physical mesh points without explicitly tracking the mapping from computational to physical coordinates,
\begin{equation}\label{eq:x-view-flow}
\frac{\partial\boldsymbol x_i}{\partial t}
=-\frac{\omega_i}{\tau}\frac{\partial I_h}{\partial\boldsymbol x_i},
\qquad i=0,\ldots,N_{\rm v}-1.
\end{equation}
The $\boldsymbol x$-view is retained as the physical-coordinate reference formulation.  During an inner epoch, \(\theta_h^{(n)}\) is fixed.  Hence $\partial G_K/\partial\boldsymbol x$ below denotes only explicit spatial dependence with $(\boldsymbol P_K,\delta_K,\boldsymbol M_K,\theta_h^{(n)})$ fixed; metric variation remains in the $\partial G_K/\partial\boldsymbol M_K$ term.

Unlike the \(\boldsymbol\xi\)-view, an \(\boldsymbol x\)-variation differentiates the whole summand \(|K|G_K\).  We therefore define the normalized physical-coordinate descent vectors by
\[
\boldsymbol v_j^K:=-\frac1{|K|}\frac{\partial(|K|G_K)}{\partial\boldsymbol x_j^K},
\qquad j=0,\ldots,d,
\]
and assemble
\begin{equation}\label{xview_equ}
\frac{\partial I_h}{\partial\boldsymbol x_i}
=\sum_{K\in S_i}\frac{\partial(|K|G_K)}{\partial\boldsymbol x_{\ell(i,K)}^K}
=-\sum_{K\in S_i}|K|\boldsymbol v_{\ell(i,K)}^K,
\qquad i=0,\ldots,N_{\rm v}-1.
\end{equation}
Thus \(\boldsymbol v_j^K\) is not \(-\partial G_K/\partial\boldsymbol x_j^K\): it also contains the variation of \(|K|\).

In the following subsections, we first reduce the element gradients in \cref{xiview_equ} and \cref{xview_equ} for a generic $T(\boldsymbol A,a)$, then specialize $T$ using \cref{GTF}.

\subsection{Geometric Discretization of Gradient}

Now, the vertex coordinates  $\boldsymbol x_0^K , \cdots, \boldsymbol x_d^K$  of cell  $K$  and the vertex coordinates $\boldsymbol \xi_0^K , \cdots, \boldsymbol \xi_d^K$ of  $\widehat K$  satisfy the following relationship:
\[\boldsymbol x_i^K - \boldsymbol x_0^K = \boldsymbol J_K(\boldsymbol\xi_i^K - \boldsymbol\xi_0^K) , i = 1,
\cdots, d\]

Let
\[
E_K=
[\boldsymbol x_1^K-\boldsymbol x_0^K,\ldots,
 \boldsymbol x_d^K-\boldsymbol x_0^K],
\qquad
\widehat E_K=
[\boldsymbol\xi_1^K-\boldsymbol\xi_0^K,\ldots,
 \boldsymbol\xi_d^K-\boldsymbol\xi_0^K]
\]
be the physical and computational edge matrices of element \(K\).  With
\(\delta_K=\det(\boldsymbol J_K^{-1})\), equivalently
\(\delta_K=\det(\widehat E_K)/\det(E_K)\), the standard geometric identity for the
\(d\) non-base vertices is
\begin{equation}\label{HuangDiscreteXi}
\frac{\partial G_K}
{\partial[\boldsymbol\xi_1^K,\ldots,\boldsymbol\xi_d^K]}
=
E_K^{-1}\frac{\partial G_K}{\partial\boldsymbol J_K^{-T}}
+
\frac{\partial G_K}{\partial\delta_K}\,
\frac{\det(\widehat E_K)}{\det(E_K)}\,\widehat E_K^{-1}.
\end{equation}
The base-vertex row is fixed by translation invariance:
\begin{equation}\label{HuangDiscreteXi0}
\frac{\partial G_K}{\partial\boldsymbol\xi_0^K}
=
-\boldsymbol e_d^T
\frac{\partial G_K}
{\partial[\boldsymbol\xi_1^K,\ldots,\boldsymbol\xi_d^K]},
\qquad
\boldsymbol e_d=(1,\ldots,1)^T\in\mathbb R^d.
\end{equation}

For the physical-coordinate formulation we use the covariant geometric differentiation of Huang--Kamenski \cite[Appendix~A]{huang2015geometric}.  Their formula differentiates the complete contribution \(|K|G_K\); in the present notation it gives the normalized descent vectors \(\boldsymbol v_j^K=-|K|^{-1}\partial(|K|G_K)/\partial\boldsymbol x_j^K\).  For the $d$ non-base vertices,
\begin{equation}\label{HuangDiscreteX}
\begin{aligned}
\begin{bmatrix}
(\boldsymbol v_1^K)^T\\ \vdots\\ (\boldsymbol v_d^K)^T
\end{bmatrix}
=&{}-G_K E_K^{-1}
+E_K^{-1}\frac{\partial G_K}{\partial\boldsymbol J_K^{-T}}
 \widehat E_K E_K^{-1}
+\frac{\partial G_K}{\partial\delta_K}\,\delta_K E_K^{-1}\\
&-\frac{1}{d+1}\boldsymbol e_d
\sum_{j=0}^{d}\operatorname{tr}\!\left(
\frac{\partial G_K}{\partial\boldsymbol M_K}\boldsymbol M_{j,K}
\right)(\nabla\phi_{j,K})^T
-\frac{1}{d+1}\boldsymbol e_d
 \left(\frac{\partial G_K}{\partial\boldsymbol x}\right)^T,
\end{aligned}
\end{equation}
where \(\boldsymbol e_d=(1,\ldots,1)^T\in\mathbb R^d\), and
\begin{equation}\label{HuangDiscreteX0}
(\boldsymbol v_0^K)^T
=-\sum_{k=1}^{d}(\boldsymbol v_k^K)^T
-\sum_{j=0}^{d}\operatorname{tr}\!\left(
\frac{\partial G_K}{\partial\boldsymbol M_K}\boldsymbol M_{j,K}
\right)(\nabla\phi_{j,K})^T
-\left(\frac{\partial G_K}{\partial\boldsymbol x}\right)^T.
\end{equation}
The first three terms are the volume, inverse-Jacobian, and determinant variations; the remaining terms account for metric interpolation and explicit spatial dependence.  For the frozen epoch data used here, \(\partial G_K/\partial\boldsymbol x=\boldsymbol0\), while the metric-interpolation term remains.

These formulas are the physical-coordinate structure used in the standard semi-discrete nonsingularity analysis \cite{huang2016nonsingularity}.  The reported solver uses the corresponding $\boldsymbol\xi$-view residual.  We next reduce the two formulas for a generic density $T(\boldsymbol A,a)$.

\subsection{Discrete Simplification of \texorpdfstring{$\boldsymbol A$}{A} Structure}

Before describing the discrete reduction, we retain the two local matrix differentiation identities used in the original derivation.

\begin{lemma}\label{lem:AMAt}
For symmetric matrices $\boldsymbol G$ and $\boldsymbol M$, and an arbitrary matrix $\boldsymbol Z$,
\[
\mathrm{tr}\!\left(\boldsymbol G \frac{\partial (\boldsymbol Z\boldsymbol M\boldsymbol Z^T)}{\partial \boldsymbol Z^T}\right)
=2\boldsymbol M\boldsymbol Z^T\boldsymbol G.
\]
\end{lemma}

\begin{proof}
Taking differentials,
\[
\delta\,\operatorname{tr}(\boldsymbol G\boldsymbol Z\boldsymbol M\boldsymbol Z^T)
=\operatorname{tr}\!\left((2\boldsymbol M\boldsymbol Z^T\boldsymbol G)^T\delta\boldsymbol Z\right),
\]
where the symmetry of \(\boldsymbol G\) and \(\boldsymbol M\) has been used.  The displayed derivative follows.
\end{proof}

\begin{lemma}\label{lem:AMinvAt}
For symmetric matrices $\boldsymbol G$ and $\boldsymbol M$, and an arbitrary matrix $\boldsymbol Z$,
\[
\mathrm{tr}\!\left(\boldsymbol G\frac{\partial (\boldsymbol Z\boldsymbol M^{-1}\boldsymbol Z^T)}{\partial\boldsymbol M}\right)
=-\boldsymbol M^{-1}\boldsymbol Z^T\boldsymbol G\boldsymbol Z\boldsymbol M^{-1}.
\]
\end{lemma}

\begin{proof}
Since
\[
\delta(\boldsymbol Z\boldsymbol M^{-1}\boldsymbol Z^T)
=-\boldsymbol Z\boldsymbol M^{-1}(\delta\boldsymbol M)\boldsymbol M^{-1}\boldsymbol Z^T,
\]
cyclic invariance of the trace gives
\[
\delta\,\operatorname{tr}(\boldsymbol G\boldsymbol Z\boldsymbol M^{-1}\boldsymbol Z^T)
=-\operatorname{tr}\!\left((\boldsymbol M^{-1}\boldsymbol Z^T\boldsymbol G\boldsymbol Z\boldsymbol M^{-1})^T\delta\boldsymbol M\right).
\]
\end{proof}

The two lemmas are retained as the local algebraic tools used below.  The purpose of the following reduction is to identify two element kernels: \(\mathcal G_K^{\xi}\) for the computational-coordinate residual used by the solver, and \(\mathcal G_K^{x}\) for the physical-coordinate gradient flow used as the geometric reference.  Let
\[
\boldsymbol\Xi_K:=
\begin{bmatrix}
(\boldsymbol\xi_0^K)^T\\ \vdots\\ (\boldsymbol\xi_d^K)^T
\end{bmatrix},
\qquad
\boldsymbol X_K:=
\begin{bmatrix}
(\boldsymbol x_0^K)^T\\ \vdots\\ (\boldsymbol x_d^K)^T
\end{bmatrix}.
\]
We define the row-stacked local kernels by
\begin{equation}\label{PairedLocalGradients}
\mathcal G_K^{\xi}:=\frac{\partial G_K}{\partial\boldsymbol\Xi_K},
\qquad
\mathcal G_K^{x}:=\frac{1}{|K|}\frac{\partial(|K|G_K)}{\partial\boldsymbol X_K}.
\end{equation}
Every row is the transpose of the derivative with respect to the associated local vertex.  Thus \(\mathcal G_K^{\xi}\) is the local contribution to \(\nabla_{\boldsymbol\xi}I_h\), while \(\mathcal G_K^x\) is the normalized local contribution to \(\nabla_{\boldsymbol x}I_h\).  The local descent vectors in \cref{xview_equ} are the rows of \(-\mathcal G_K^x\).

Write
\[
\boldsymbol T_{A,K}:=\frac{\partial T_K}{\partial\boldsymbol A_K},
\qquad
T_{a,K}:=\frac{\partial T_K}{\partial a_K},
\]
where the matrix derivative is the symmetric-Frobenius derivative on \(\operatorname{Sym}(d)\).  The same convention is used for derivatives with respect to \(\boldsymbol M_K\).  For the passage from the $d$ non-base vertices to all $d+1$ vertices, set
\[
\boldsymbol R:=
\begin{bmatrix}
-1 & -1 & \cdots & -1\\
1 & 0 & \cdots & 0\\
0 & 1 & \cdots & 0\\
\vdots & \vdots & \ddots & \vdots\\
0 & 0 & \cdots & 1
\end{bmatrix}\in\mathbb R^{(d+1)\times d}.
\]

\begin{proposition}[Paired local gradient representations]
\label{prop:discrete-reduction}
Let \(T_K=T(\boldsymbol A_K,a_K)\) be differentiable.  In one
\(\boldsymbol\xi\)-view variation, keep \(E_K\), \(\boldsymbol M_K\), and the
epochwise scalar data fixed.  With the derivative data
\(\boldsymbol T_{A,K}=\partial T_K/\partial\boldsymbol A_K\) and
\(T_{a,K}=\partial T_K/\partial a_K\), define
\[
\boldsymbol Q_K:=\boldsymbol A_K\boldsymbol T_{A,K}
+a_KT_{a,K}\boldsymbol I \quad\text{and} \quad
\boldsymbol U_K
:=
\boldsymbol J_K
\left(-\frac12T_K\boldsymbol I+\boldsymbol Q_K\right)
\boldsymbol J_K^{-1}.
\]
Then the normalized row-stacked local \(\boldsymbol\xi\)-gradient is
\begin{equation}\label{Axiview}
\mathcal G_K^{\xi}
=\frac{\partial G_K}{\partial[\boldsymbol\xi_0^K,\boldsymbol\xi_1^K,\ldots,\boldsymbol\xi_d^K]}
=2\rho_K\boldsymbol R\widehat E_K^{-1}\boldsymbol Q_K.
\end{equation}

For the physical-coordinate representation, which differentiates the complete contribution \(|K|G_K\), let
\[
\boldsymbol V_K:=\boldsymbol R E_K^{-1}
=
\begin{bmatrix}
(\nabla\phi_{0,K})^T\\ \vdots\\ (\nabla\phi_{d,K})^T
\end{bmatrix},
\qquad
\boldsymbol e_{d+1}:=(1,\ldots,1)^T\in\mathbb R^{d+1}.
\]
Then
\begin{equation}\label{AXview}
\mathcal G_K^x
=-2\rho_K\boldsymbol V_K\boldsymbol U_K
+\frac{1}{d+1}\boldsymbol e_{d+1}
\begin{bmatrix}
\operatorname{tr}\!\left(\frac{\partial G_K}{\partial\boldsymbol M_K}\boldsymbol M_{0,K}\right)&\cdots&
\operatorname{tr}\!\left(\frac{\partial G_K}{\partial\boldsymbol M_K}\boldsymbol M_{d,K}\right)
\end{bmatrix}
\boldsymbol V_K.
\end{equation}
Hence the local descent vectors are the rows of \(-\mathcal G_K^x\).  The two kernels are assembled globally as
\[
\frac{\partial I_h}{\partial\boldsymbol\xi_i}
=\sum_{K\in S_i}|K|\,(\mathcal G_K^\xi)_{\ell(i,K),:},
\qquad
\frac{\partial I_h}{\partial\boldsymbol x_i}
=\sum_{K\in S_i}|K|\,(\mathcal G_K^x)_{\ell(i,K),:}.
\]
\end{proposition}

\begin{proof}
We first derive the \(\boldsymbol\xi\)-view formula row by row.  By
Lemma~\ref{lem:AMAt},
\[
\frac{\partial G_K}{\partial\boldsymbol J_K^{-T}}
=
2\rho_K\boldsymbol M_K^{-1}E_K^{-T}\widehat E_K^T
\boldsymbol T_{A,K}.
\]
Moreover,
\[
a_K=\det(\boldsymbol M_K^{-1})\delta_K^2,
\qquad
\frac{\partial G_K}{\partial\delta_K}
=
2\rho_K\det(\boldsymbol M_K^{-1})\delta_K\,T_{a,K}.
\]
For the \(d\) non-base vertices, \cref{HuangDiscreteXi} therefore gives
\[
\begin{aligned}
\frac{\partial G_K}
{\partial[\boldsymbol\xi_1^K,\ldots,\boldsymbol\xi_d^K]}
={}&
2\rho_K
E_K^{-1}\boldsymbol M_K^{-1}E_K^{-T}\widehat E_K^T
\boldsymbol T_{A,K}
+2\rho_K\det(\boldsymbol M_K^{-1})\delta_K^2
T_{a,K}\widehat E_K^{-1}\\
={}&
2\rho_K\widehat E_K^{-1}
\left(
\boldsymbol A_K\boldsymbol T_{A,K}
+a_KT_{a,K}\boldsymbol I
\right)
=
2\rho_K\widehat E_K^{-1}\boldsymbol Q_K.
\end{aligned}
\]
Here we used
\[
E_K^{-1}\boldsymbol M_K^{-1}E_K^{-T}\widehat E_K^T
=\widehat E_K^{-1}\boldsymbol A_K,
\qquad
\det(\boldsymbol M_K^{-1})\delta_K^2=a_K.
\]
The base-vertex identity \cref{HuangDiscreteXi0} then stacks these \(d\)
rows with their negative row sum, which is exactly
\(\mathcal G_K^\xi=2\rho_K\boldsymbol R\widehat E_K^{-1}\boldsymbol Q_K\).

For the \(\boldsymbol x\)-view, differentiating the complete contribution \(|K|G_K\), including its volume factor and the piecewise-affine metric
interpolation, yields
\[
\begin{aligned}
-\mathcal G_K^x
={}&\rho_K\boldsymbol R\Big(
-T_KE_K^{-1}
+2\widehat E_K^{-1}\boldsymbol A_K\boldsymbol T_{A,K}
\widehat E_KE_K^{-1}
+2a_KT_{a,K}E_K^{-1}
\Big)\\
&\quad-\frac{1}{d+1}\boldsymbol e_{d+1}
\begin{bmatrix}
\operatorname{tr}\!\left(\frac{\partial G_K}{\partial\boldsymbol M_K}\boldsymbol M_{0,K}\right)&\cdots&
\operatorname{tr}\!\left(\frac{\partial G_K}{\partial\boldsymbol M_K}\boldsymbol M_{d,K}\right)
\end{bmatrix}
\boldsymbol V_K.
\end{aligned}
\]
Using
\(\boldsymbol V_K=\boldsymbol R E_K^{-1}\),
\(\boldsymbol J_K=E_K\widehat E_K^{-1}\), and
\(\boldsymbol J_K^{-1}=\widehat E_KE_K^{-1}\), it transforms directly as
\[
\begin{aligned}
&\rho_K\boldsymbol R\Big(
-T_KE_K^{-1}
+2\widehat E_K^{-1}\boldsymbol A_K\boldsymbol T_{A,K}
\widehat E_KE_K^{-1}
+2a_KT_{a,K}E_K^{-1}
\Big)\\
&\qquad
=\rho_K\boldsymbol V_K
\left(
-T_K\boldsymbol I
+2\boldsymbol J_K\boldsymbol A_K\boldsymbol T_{A,K}\boldsymbol J_K^{-1}
+2a_KT_{a,K}\boldsymbol I
\right)\\
&\qquad
=2\rho_K\boldsymbol V_K\boldsymbol J_K
\left(
-\frac12T_K\boldsymbol I
+\boldsymbol A_K\boldsymbol T_{A,K}
+a_KT_{a,K}\boldsymbol I
\right)\boldsymbol J_K^{-1}
=2\rho_K\boldsymbol V_K\boldsymbol U_K.
\end{aligned}
\]
Substitution gives \cref{AXview}, and the relation with
\(\boldsymbol v_j^K\) follows from \cref{xview_equ}.
\end{proof}

The metric derivative in the interpolation term is taken on the symmetric matrix space.  With \(\boldsymbol P_K\) (and hence \(\delta_K\)) fixed, Lemma~\ref{lem:AMinvAt} gives the symmetric-Frobenius representative
\begin{equation}\label{GMKPartial}
D_{\boldsymbol M_K}G_K[\delta\boldsymbol M_K]
=
\left\langle
\operatorname{sym}\!\left(-\rho_K\boldsymbol U_K\boldsymbol M_K^{-1}\right),
\delta\boldsymbol M_K
\right\rangle_F,
\qquad
\delta\boldsymbol M_K\in\operatorname{Sym}(d).
\end{equation}
In \cref{AXview}, \(\partial G_K/\partial\boldsymbol M_K\) denotes this symmetric representative.  The compact form follows from
\(\partial\rho_K/\partial\boldsymbol M_K=\tfrac12\rho_K\boldsymbol M_K^{-1}\),
\(\partial a_K/\partial\boldsymbol M_K=-a_K\boldsymbol M_K^{-1}\), and
\(\boldsymbol J_K\boldsymbol A_K
=\boldsymbol M_K^{-1}\boldsymbol J_K^{-T}\).

\begin{corollary}[Trace--logarithmic specialization]
\label{cor:trace-log-discrete}
For the proposed density, the epochwise frozen scalar \(\theta_h\) gives
\begin{equation}\label{ConstructorPF}
\begin{aligned}
T_K&=[\operatorname{tr}(\boldsymbol A_K)]^p
-\kappa_{\theta_h}\log a_K,
&\kappa_{\theta_h}&=p d^{p-1}\theta_h^p,\\
\frac{\partial T_K}{\partial \boldsymbol A_K}
&=p[\operatorname{tr}(\boldsymbol A_K)]^{p-1}\boldsymbol I,
&\frac{\partial T_K}{\partial a_K}
&=-\kappa_{\theta_h}a_K^{-1}.
\end{aligned}
\end{equation}
Consequently,
\begin{equation}\label{TraceLogQ}
\boldsymbol Q_K
=p[\operatorname{tr}(\boldsymbol A_K)]^{p-1}\boldsymbol A_K
-\kappa_{\theta_h}\boldsymbol I,
\end{equation}
and
\begin{equation}\label{TraceLogU}
\boldsymbol U_K
=\boldsymbol J_K\left[
p[\operatorname{tr}(\boldsymbol A_K)]^{p-1}\boldsymbol A_K
-\left(\frac12T_K+\kappa_{\theta_h}\right)\boldsymbol I
\right]\boldsymbol J_K^{-1}.
\end{equation}
In particular, the implementation-oriented residual \cref{Axiview} becomes
\begin{equation}\label{TraceLogXi}
\mathcal G_K^{\xi}
=\frac{\partial G_K}{\partial [\boldsymbol\xi_0^K,\ldots,\boldsymbol\xi_d^K]}
=2\rho_K\boldsymbol R\widehat E_K^{-1}
\left[
p[\operatorname{tr}(\boldsymbol A_K)]^{p-1}\boldsymbol A_K
-\kappa_{\theta_h}\boldsymbol I
\right].
\end{equation}
\end{corollary}

\paragraph{Huang benchmark.}
The same local kernel is used for the Huang comparisons reported in
Section~5. With \(p=d\gamma/2\), the Huang density in \eqref{AHF}
satisfies
\[
T_{A,K}^{\rm H}
=
\mu p[\operatorname{tr}(\boldsymbol A_K)]^{p-1}\boldsymbol I,
\qquad
T_{a,K}^{\rm H}
=
(1-2\mu)d^p\frac{\gamma}{2}
a_K^{\gamma/2-1}.
\]
Consequently,
\[
\boldsymbol Q_K^{\rm H}
=
\mu p[\operatorname{tr}(\boldsymbol A_K)]^{p-1}\boldsymbol A_K
+
(1-2\mu)d^p\frac{\gamma}{2}
a_K^{\gamma/2}\boldsymbol I.
\]
Thus the proposed and Huang densities use the same geometric assembly
and nonlinear solver; only the density-specific element data differ.

\subsection{Moving Mesh Algorithms}

The reported implementation advances the computational-coordinate gradient
flow \cref{eq:xi-view-flow}.  After imposing the prescribed boundary
correspondence, let \(\boldsymbol\xi_f\) collect the independent coordinates.
At outer epoch \(n\), the physical mesh, metric data, and nodal factors
\(\omega_i^{\xi,(n)}\) are frozen; \(\widehat E_K\), \(\boldsymbol A_K\), and
\(\boldsymbol Q_K\) are recomputed at each residual evaluation.  Thus the
balancing factors are fixed positive row scalings during the inner
artificial-time solve.

For a \(k\)-step backward differentiation formula (BDF) stage \(m\to m+1\) \cite{hairer1996solving}, the residual is
\begin{equation}\label{eq:xi-gradient-flow}
\mathcal R^{(n,m+1)}(\boldsymbol\xi_f)
=
\alpha_0^{(m)}\boldsymbol\xi_f-\boldsymbol b_m+
h_m\mathcal V_h^{(n)}(\boldsymbol\xi_f),
\quad
\boldsymbol b_m=-\sum_{j=1}^{k}\alpha_j^{(m)}
\boldsymbol\xi_f^{m+1-j},
\quad
h_m=\frac{\Delta t_m}{\tau},
\end{equation}
where \(\mathcal V_h^{(n)}\) is assembled from \cref{TraceLogXi} with the
frozen factors \(\omega_i^{\xi,(n)}\).  At the start of a compact
Direct-Secant cycle, we assemble and factorize the base tangent
\begin{equation}\label{eq:base-tangent}
\boldsymbol{\mathcal B}_0
=D\mathcal R^{(n,m+1)}(\boldsymbol\xi_f^{(0)})
=\alpha_0^{(m)}\boldsymbol I+
h_mD\mathcal V_h^{(n)}(\boldsymbol\xi_f^{(0)}).
\end{equation}
The elementwise differentiation used to assemble
\(D\mathcal V_h^{(n)}\) is recorded in Appendix~\ref{app:base-tangent}.

Within a compact cycle, the tangent is represented as
\[
\boldsymbol{\mathcal B}_\ell
=\boldsymbol{\mathcal B}_0+\mathcal C_\ell\mathcal D_\ell^T.
\]
For \(\boldsymbol r_\ell=\mathcal R^{(n,m+1)}(\boldsymbol\xi_f^{(\ell)})\),
the search direction is recovered from the base factorization by
\begin{equation}\label{eq:smw-direction}
\boldsymbol p_\ell
=-\boldsymbol q_\ell+
\mathcal W_\ell\mathcal S_\ell^{-1}
\mathcal D_\ell^T\boldsymbol q_\ell,
\quad
\boldsymbol q_\ell=\boldsymbol{\mathcal B}_0^{-1}\boldsymbol r_\ell,
\quad
\mathcal W_\ell=\boldsymbol{\mathcal B}_0^{-1}\mathcal C_\ell,
\quad
\mathcal S_\ell=\boldsymbol I+\mathcal D_\ell^T\mathcal W_\ell.
\end{equation}
Accepted steps are required to preserve the boundary correspondence, keep \(\det\widehat E_K>0\), and decrease the residual.  Compact records are restarted from the current accepted iterate when the memory cap is reached or \(\mathcal S_\ell\) fails the conditioning check.  The low-rank update below is a compact secant correction \cite{broyden1965class,fletcher1963rapidly,nocedal2006numerical}.

For an accepted step, set
\[
\boldsymbol s_\ell=\boldsymbol\xi_f^{(\ell+1)}-\boldsymbol\xi_f^{(\ell)},
\qquad
\boldsymbol t_\ell=\boldsymbol r_{\ell+1}-\boldsymbol r_\ell,
\qquad
\boldsymbol w_\ell=\boldsymbol{\mathcal B}_\ell\boldsymbol s_\ell.
\]
Let \(\varepsilon_{\rm sec}>0\), \(\kappa_{\max}>1\), and \(\ell_{\max}\) be fixed solver controls.  
We use \(\varepsilon_{\rm sec}=10^{-12}\) and \(\ell_{\max}=12\); \(\kappa_{\max}\) is only a large protective threshold.  A compact cycle is restarted when
\begin{equation}\label{eq:cycle-admissibility}
\operatorname{cond}_2(\mathcal S_\ell)>\kappa_{\max}
\qquad\text{or}\qquad
\ell\ge\ell_{\max}.
\end{equation}
Within an admissible cycle, the rank-two update is used when
\begin{equation}\label{eq:secant-admissibility}
|\boldsymbol s_\ell^T\boldsymbol t_\ell|
\ge\varepsilon_{\rm sec}\|\boldsymbol s_\ell\|_2\|\boldsymbol t_\ell\|_2,\quad
|\boldsymbol s_\ell^T\boldsymbol w_\ell|
\ge\varepsilon_{\rm sec}\|\boldsymbol s_\ell\|_2\|\boldsymbol w_\ell\|_2.
\end{equation}
Otherwise, the rank-one fallback is used.  The update is
\begin{equation}\label{eq:direct-secant-update}
\boldsymbol{\mathcal B}_{\ell+1}
=
\boldsymbol{\mathcal B}_\ell+
\begin{cases}
\displaystyle
\frac{\boldsymbol t_\ell\boldsymbol t_\ell^T}
{\boldsymbol s_\ell^T\boldsymbol t_\ell}
-\frac{\boldsymbol w_\ell\boldsymbol w_\ell^T}
{\boldsymbol s_\ell^T\boldsymbol w_\ell},
&\text{rank-two},\\[1.1em]
\displaystyle
\frac{(\boldsymbol t_\ell-\boldsymbol w_\ell)\boldsymbol s_\ell^T}
{\boldsymbol s_\ell^T\boldsymbol s_\ell},
&\text{rank-one fallback}.
\end{cases}
\end{equation}
Both updates satisfy
\(\boldsymbol{\mathcal B}_{\ell+1}\boldsymbol s_\ell=\boldsymbol t_\ell\).
The compact representation is closed by appending
\[
(\mathcal C_{\ell+1},\mathcal D_{\ell+1})=
\begin{cases}
\left([
\mathcal C_\ell,\boldsymbol t_\ell,-\boldsymbol w_\ell],
[
\mathcal D_\ell,
\boldsymbol t_\ell/(\boldsymbol s_\ell^T\boldsymbol t_\ell),
\boldsymbol w_\ell/(\boldsymbol s_\ell^T\boldsymbol w_\ell)]\right),
&\text{rank-two},\\[.4em]
\left([
\mathcal C_\ell,\boldsymbol t_\ell-\boldsymbol w_\ell],
[
\mathcal D_\ell,
\boldsymbol s_\ell/(\boldsymbol s_\ell^T\boldsymbol s_\ell)]\right),
&\text{rank-one fallback}.
\end{cases}
\]
Hence \(\mathcal C_{\ell+1}\mathcal D_{\ell+1}^T\) reproduces
\cref{eq:direct-secant-update} exactly; only the nonlinear stage model is updated.

If no feasible residual-decreasing step is found after the permitted restart,
the BDF stage is rejected, the preceding accepted state is restored, and
\(\Delta t_m\) is reduced.  A stage is accepted when
\begin{equation}\label{eq:secant-stop}
\|\mathcal R^{(n,m+1)}(\boldsymbol\xi_f^{(\ell)})\|_2
\le
a_{\rm tol}+r_{\rm tol}
\|\mathcal R^{(n,m+1)}(\boldsymbol\xi_f^{(0)})\|_2.
\end{equation}

\begin{algorithm}[!htbp]\small
\caption{\(\boldsymbol\xi\)-view BDF Direct-Secant Update Algorithm}
\label{alg:xi-update}
\begin{algorithmic}[1]
\Require Initial \(\boldsymbol\xi^0\), outer epoch count, BDF controls, \(\varepsilon_{\rm sec}\), \(\kappa_{\max}\), \(\ell_{\max}\), and boundary correspondence
\Ensure Updated physical mesh \(\mathcal T_h\)
\For{each outer epoch \(n\)}
    \State Construct and freeze the metric data, \(\theta_h^{(n)}\), \(\rho_K\), and \(\omega_i^{\xi,(n)}\)
    \For{each BDF stage \(m\to m+1\)}
        \State Construct a feasible predictor and assemble \(\boldsymbol{\mathcal B}_0\)
        \While{the criterion \cref{eq:secant-stop} is not satisfied}
            \State For residual evaluations, compute \(\boldsymbol A_K\), \(\boldsymbol Q_K\), and \(\mathcal G_K^\xi\) from \cref{TraceLogQ,TraceLogXi} and assemble \(\mathcal V_h^{(n)}\)
            \State Recover \(\boldsymbol p_\ell\) from \cref{eq:smw-direction}; restart the compact cycle when required
            \State Perform the feasible residual-decreasing line search
            \State On acceptance, append the direct-secant update by \cref{eq:direct-secant-update}
        \EndWhile
        \State Accept the stage; on rejection restore the preceding state and reduce \(\Delta t_m\)
    \EndFor
    \State Recover the physical mesh and transfer the physical quantities by the standard piecewise-linear inverse-map operation \cite{huang2015geometric}
\EndFor
\State \Return \(\mathcal T_h\)
\end{algorithmic}
\end{algorithm}

\section{Numerical Experiments}

\subsection{Experimental Setup}

The computational and physical domains coincide \((\Omega_c=\Omega_p)\) and use simplicial meshes.  
Within each functional--metric run, the initial mesh, boundary treatment, metric protocol, BDF controller, and nonlinear stopping rule are fixed across resolutions; the tolerance in \cref{eq:secant-stop} is \(10^{-6}\) unless stated otherwise.  
With the stated functional-specific artificial-time scales, wall-clock times and accepted nonlinear updates provide implementation diagnostics.

Following \cite{kolasinski2018new}, Huang's functional \cref{HF} is the benchmark.

Metrics are built at the start of each outer epoch, smoothed as stated below, and frozen during the inner artificial-time solve.
\begin{itemize}
\item \emph{Hessian-based metric} \cite{huang2005metric}:
\begin{equation}\label{HessianMetric}
  \boldsymbol B_K=\boldsymbol I+\alpha_h^{-1}|\boldsymbol H_K(u)|,
  \qquad
  \boldsymbol M_K=\det(\boldsymbol B_K)^{-\frac{1}{d+4}}\boldsymbol B_K,
\end{equation}
where \(\boldsymbol H_K(u)\) is the cellwise recovered Hessian and
\(|\boldsymbol H_K|\) is its spectral absolute value.  In the two-dimensional
experiments, \(\alpha_h\) is chosen so that
\begin{equation}\label{eq:alpha-h}
\sum_{K\in\mathcal T_h}|K|\sqrt{\det(\boldsymbol M_K)}=2|\Omega_p|,
\end{equation}
which fixes the metric-volume scale without introducing a further mesh-control
parameter.

\item \emph{Arc-length-like metric} \cite{shen2009efficient,li2001moving}:
\begin{equation}\label{ArcLengthMetric}
\boldsymbol M_K = \sqrt{1 + \beta |\nabla u_h|_K^2}\,\boldsymbol I,
\end{equation}
with \(\beta=0.5\) in the Burgers test.

\item \emph{Eigen-decomposition metric} \cite{huang2001practical}:
\begin{equation}\label{EigenDecompositionMetric}
\boldsymbol M_K = \lambda_{n,K} \boldsymbol v_K \boldsymbol v_K^{T}
+ \lambda_{t,K} \boldsymbol v_K^{\perp}(\boldsymbol v_K^{\perp})^{T},
\end{equation}
where \(\boldsymbol v_K\) is the phase-field normal in the Rayleigh--Taylor test.
\end{itemize}
Examples~\ref{exp:1}--\ref{exp:2} use \cref{HessianMetric}; the woven-junction test uses \cref{eq:woven-raw-metric}; Burgers uses \cref{HessianMetric,ArcLengthMetric}; and Rayleigh--Taylor uses \cref{HessianMetric,EigenDecompositionMetric}.  On the $80\times80$ metric grids used for the mesh-quality visualizations and stress tests, eight nodal-averaging sweeps are applied before mesh adaptation.

Unless otherwise specified, the proposed functional uses $\gamma=1.25$ and
$\tau=0.004$, while Huang's functional employs $\gamma=3/2$ and $\tau=0.01$,
following standard choices in the literature.  In the matched finite-budget
woven-junction stress test of \cref{sec:woven-junction}, each density retains
its functional-specific exponent---$\gamma=1.25$ for the trace--logarithmic
density and $\gamma=3/2$ for Huang's functional---while both use the common
artificial-time scale $\tau=0.02$.

All experiments use FEALPy \cite{zheng2025fealpy} and ParaView and were run on a 12th Gen Intel(R) Core(TM) i7-12700 CPU with 16.0~GB of DDR4 memory at 3200~MT/s.  Every numerical figure and table identifies its metric in the caption, legend, or column label.

\subsection{Mesh Distribution Induced by Two-Dimensional Functions}
\label{with_functions}

To assess mesh quality, we employ three standard metrics
\cite{huang2015geometric,kolasinski2018new}:
the equidistribution measure $Q_{eq}$,
the alignment measure $Q_{ali}$,
and the geometric quality measure $Q_{geo}$.
The global indicator is defined as
\[
Q = \sqrt{\frac{1}{NC}\sum_{K\in \mathcal{T}_c} Q_K^2},
\]
where $NC$ is the number of elements.
The element-wise measures are
\[
Q_{eq,K} =
\frac{|K| \det(\boldsymbol M_K)^{\frac{1}{2}}}{\sigma_h/NC}, 
Q_{ali,K} =
\frac{\mathrm{tr}(\boldsymbol J_K^{T} \boldsymbol M_K \boldsymbol J_K)}
{d\,\det(\boldsymbol J_K^{T} \boldsymbol M_K \boldsymbol J_K)^{\frac{1}{d}}},
Q_{geo,K} =
\frac{\mathrm{tr}(\boldsymbol J_K^{T} \boldsymbol J_K)}
{d\,\det(\boldsymbol J_K^{T} \boldsymbol J_K)^{\frac{1}{d}}}.
\]
The ideal values are one.  The local distributions and the tables use the standard indicators directly; in particular, $Q_{ali,K}\ge1$ and smaller values mean better local alignment.

\begin{example}\label{exp:1}
A sinusoidally modulated steep-gradient function:
\[
f(x,y) = \tanh\!\bigl(-100(y-0.5-0.25\sin(2\pi x))\bigr),
\quad (x,y)\in[0,1]^2.
\]
\end{example}

\begin{example}\label{exp:2}
An X-shaped anisotropic test function:
\[
f(x,y) = \tanh\!\bigl(100(1-x-y)\bigr)
- \tanh\!\bigl(100(x-y)\bigr),
\quad (x,y)\in[0,1]^2.
\]
\end{example}

Both mesh flows are advanced to the artificial time $t_{\rm mesh}=1$.
For Example~\ref{exp:1}, $\Delta t=0.05$,
while for Example~\ref{exp:2}, $\Delta t=0.1$.
We compare the proposed functional with Huang's functional.

\begin{figure}[!htbp]
\centering
\subfloat[Exp~\ref{exp:1} Trace-log]{
\includegraphics[width=0.31\textwidth]{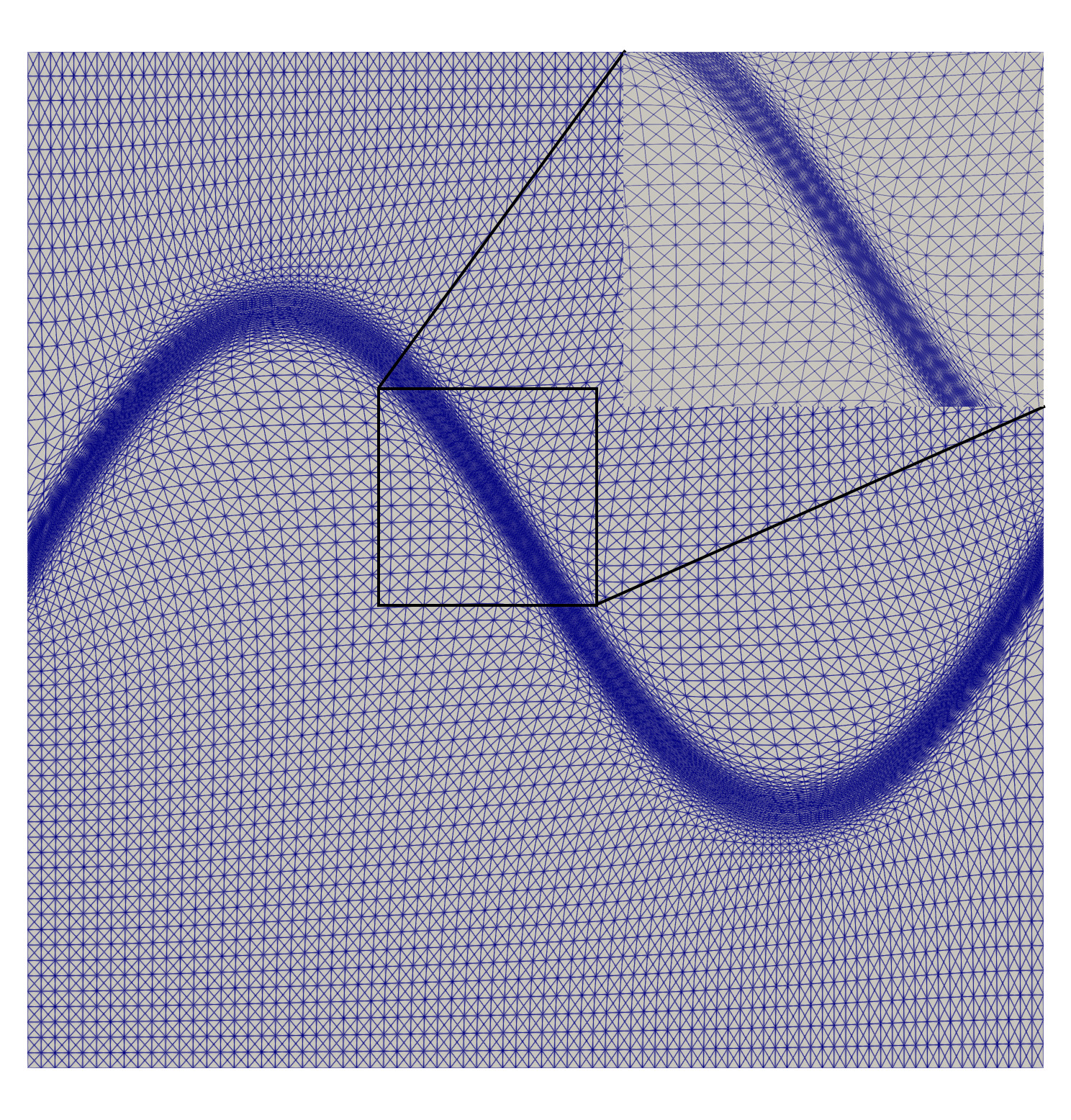}}
\subfloat[Exp~\ref{exp:1} Huang]{
\includegraphics[width=0.31\textwidth]{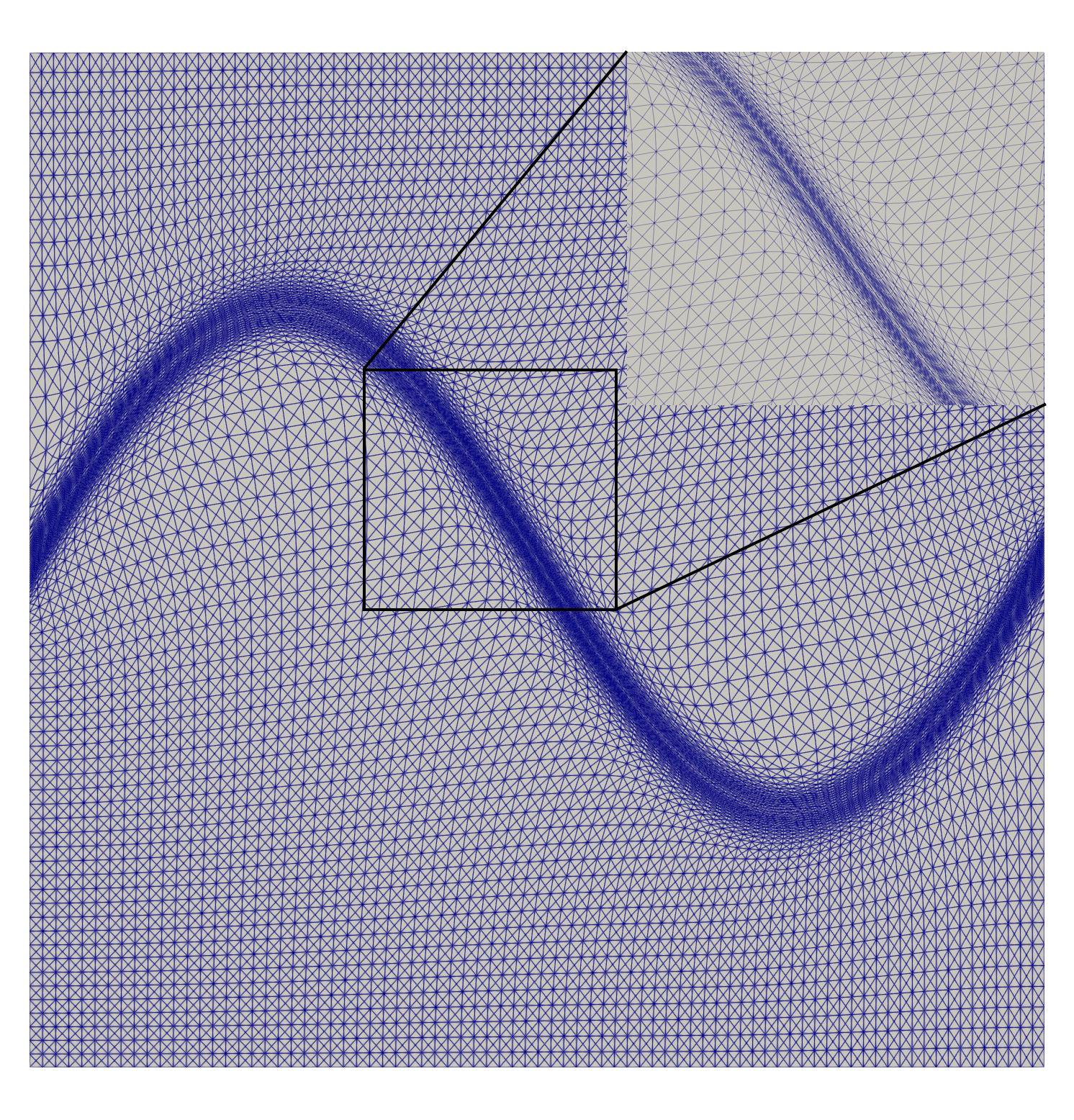}}\\[-0.3ex]
\subfloat[Exp~\ref{exp:2} Trace-log]{
\includegraphics[width=0.31\textwidth]{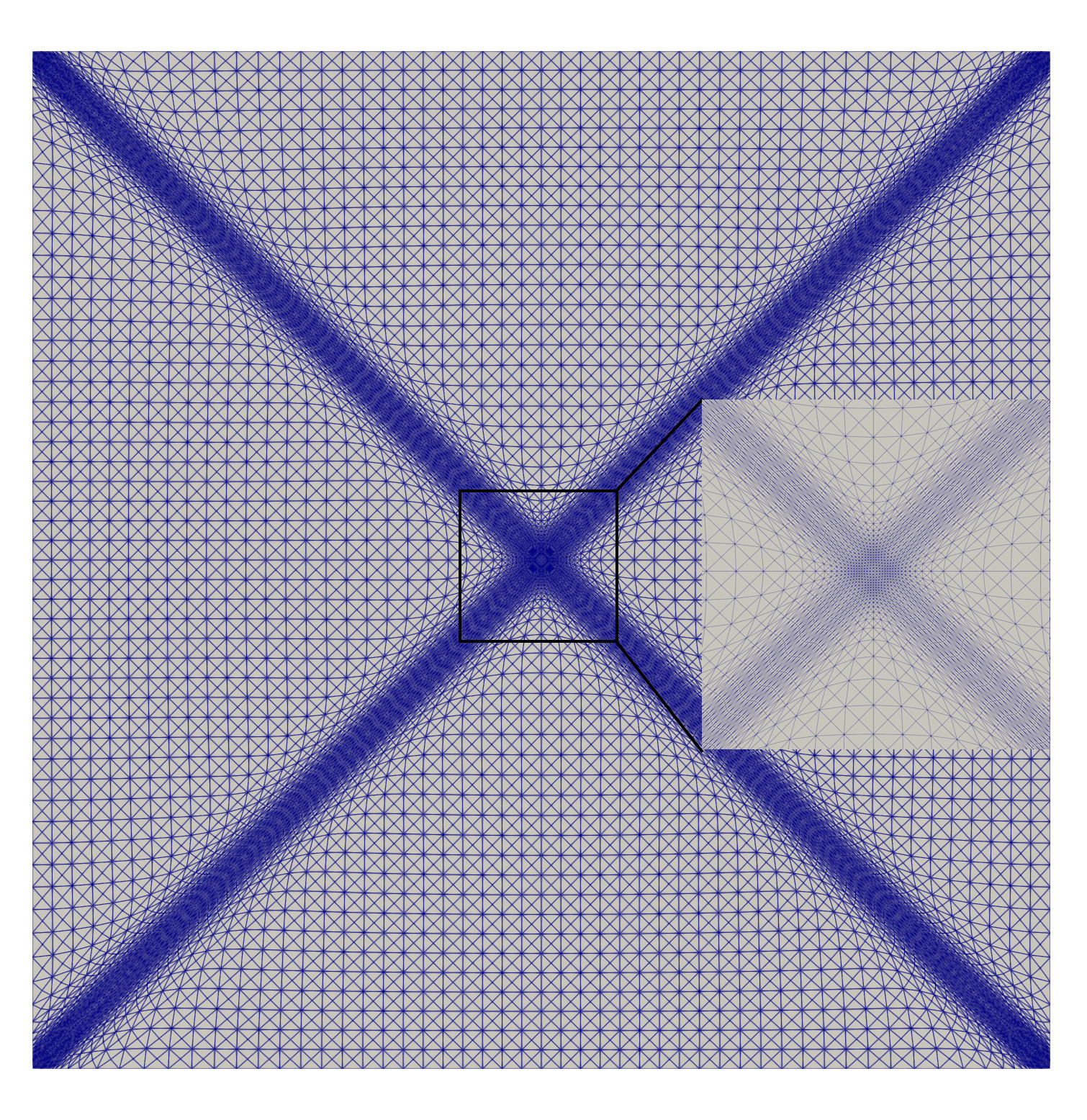}}
\subfloat[Exp~\ref{exp:2} Huang]{
\includegraphics[width=0.31\textwidth]{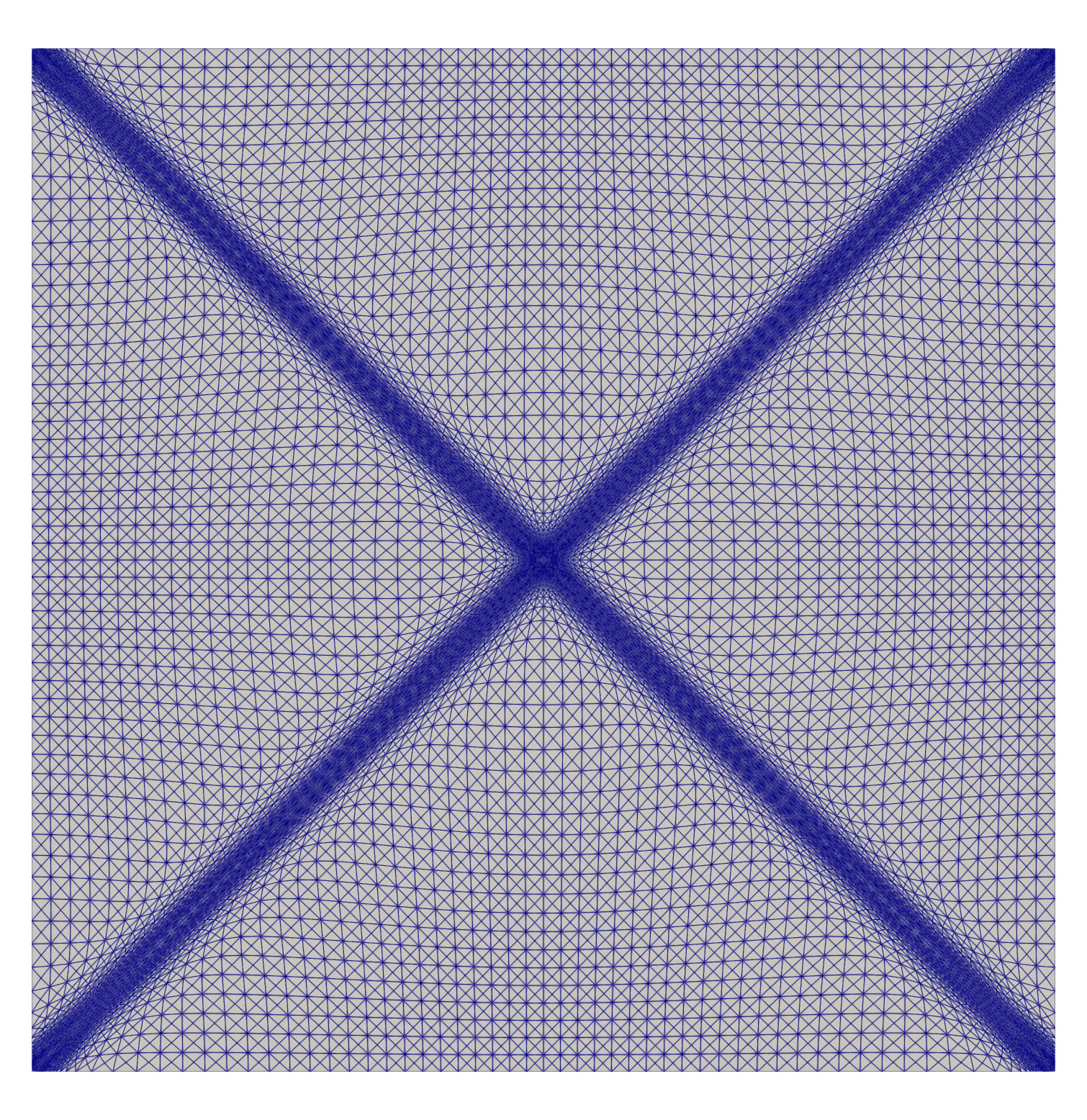}}
\caption{Mesh distributions for Examples~\ref{exp:1} and~\ref{exp:2}, both using the Hessian-based metric \cref{HessianMetric}, with 25600 elements.}
\label{Fig1.main}
\end{figure}

The resulting mesh distributions are shown in \cref{Fig1.main}.
For both examples, the meshes align well with the prescribed metric.
In Example~\ref{exp:1}, strong refinement is observed along the sinusoidally
modulated steep-gradient band, while elements remain relatively coarse
along directions of small curvature.
In Example~\ref{exp:2}, the mesh aligns naturally with the two intersecting
anisotropic lines and transitions smoothly toward nearly isotropic elements
near the intersection, without visible distortion.

\begin{figure}[!htbp]
\centering
\subfloat[Exp~\ref{exp:1} Trace-log and Huang]{
\includegraphics[width=0.88\textwidth]{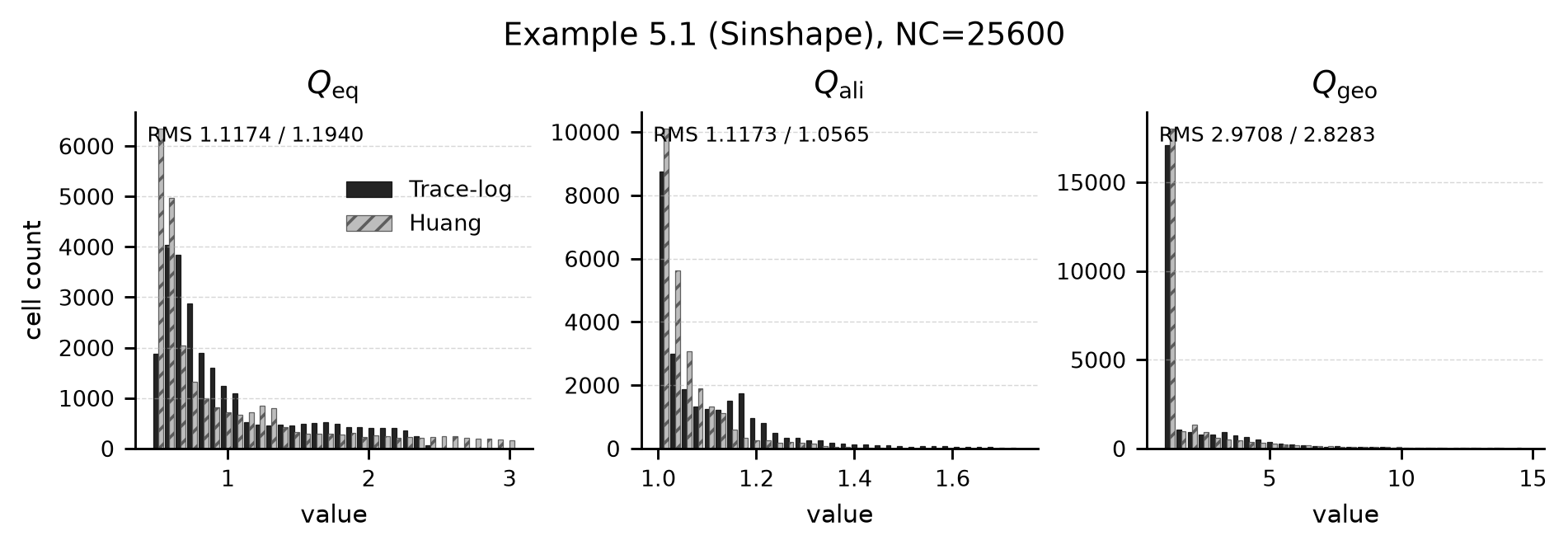}}\\[-0.4ex]
\subfloat[Exp~\ref{exp:2} Trace-log and Huang]{
\includegraphics[width=0.88\textwidth]{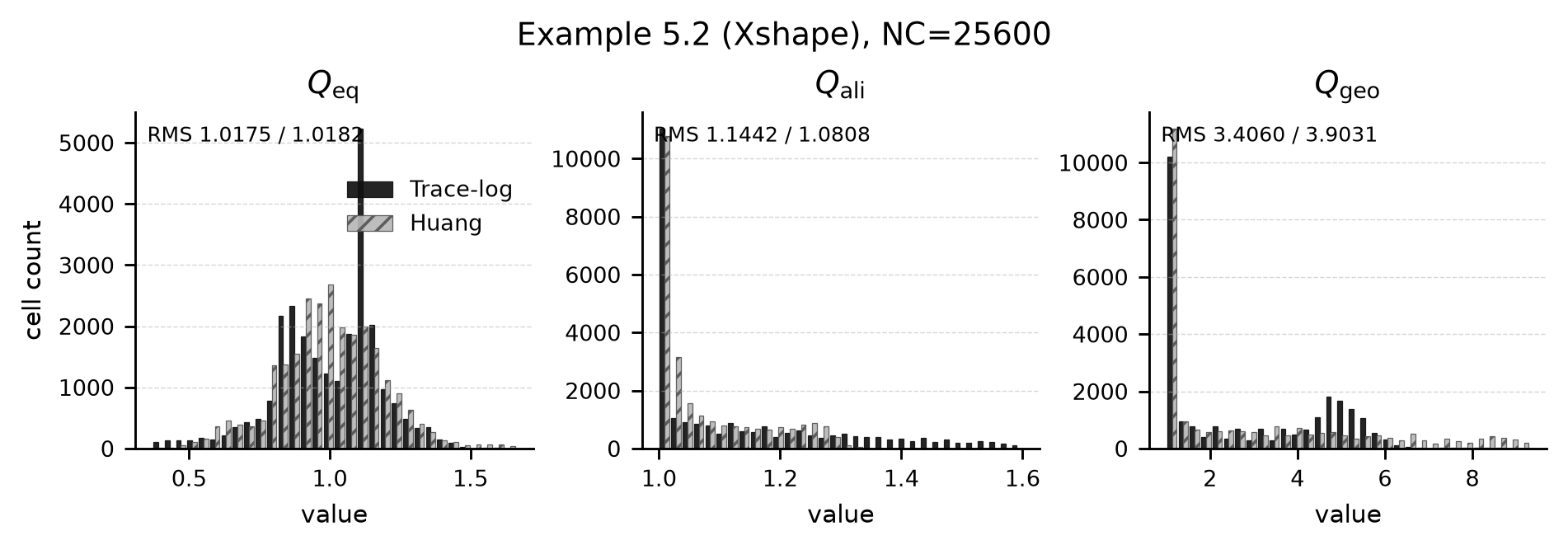}}
\caption{Element-wise mesh-quality distributions for Examples~\ref{exp:1} and~\ref{exp:2}, both using the Hessian-based metric \cref{HessianMetric}.  The distributions use the standard local indicators $Q_{eq,K}$, $Q_{ali,K}$, and $Q_{geo,K}$; values closer to one indicate better local quality.}
\label{Fig2.main}
\end{figure}

The distributions in \cref{Fig2.main} show tighter clustering of $Q_{eq}$ around one for the proposed functional.  Huang's functional is slightly better aligned in some fine-mesh cases, but the alignment and geometric-quality levels remain comparable.  The global data in \cref{tab:metrics-sinshape,tab:metrics-xshape} show consistently smaller $Q_{eq}$ for the proposed functional and $L_2$ errors of the same order; neither method has a uniform interpolation-error advantage at every resolution.

\begin{figure}[!htbp]
\centering
\subfloat[Exp.~\ref{exp:1}: $I_h$]{
\includegraphics[width=0.4\textwidth]{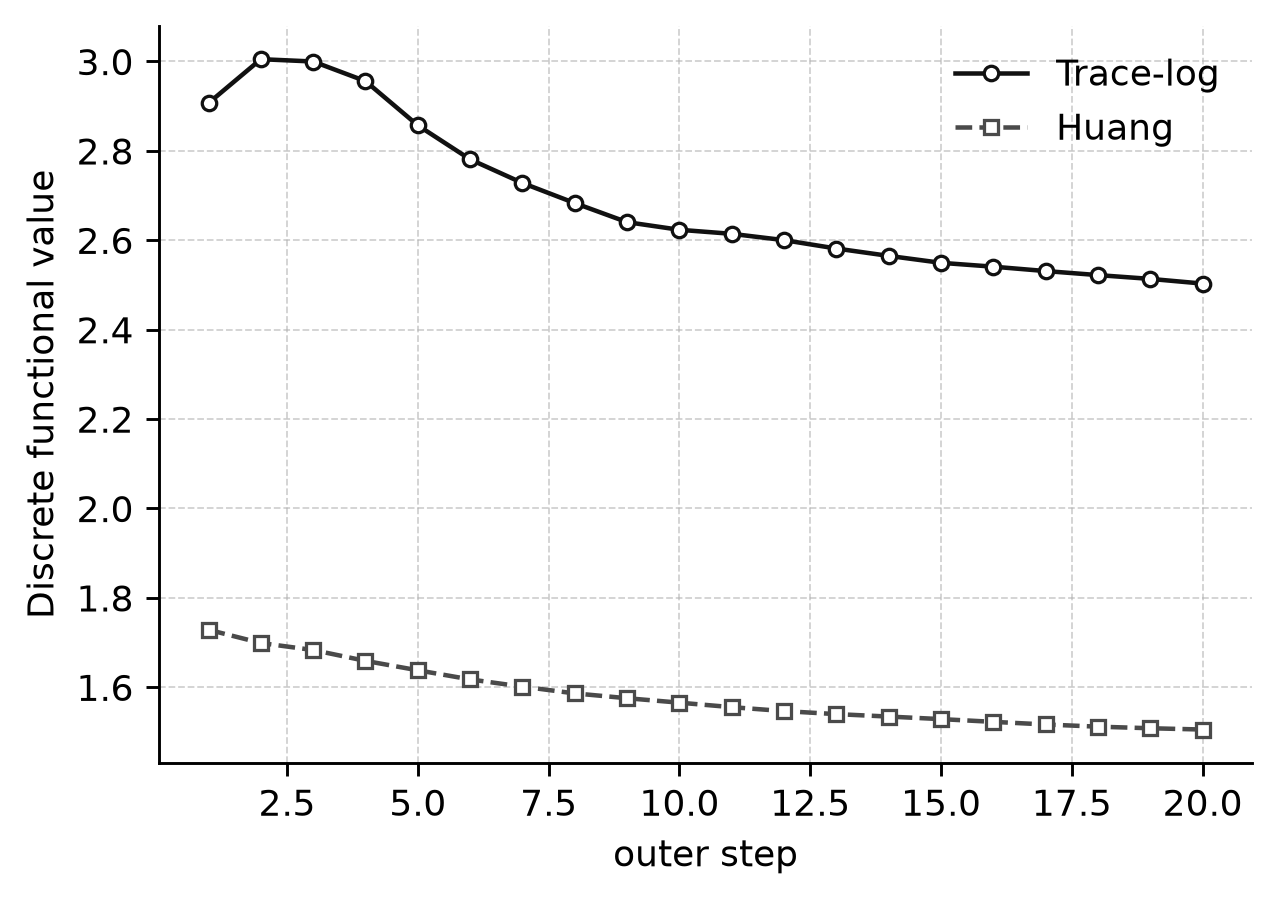}}
\subfloat[Exp.~\ref{exp:1}: $\min |K|$]{
\includegraphics[width=0.4\textwidth]{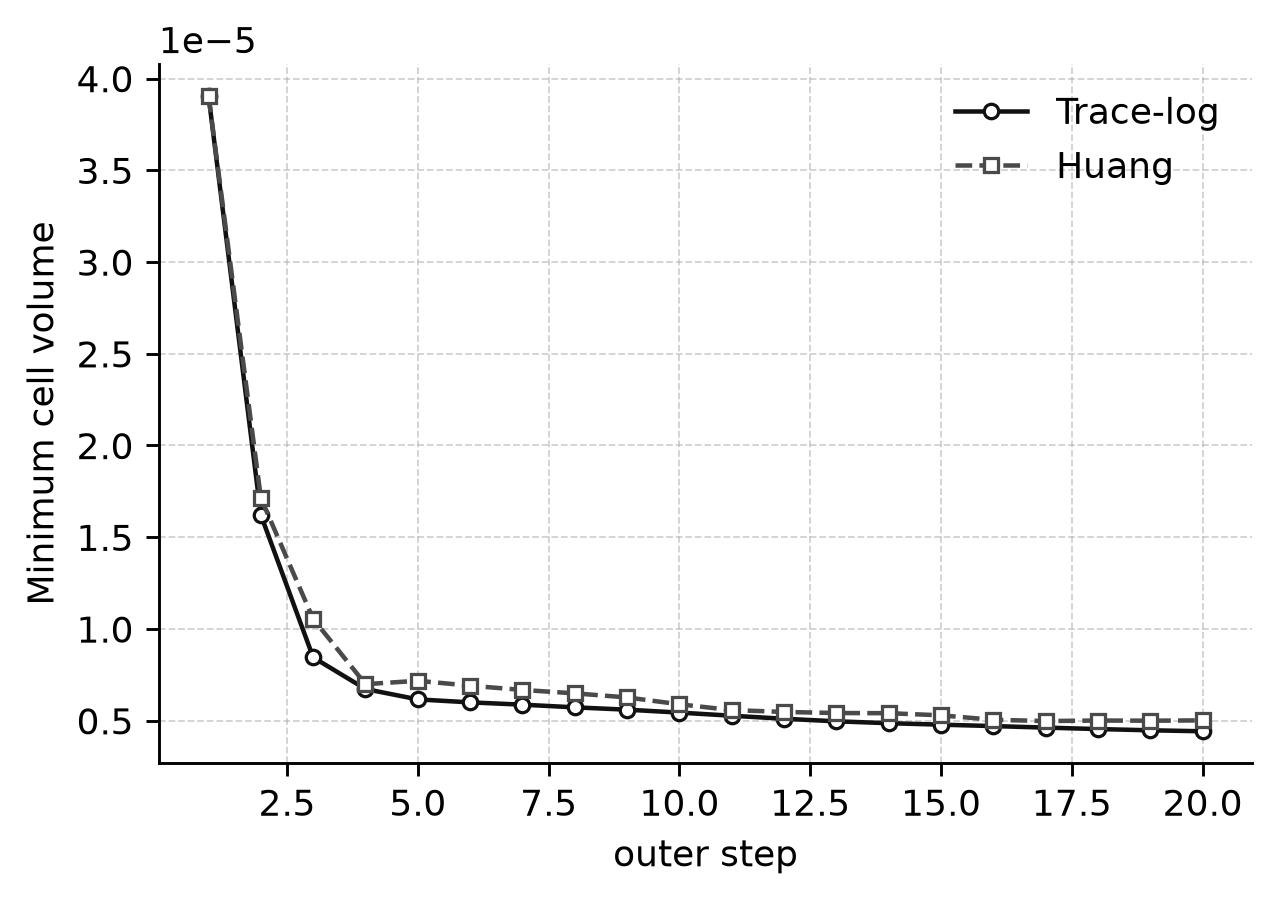}}\\[-0.2ex]
\subfloat[Exp.~\ref{exp:2}: $I_h$]{
\includegraphics[width=0.4\textwidth]{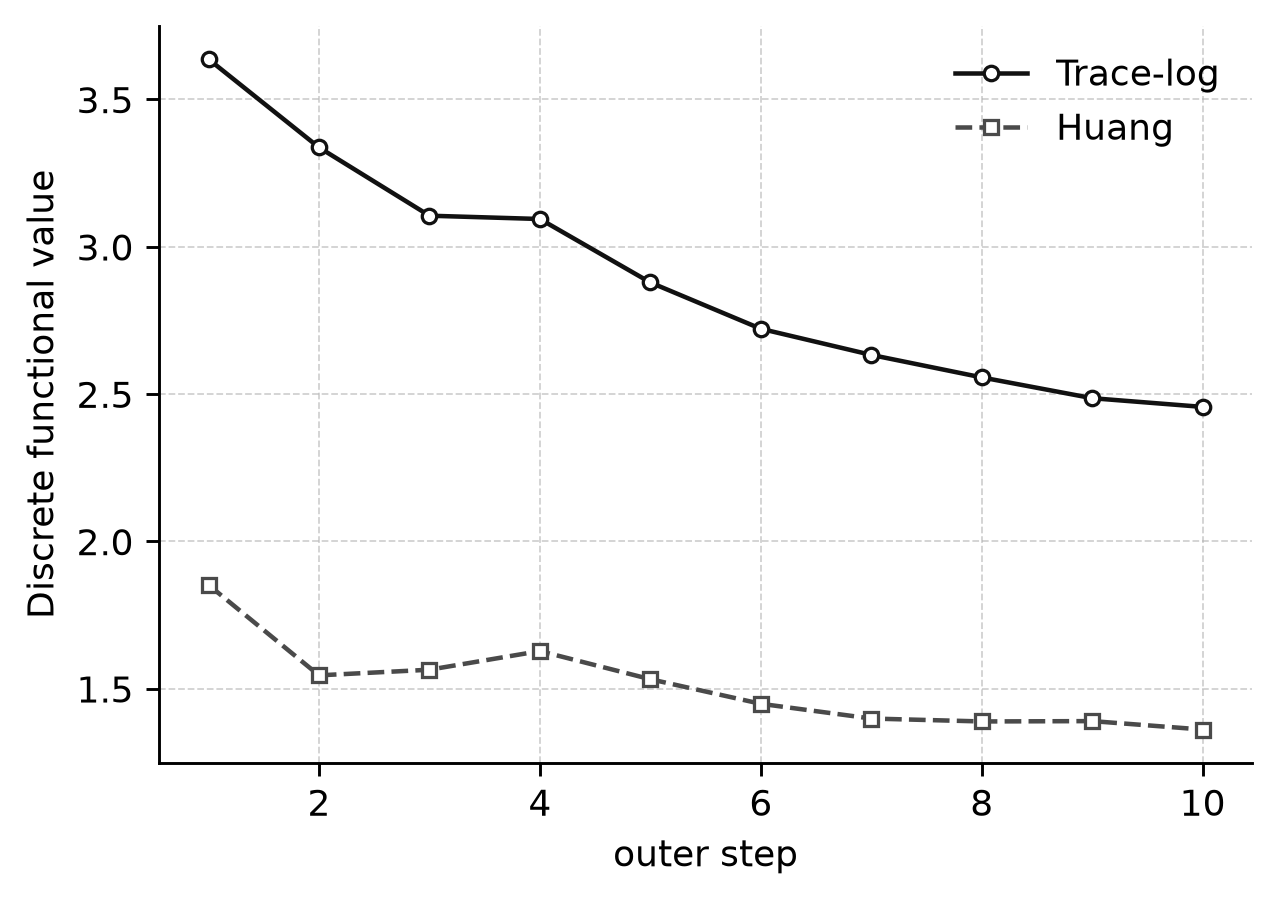}}
\subfloat[Exp.~\ref{exp:2}: $\min |K|$]{
\includegraphics[width=0.4\textwidth]{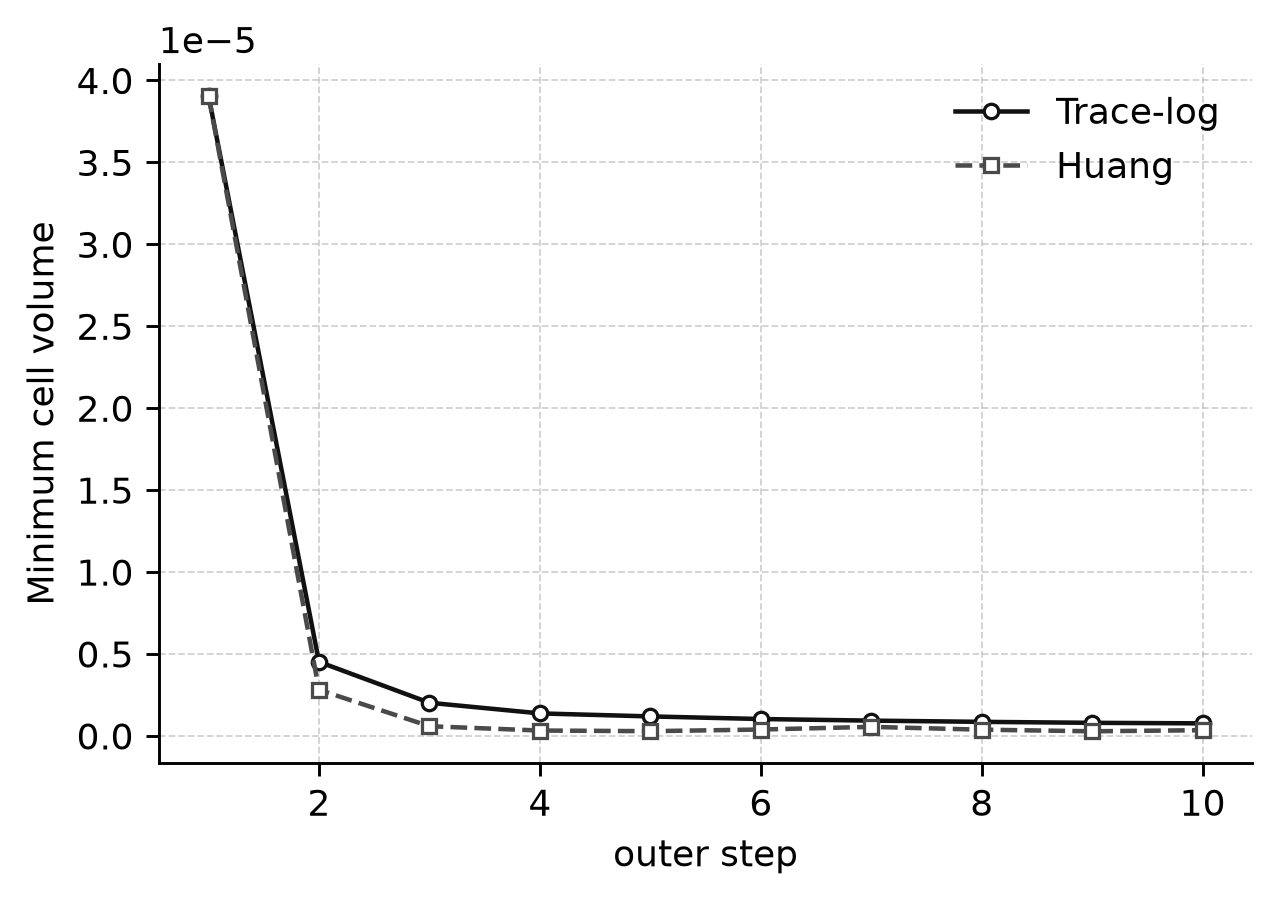}}
\caption{Evolution of the epochwise discrete functional and minimum element volume for Examples~\ref{exp:1} and~\ref{exp:2} under the Hessian-based metric \cref{HessianMetric}.}
\label{Fig3.main}
\end{figure}

For the accepted BDF iterates in \cref{Fig3.main}, the epochwise energy decreases and the minimum element volume stays positive.  These are numerical diagnostics; the result in Section~3 concerns the semi-discrete physical-coordinate MMPDE.

\begin{table}[!htbp]
\centering
\setlength{\tabcolsep}{3pt}
\renewcommand{\arraystretch}{0.9}
\caption{Example~\ref{exp:1} under the Hessian-based metric \cref{HessianMetric}: global mesh-quality measures, interpolation errors, wall-clock times, and accepted nonlinear updates.}
\label{tab:metrics-sinshape}
\begin{tabular}{c c c c c c c c}
\hline
functional & size(NC) & $Q_{eq}$ & $Q_{ali}$ & $Q_{geo}$ & $e_{L_2}$ & time(s) & nonlinear iters \\
\hline
\multirow{3}{*}{Huang} & 1600 & 1.09460 & 1.08452 & 2.27007 & 0.01384 & 0.379 & 136 \\
 & 6400 & 1.12344 & 1.07156 & 3.33749 & 0.00212 & 2.033 & 204 \\
 & 25600 & 1.19395 & 1.05653 & 2.82831 & 0.00055 & 9.058 & 229 \\
\hline
\multirow{3}{*}{Trace-log} & 1600 & 1.07685 & 1.11098 & 2.51211 & 0.01223 & 0.533 & 213\\
 & 6400 & 1.08087 & 1.11575 & 3.56923 & 0.00190 & 2.001 & 236\\
 & 25600 & 1.11740 & 1.11727 & 2.97077 & 0.00048 & 9.138 & 253\\
\hline
\end{tabular}
\end{table}

\begin{table}[!htbp]
\centering
\setlength{\tabcolsep}{3pt}
\renewcommand{\arraystretch}{0.9}
\caption{Example~\ref{exp:2} under the Hessian-based metric \cref{HessianMetric}: global mesh-quality measures, interpolation errors, wall-clock times, and accepted nonlinear updates.}
\label{tab:metrics-xshape}
\begin{tabular}{c c c c c c c c}
\hline
functional & size(NC) & $Q_{eq}$ & $Q_{ali}$ & $Q_{geo}$ & $e_{L_2}$ & time(s) & nonlinear iters \\
\hline
\multirow{3}{*}{Huang} & 1600 & 1.03923 & 1.06254 & 2.64811 & 0.01522 & 0.283 & 106 \\
 & 6400 & 1.06361 & 1.09446 & 5.00015 & 0.00125 & 1.682 & 183 \\
 & 25600 & 1.01816 & 1.08076 & 3.90312 & 0.00044 & 10.398 & 228 \\
\hline
\multirow{3}{*}{Trace-log} & 1600 & 1.02897 & 1.06436 & 2.44556 & 0.01807 & 0.347 & 157\\
 & 6400 & 1.04535 & 1.09926 & 4.06382 & 0.00150 & 1.705 & 200\\
 & 25600 & 1.01753 & 1.14424 & 3.40602 & 0.00043 & 8.186 & 216\\
\hline
\end{tabular}
\end{table}


\subsection{High-Anisotropy Woven-Junction Stress Test}
\label{sec:woven-junction}

The preceding examples use solution-based metrics.  We next apply the same finite-budget artificial-time protocol to a densely coupled woven metric and increase the prescribed compression level.  This probes robustness under localized, competing anisotropy within a fixed protocol.

Let \(n_J\in\{2,3,4\}\), with junctions at
\[
 x_i=y_i=\frac{i}{n_J+1},\qquad i=1,\ldots,n_J,
\]
so that \(n_J^2=4,9,16\) intersections are present.  For a compression parameter
\(\Lambda>1\), define
\[
 w_\Lambda=\frac{0.8}{\Lambda},\qquad
 \chi_\Lambda(s)=\exp\!\left[-\left(\frac{s}{w_\Lambda}\right)^4\right],
\]
\[
 \Phi_h(x,y)=\sum_{j=1}^{n_J}\chi_\Lambda(y-y_j),\qquad
 \Phi_v(x,y)=\sum_{i=1}^{n_J}\chi_\Lambda(x-x_i),
\]
and prescribe
\begin{equation}\label{eq:woven-raw-metric}
 \boldsymbol M^{\rm raw}_{n_J,\Lambda}(x,y)
 =
 \begin{bmatrix}
 1+(\Lambda^2-1)\Phi_v(x,y)&0\\
 0&1+(\Lambda^2-1)\Phi_h(x,y)
 \end{bmatrix}.
\end{equation}
The horizontal and vertical arms request orthogonal anisotropy, while each
intersection creates a narrow nearly isotropic core.  For every
\((n_J,\Lambda)\) we apply a scalar normalization
\[
 \int_{(0,1)^2}\sqrt{\det\boldsymbol M_{n_J,\Lambda}}\,d\boldsymbol x=1.
\]
Thus increasing \(\Lambda\) sharpens localization and compatibility demands
without merely increasing total metric volume.

All boundary nodes are fixed.  Both densities use the same \(80\times80\)
checkerboard triangulation, quadrature, normalization,
\(\tau=0.02\), line search, nonlinear tolerance, and
BDF/direct-secant realization.  A run is counted as completed only if it
reaches all twenty stages of length \(0.05\) and preserves
\(\det\widehat E_K>0\) at every accepted state; internal BDF steps may decrease
to \(2\times10^{-8}\), with at most \(500\) attempts per stage.

\begin{figure}[!htbp]
\centering
\subfloat[mesh, \(\Lambda=32\)]{
\includegraphics[width=0.28\textwidth]{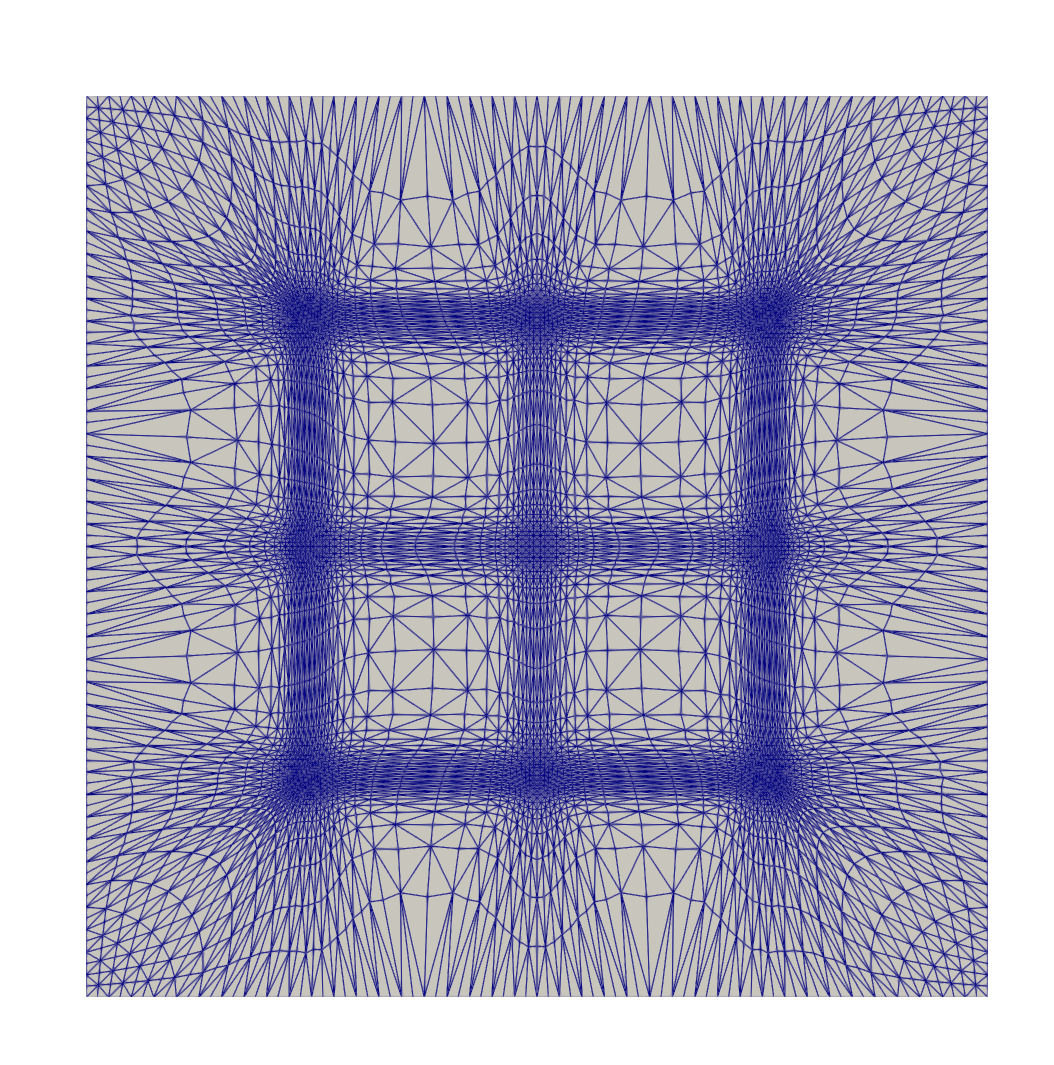}}
\subfloat[mesh, \(\Lambda=96\)]{
\includegraphics[width=0.28\textwidth]{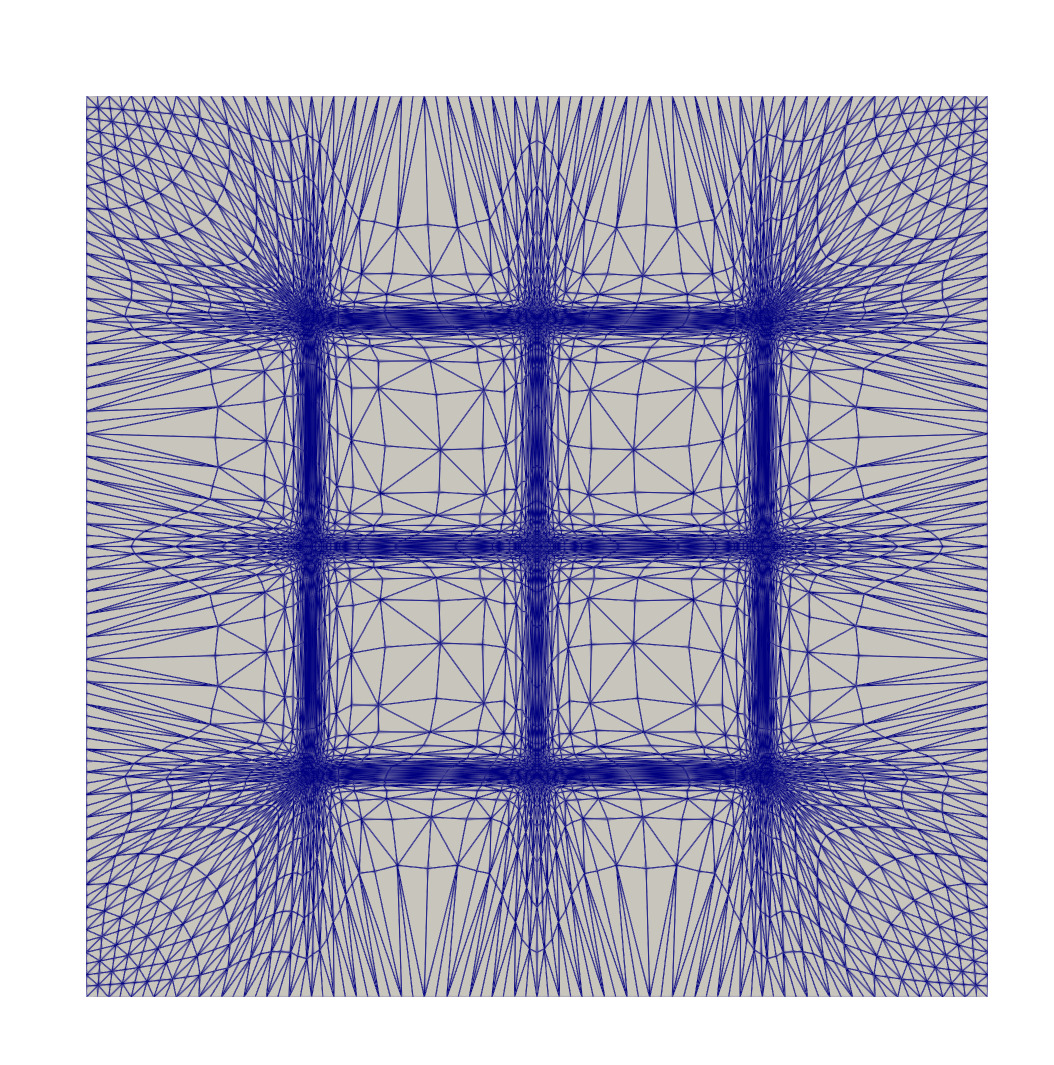}}
\subfloat[mesh, \(\Lambda=160\)]{
\includegraphics[width=0.28\textwidth]{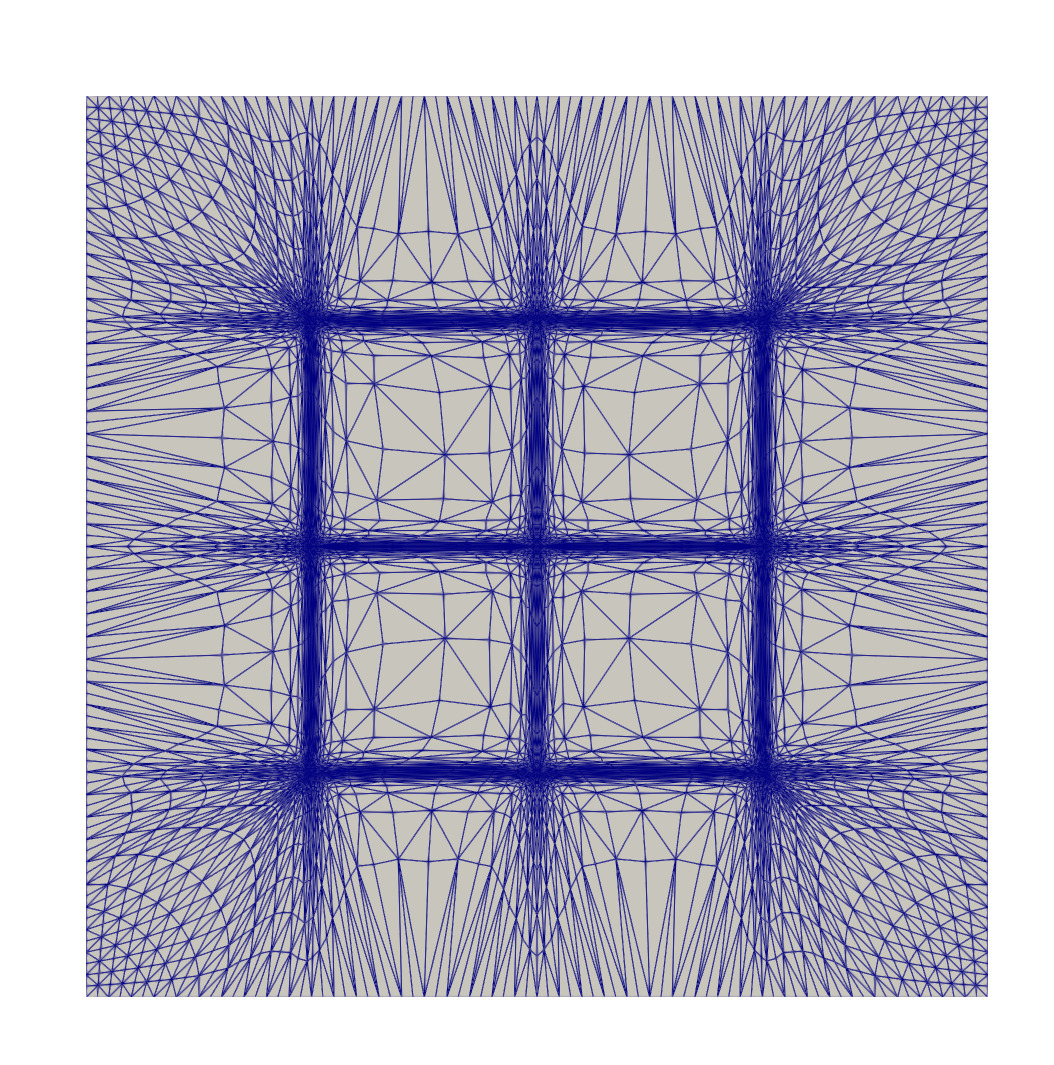}}
\caption{Trace--log meshes for the nine-junction (\(n_J^2=9\)) woven metric \cref{eq:woven-raw-metric} on the \(80\times80\) triangulation.  The panels show a baseline state (\(\Lambda=32\)), the compression level at which no scanned Huang parameter completes for this junction count (\(\Lambda=96\)), and the largest completed trace--log compression level (\(\Lambda=160\)), all under the common fixed protocol.}
\label{fig:woven-j3-visuals}
\end{figure}

Figure~\ref{fig:woven-j3-visuals} displays the nine-junction case at three representative compression levels.

For Huang's functional, the complete scan
\[
 \mathcal U_0=\{0.05,0.10,0.20,\tfrac13,0.45,0.49\}
\]
is performed through \(\Lambda=96\).  At \(\Lambda=128,160,192\), the
four larger \(\mu\)-values have already failed on the preceding continuation
levels.  The compression scan therefore continues with
\(\mathcal U_1=\{0.05,0.10\}\), the only two values that still provide
reference trajectories at lower compression levels.

\begin{table}[!htbp]
\centering
\small
\caption{Completion map for the woven metric \cref{eq:woven-raw-metric}.  A Huang entry \(r/s\) is the number of completed \(\mu\)-values among those tested: \(s=6\) for \(\Lambda\le96\) and \(s=2\) for \(\Lambda\ge128\).  A check mark/cross denotes completion/non-completion for the trace--log density under the same protocol.}
\label{tab:woven-completion-map}
\setlength{\tabcolsep}{5.0pt}
\renewcommand{\arraystretch}{1.10}
\begin{tabular}{c cc cc cc}
\hline
& \multicolumn{2}{c}{\(n_J^2=4\)}
& \multicolumn{2}{c}{\(n_J^2=9\)}
& \multicolumn{2}{c}{\(n_J^2=16\)}\\
\cline{2-7}
\(\Lambda\) & Huang & TL & Huang & TL & Huang & TL\\
\hline
16  & 6/6 & \checkmark & 6/6 & \checkmark & 6/6 & \checkmark \\
32  & 3/6 & \checkmark & 5/6 & \checkmark & 6/6 & \checkmark \\
64  & 1/6 & \checkmark & 3/6 & \checkmark & 6/6 & \checkmark \\
96  & 0/6 & \checkmark & 0/6 & \checkmark & 2/6 & \checkmark \\
128 & 0/2 & \checkmark & 0/2 & \checkmark & 0/2 & \checkmark \\
160 & 0/2 & \checkmark & 0/2 & \checkmark & 0/2 & \checkmark \\
192 & 0/2 & \checkmark & 0/2 & \(\times\) & 0/2 & \(\times\) \\
\hline
\end{tabular}
\end{table}

Among the scanned Huang parameters, the largest completed compression levels are
\(64,64,96\).  The trace--log density completes up to \(192,160,160\) for
\(4,9,16\) junctions, respectively.  The corresponding ratios of the largest completed levels are \(3.00,2.50,1.67\).

\paragraph{Compression-response diagnostic}
The completion map is supplemented by a local diagnostic for Huang runs that reach the per-stage attempt cap while their final accepted states remain positively oriented.  Let \(\boldsymbol x_K\) be the barycenter of \(K\) and define the junction-core set
\[
\mathcal C_{\rm core}:=\left\{K:\ \exists\,i,j\ \text{such that}\ |x_K-x_i|\le w_\Lambda\ \text{and}\ |y_K-y_j|\le w_\Lambda\right\}.
\]
For \(K\in\mathcal C_{\rm core}\), set
\[
 q_K=\frac{\theta_h^{d/2}}{\sqrt{a_K}},
 \qquad a_K=\det\boldsymbol A_K.
\]
Thus \(q_K>1\) indicates compression below the metric-calibrated scale.  Under the virtual isotropic compression \(\widehat E_K(\varepsilon)=e^{-\varepsilon/d}\widehat E_K\), define \(\Pi_{F,K}=\left.dT_F/d\varepsilon\right|_{\varepsilon=0}\).  For \(p=d\gamma/2\),
\begin{equation}\label{eq:woven-compression-slopes}
 \Pi_{{\rm TL},K}
 =\gamma\!\left[(d\theta_h)^p-[\operatorname{tr}(\boldsymbol A_K)]^p\right],
 \qquad
 \Pi_{{\rm H},K}=-\gamma T_{{\rm H},K}<0,
\end{equation}
where \(T_{{\rm H},K}=\mu[\operatorname{tr}(\boldsymbol A_K)]^p+(1-2\mu)d^p a_K^{\gamma/2}\).  The two densities are evaluated at the same local state; the compression ray is used only to identify their local response.

\begin{table}[!htbp]
\centering
\small
\caption{Core-compression evidence for Huang runs terminating at the per-stage attempt cap. \(\Lambda^-\) is the largest completed compression level and \(\Lambda^+\) the first subsequently tested attempt-cap event for the same \(n_J\) and \(\mu\). \(q_{.95}^{\rm last}\) is the 95th percentile over \(\mathcal C_{\rm core}\) at the final accepted state.  \(N_{\rm rej}^{\rm last}\) and \(h_{\rm last}\) refer to the terminal BDF stage; all listed paths retain \(\det\widehat E_K>0\) at that state.}
\label{tab:woven-response}
\setlength{\tabcolsep}{5.2pt}
\renewcommand{\arraystretch}{1.10}
\begin{tabular}{c c c c c c}
\hline
\(n_J^2\) & \(\mu\) & \(\Lambda^-\!\to\!\Lambda^+\)
& \(q_{.95}^{\rm last}\)
& \(N_{\rm rej}^{\rm last}\) & \(h_{\rm last}\) \\
\hline
4  & .05 & \(64\to96\)  & \(1.89\to5.08\) & 73  & \(5.7\!\times\!10^{-7}\) \\
9  & .05 & \(64\to96\)  & \(2.70\to5.46\) & 496 & \(2.0\!\times\!10^{-8}\) \\
16 & .05 & \(96\to128\) & \(4.05\to7.09\) & 489 & \(2.0\!\times\!10^{-8}\) \\
\hline
\end{tabular}
\end{table}

Along these paths, the core compression quantile rises while orientation is
retained.  Together with \cref{eq:woven-compression-slopes}, this gives a
local response explanation for the broader completion range of the trace--log
density in the compression parameter under the common finite-budget protocol.


\subsection{Numerical Solution of the Two-Dimensional Burgers' Equation}\label{burgers}
\begin{example}\label{exp:3}
We consider the two-dimensional scalar Burgers' equation
\begin{equation}
\left\{
\begin{aligned}
&\frac{\partial u}{\partial t} + u u_x + u u_y = \frac{1}{Re} \Delta u,
&&\text{in } \Omega=[0,1]^2,\\
&u=g, &&\text{on } \partial\Omega,\\
&u(\boldsymbol{x},0)=u_0,
\end{aligned}
\right.
\end{equation}
which admits the exact solution
\[
 u_{\rm ex}(x,y,t)
 =\frac{1}{1+\exp\!\left(\frac{Re}{2}(x+y-t)\right)}.
\]
We set $Re=200$, $g=u_{\rm ex}|_{\partial\Omega}$, and
$u_0=u_{\rm ex}(\cdot,0)$; direct substitution verifies the equation.
Time integration is performed on $T=[0,2]$ using a second-order singly diagonally implicit Runge--Kutta scheme \cite{ascher1998computer}
with uniform step size $0.002$ and a moving mesh of $3600$ elements.
For each of the Hessian-based metric \cref{HessianMetric} and the
arc-length-type metric \cref{ArcLengthMetric}, we run both the proposed and
Huang functionals.  The PDE discretization, physical time step, initial physical mesh, boundary
treatment, metric definition, common metric-smoothing schedule, and remeshing
schedule are identical within these four Burgers configurations.  The inner
artificial-time parameters remain the functional-specific values stated in
the experimental setup.
\end{example}

\begin{figure}[!htbp]
    \centering
    \subfloat[$t=0.5$s, Hessian]{
        \includegraphics[width=0.225\textwidth]{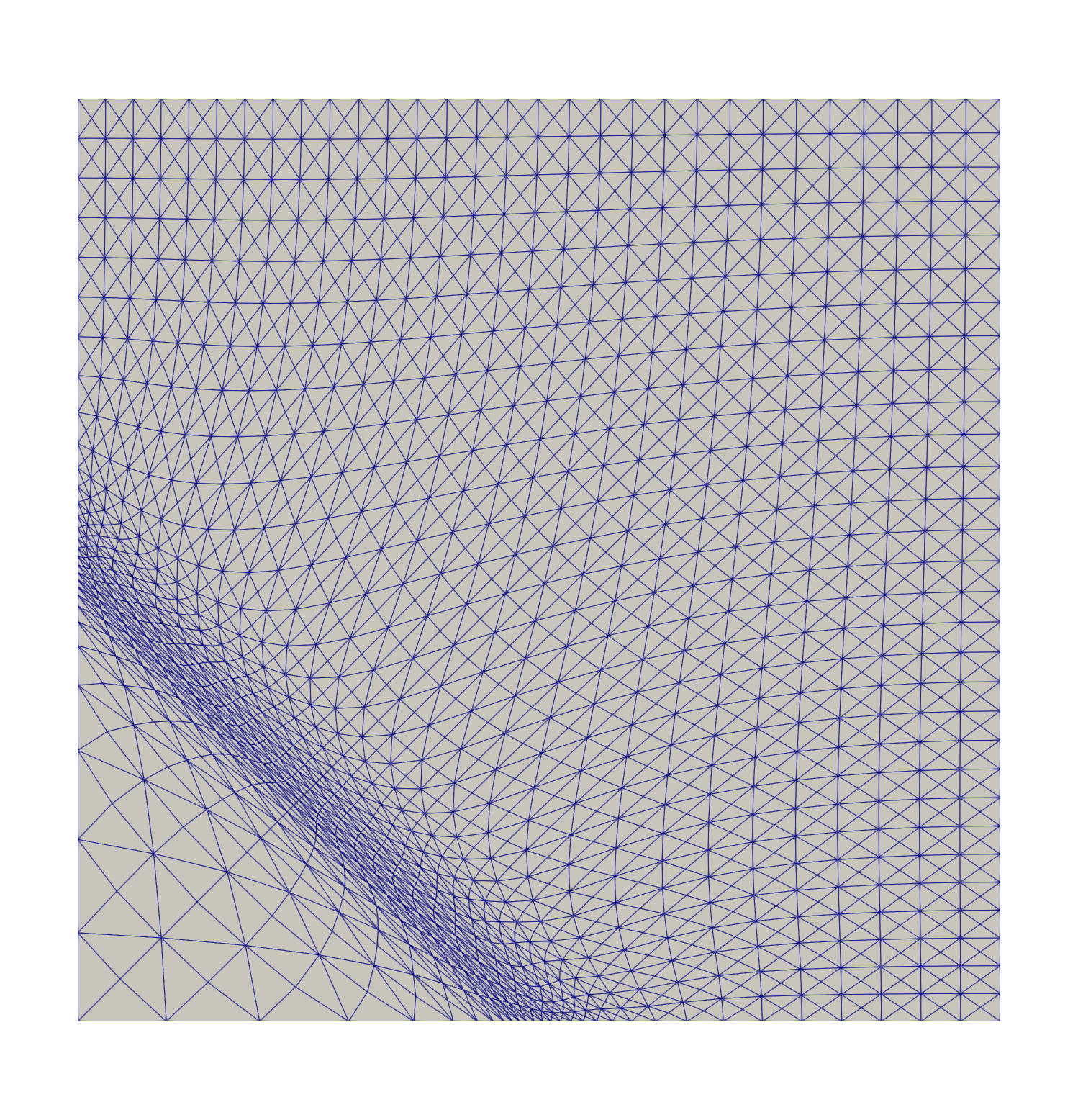}
        \label{Fig7.sub1}
    }
    \subfloat[$t=1.0$s, Hessian]{
        \includegraphics[width=0.225\textwidth]{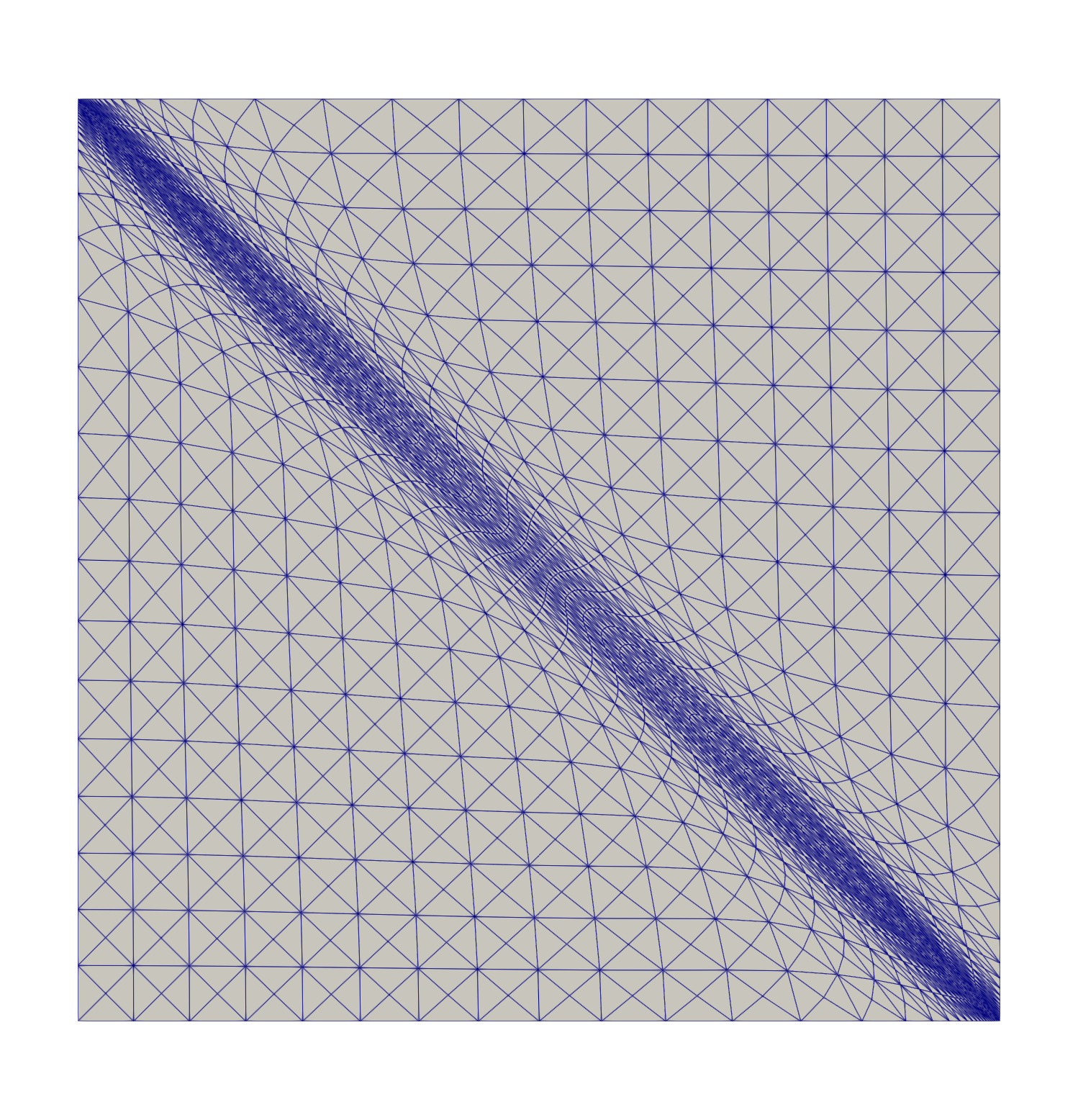}
        \label{Fig7.sub2}
    }
    \subfloat[$t=1.5$s, Hessian]{
        \includegraphics[width=0.225\textwidth]{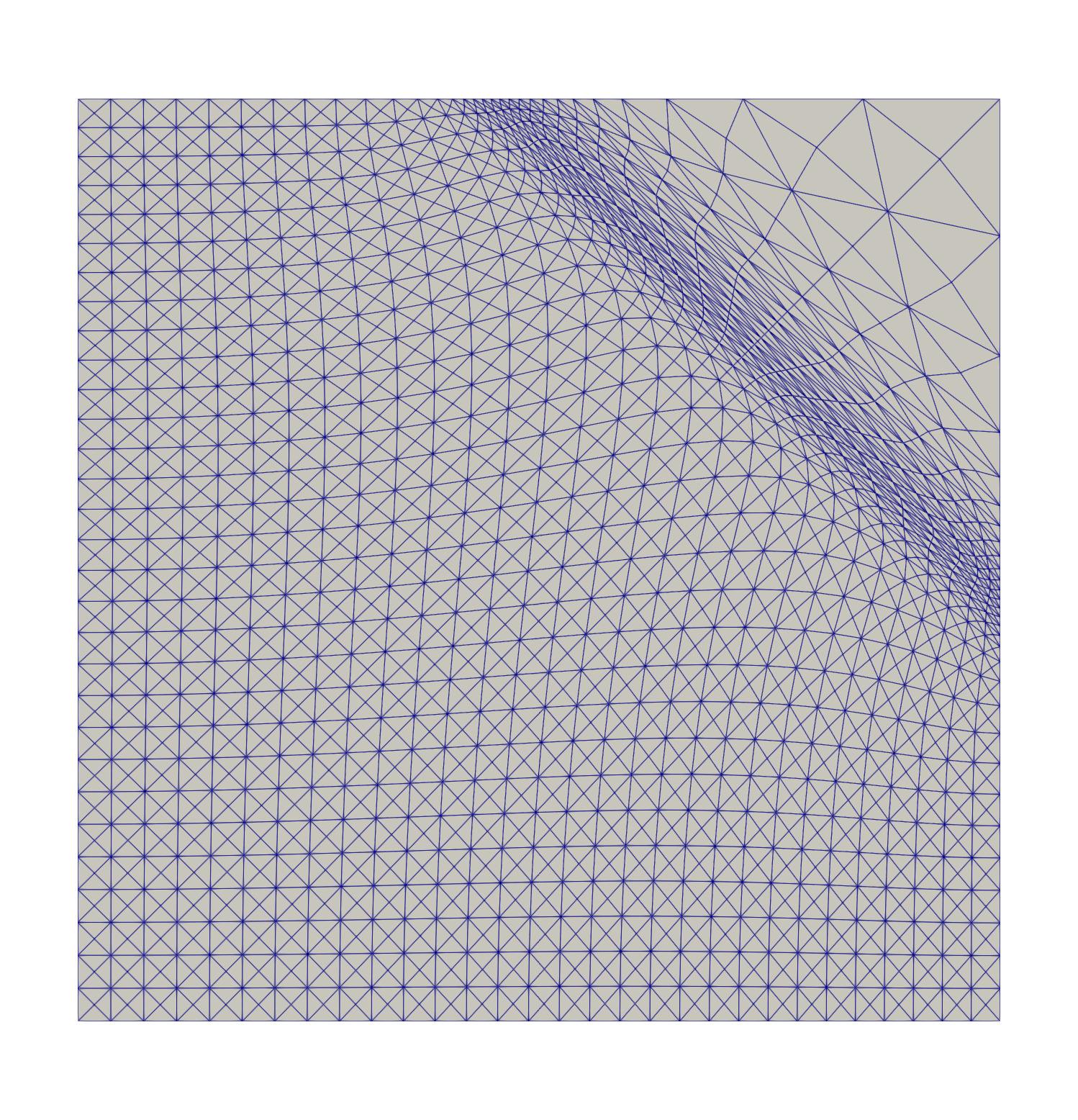}
        \label{Fig7.sub3}
    }

    \vspace{0.1em}

    \subfloat[$t=0.5$s, Arc-length]{
        \includegraphics[width=0.225\textwidth]{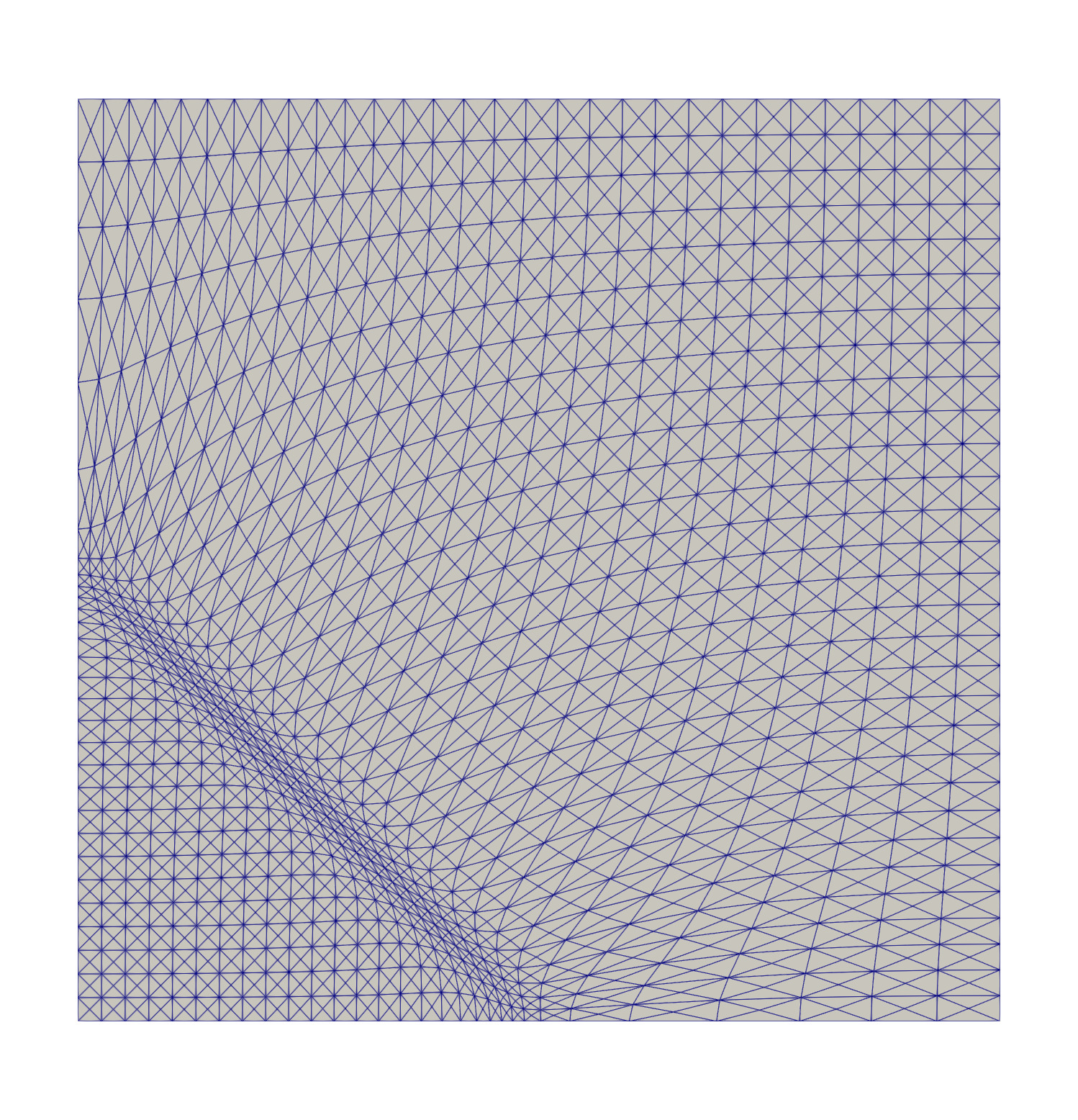}
        \label{Fig7.sub4} 
    }
    \subfloat[$t=1.0$s, Arc-length]{
        \includegraphics[width=0.225\textwidth]{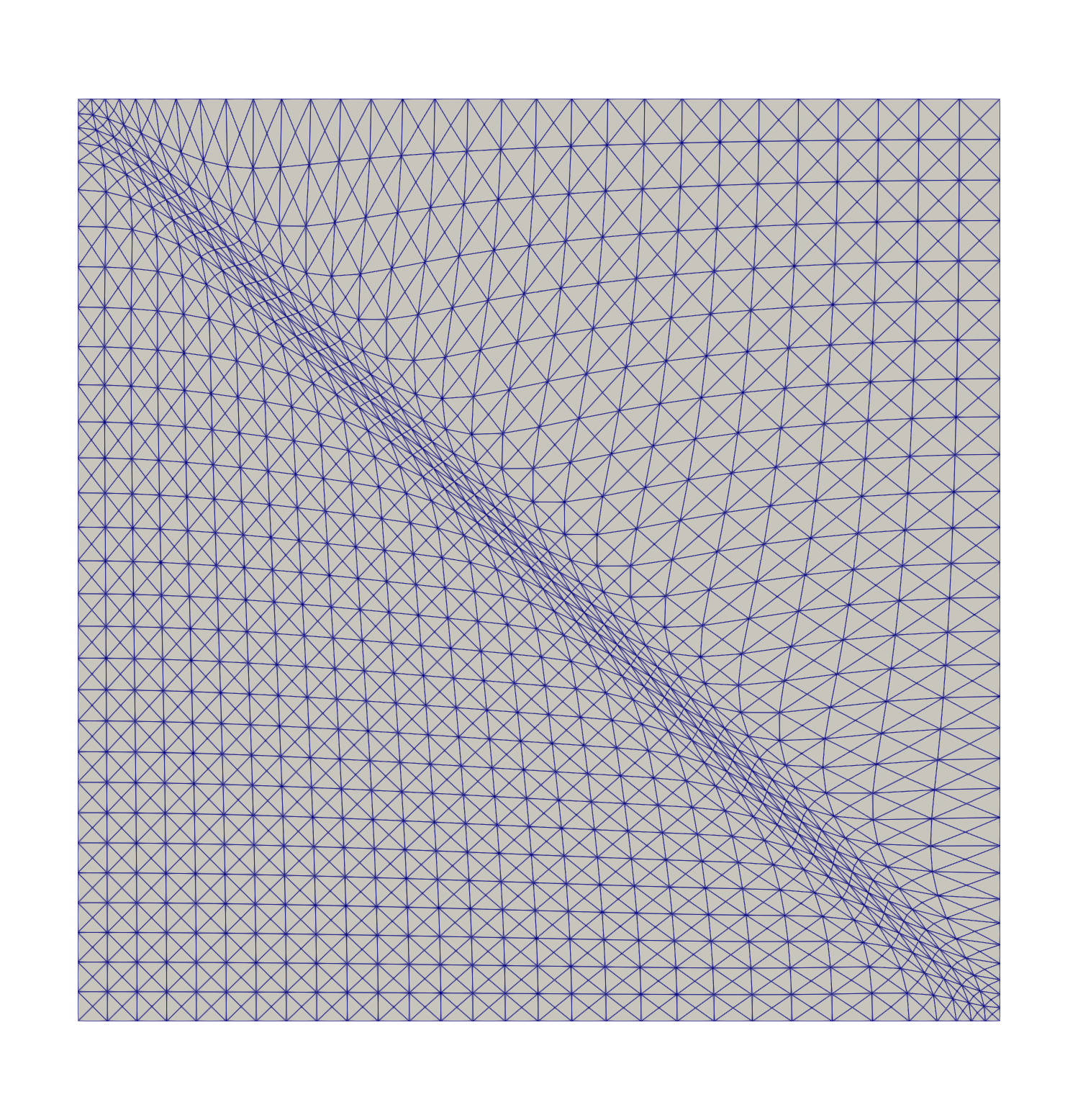}
        \label{Fig7.sub5} 
    }
    \subfloat[$t=1.5$s, Arc-length]{
        \includegraphics[width=0.225\textwidth]{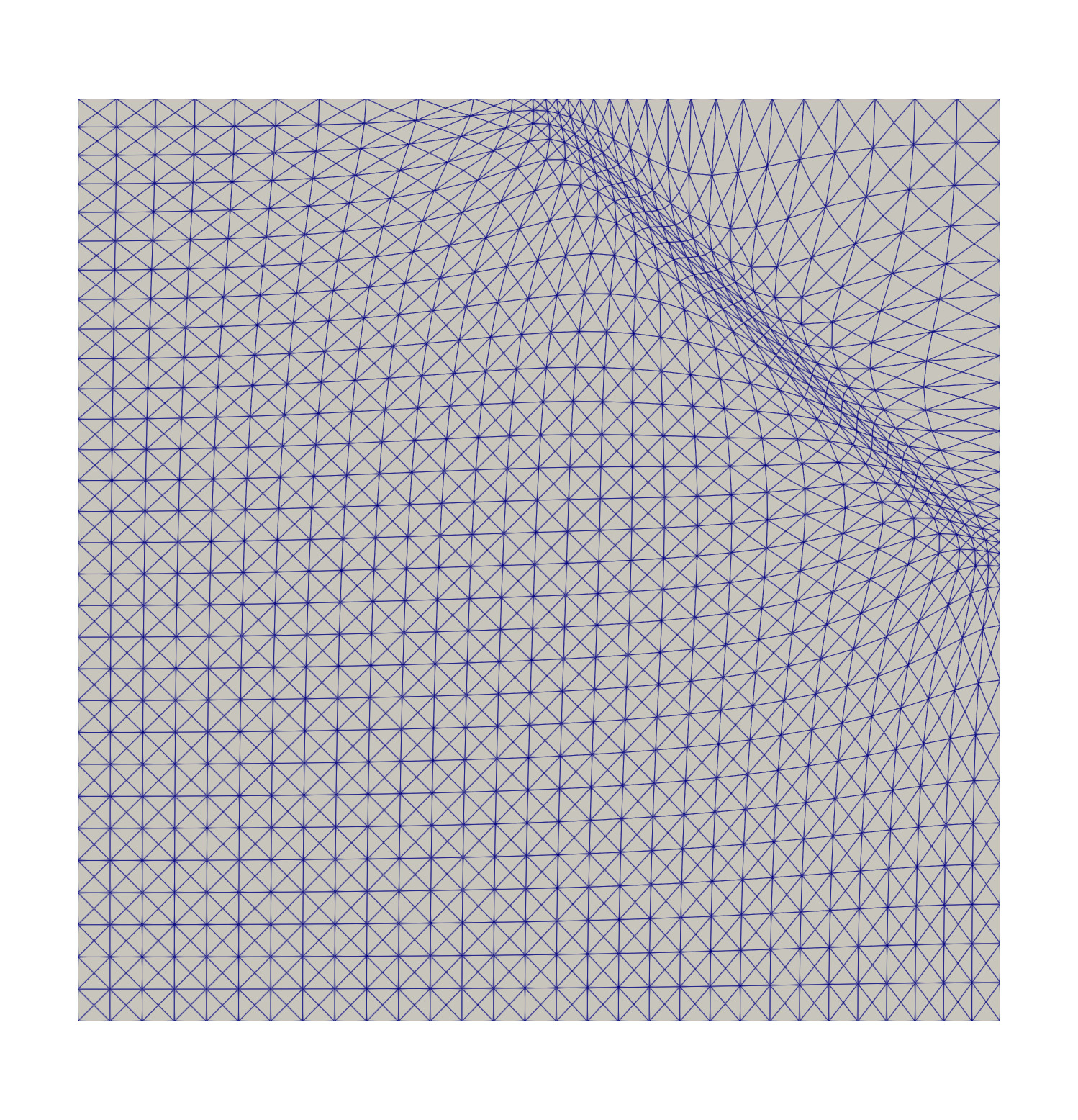}\label{Fig7.sub6} 
    }

    \caption{
    Mesh evolution produced by the trace--logarithmic functional for the Burgers' equation.
    Top row: Hessian-based metric.
    Bottom row: Arc-length-based metric.
    From left to right: $t=0.5$s, $1.0$s, and $1.5$s.
    }
    \label{Fig7.main} 
\end{figure}

\cref{Fig7.main} shows representative mesh distributions produced by the
proposed functional at $t=0.5$, $1.0$, and $1.5$.
Both metrics align and concentrate the mesh along the characteristic line
$x+y=t$, with tangential stretching and normal compression.  The
Hessian-based metric produces the more concentrated anisotropic mesh.

\begin{figure}[!htbp]
        \centering 
        \subfloat[$L_2$ error comparison.]{
            \includegraphics[width=0.38\textwidth]{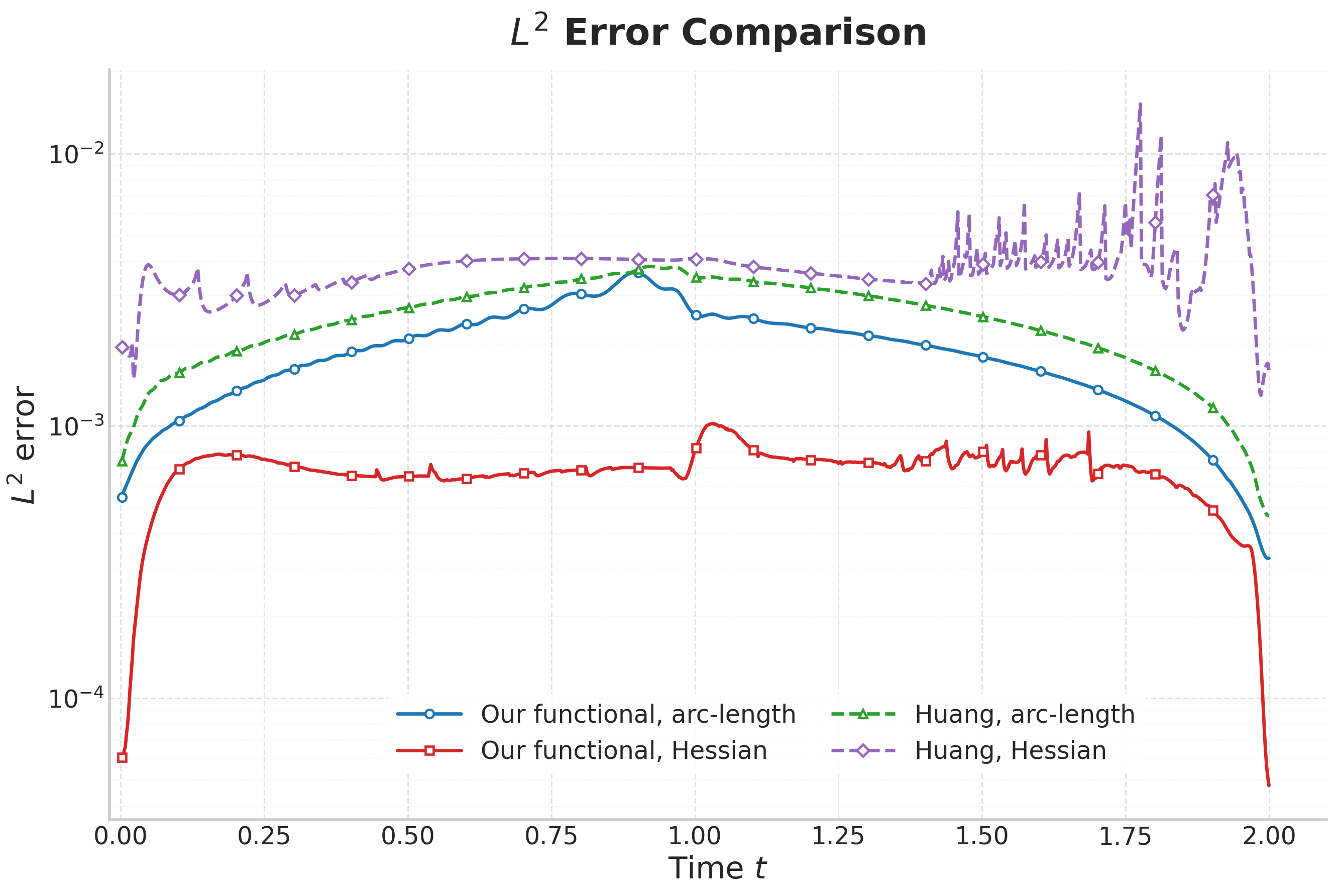}
            \label{Fig9.sub1} 
        }
		\subfloat[$H_1$ error comparison.]{
            \includegraphics[width=0.38\textwidth]{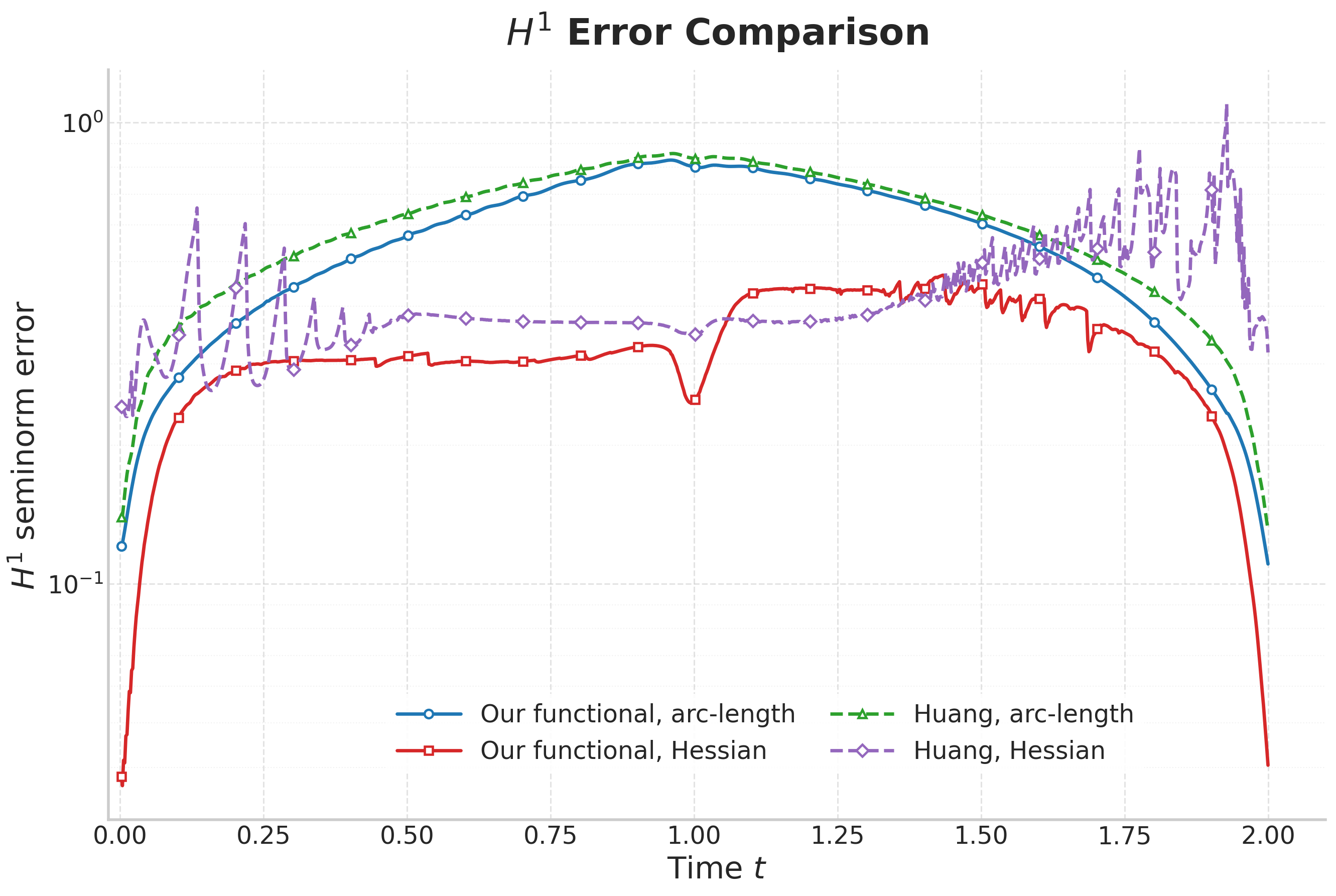}
            \label{Fig9.sub2}
        }
        \caption{Time histories of the \(L_2\) (a) and \(H_1\) (b) errors for the four functional--metric configurations.  The legend identifies both the functional and the metric used by each curve.}
        \label{Fig9.main} 
\end{figure}

\cref{Fig9.main} compares all four functional--metric configurations.
For the arc-length metric, the proposed functional remains below Huang's
functional in both error norms over the reported time interval.  The Hessian comparison applies the stronger monitor to both densities under the same physical discretization and remeshing schedule.  In the reported implementation, the trace--log/Hessian branch has the smallest and smoothest error histories, whereas the Huang--Hessian branch has larger and more nonmonotone errors.  This is a consistent implementation-level separation under the stated functional-specific artificial-time calibrations.

The quantitative data in \cref{tab:burgers-error-history} support the same
implementation-level conclusion.  At every tabulated time, the trace--log/Hessian configuration gives the smallest \(L_2\) and \(H_1\) errors among the four reported runs; the Huang--Hessian branch remains less accurate throughout the reported interval.

\begingroup
\setlength{\tabcolsep}{2.8pt}
\renewcommand{\arraystretch}{0.8}
\begin{table}[!htbp]
\centering
\caption{$L_2$ and $H_1$ errors for the four functional--metric Burgers configurations (Example 5.3).}
\label{tab:burgers-error-history}
\begin{tabular}{l c c c c c c}
\hline
functional & Error type & \multicolumn{5}{c}{time \(t\)} \\
\cline{3-7}
 &  & 0 & 0.5 & 1.0 & 1.5 & 2.0 \\
\hline
Trace-log (Arc-Length)    & \multirow{4}{*}{$L_2$} & 5.46e-04 & 2.08e-03 & 2.56e-03 & 1.79e-03 & 3.26e-04 \\
Huang's (Arc-Length) &                         & 7.43e-04 & 2.71e-03 & 3.51e-03 & 2.52e-03 & 4.66e-04 \\
Trace-log (Hessian)       &                         & 6.06e-05 & 6.52e-04 & 8.05e-04 & 7.94e-04 & 4.76e-05 \\
Huang's (Hessian)     &                        & 1.94e-03 & 3.78e-03 & 4.10e-03 & 3.81e-03 & 1.60e-03 \\
\hline
Trace-log (Arc-Length)    & \multirow{4}{*}{$H_1$} & 1.20e-01 & 5.67e-01 & 8.02e-01 & 6.04e-01 & 1.10e-01 \\
Huang's (Arc-Length) &                         & 1.39e-01 & 6.33e-01 & 8.36e-01 & 6.31e-01 & 1.31e-01 \\
Trace-log (Hessian)       &                         & 3.82e-02 & 3.11e-01 & 2.48e-01 & 4.43e-01 & 4.04e-02 \\
Huang's (Hessian)     &                        & 2.43e-01 & 3.84e-01 & 3.48e-01 & 4.81e-01 & 3.17e-01 \\
\hline
\end{tabular}
\end{table}
\endgroup

\subsection{Rayleigh--Taylor Instability Test}
\label{rayleightaylor}

The Rayleigh--Taylor instability test examines mesh adaptation in a flow with multiscale interface roll-up.  We use the Hessian and eigen-decomposition metrics in \cref{HessianMetric,EigenDecompositionMetric}.

\begin{example}\label{exp:rti}
Following \cite{lee2011long}, we consider the nondimensional
diffuse-interface phase-field model for Rayleigh--Taylor instability,
\[
\left\{
\begin{aligned}
&\rho(\phi)\bigl(\partial_t\boldsymbol u
+(\boldsymbol u\cdot\nabla)\boldsymbol u\bigr)
=-\nabla p+Re^{-1}\Delta\boldsymbol u
+Fr^{-1}\rho(\phi)\boldsymbol g,\\
&\nabla\cdot\boldsymbol u=0,\\
&\partial_t\phi+\nabla\cdot(\phi\boldsymbol u)
=Pe^{-1}\Delta\mu,\\
&\mu=\phi^3-\phi-\varepsilon^2\Delta\phi .
\end{aligned}\right.
\]
Here
\[
\rho(\phi)
=\frac{\rho_\ell(1-\phi)+\rho_h(1+\phi)}{2},
\qquad
\rho_h=2,\quad \rho_\ell=1,\quad
\boldsymbol g=(0,-1)^T .
\]
On \(\Omega_p=[0,1]\times[0,2]\), the initial data are
\[
\boldsymbol u_0=\boldsymbol0,\qquad
\phi_0(x,y)=
\tanh\!\left(
\frac{y-1-0.1\cos(2\pi x)}
{\sqrt2\,\varepsilon}\right).
\]
The boundary conditions are taken as in \cite{lee2011long}.
We use
\(Re=3000\), \(Fr=1\), \(\varepsilon=0.005\), and
\(Pe=100/\varepsilon\).
The equations are discretized by BDF2 with
\(\Delta t=0.005\times0.16\sqrt2\) on an initial mesh of
\(14{,}112\) triangles. Quadratic elements are used for
\(\phi\), \(\mu\), and \(\boldsymbol u\), and linear elements for \(p\).
After the arbitrary Lagrangian--Eulerian reformulation, the mesh velocity in the transport terms is
approximated by
\[
\boldsymbol w_h^{n+1}
=
\frac{
3\boldsymbol x_h^{n+1}
-4\boldsymbol x_h^n
+\boldsymbol x_h^{n-1}}
{2\Delta t},
\]
where \(\boldsymbol x_h^n\) denotes the mesh-node vector at time \(t_n\).
\end{example}

\begin{figure}[!htbp]
    \centering
    \includegraphics[width=0.155\textwidth]{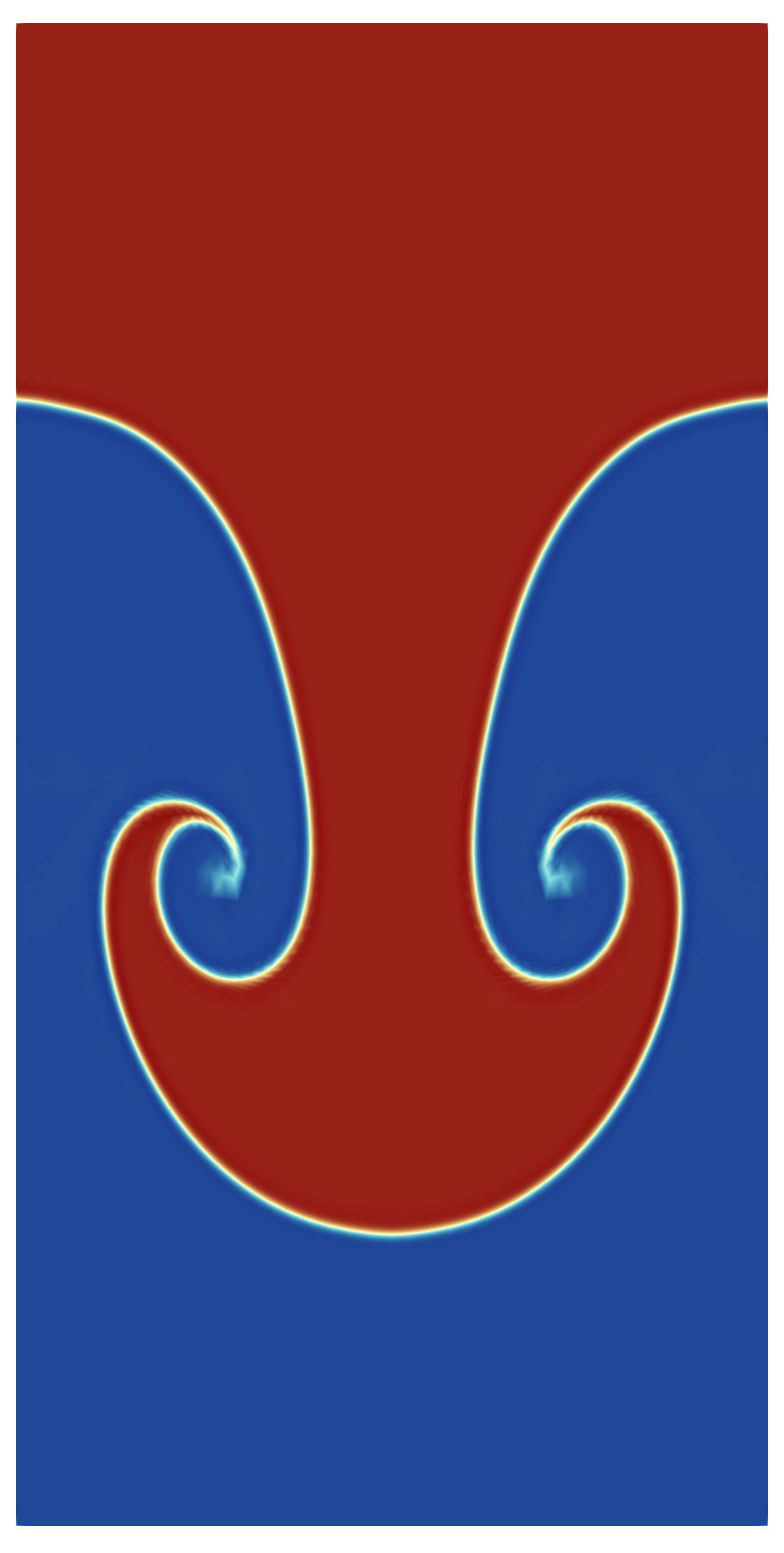}
    \includegraphics[width=0.155\textwidth]{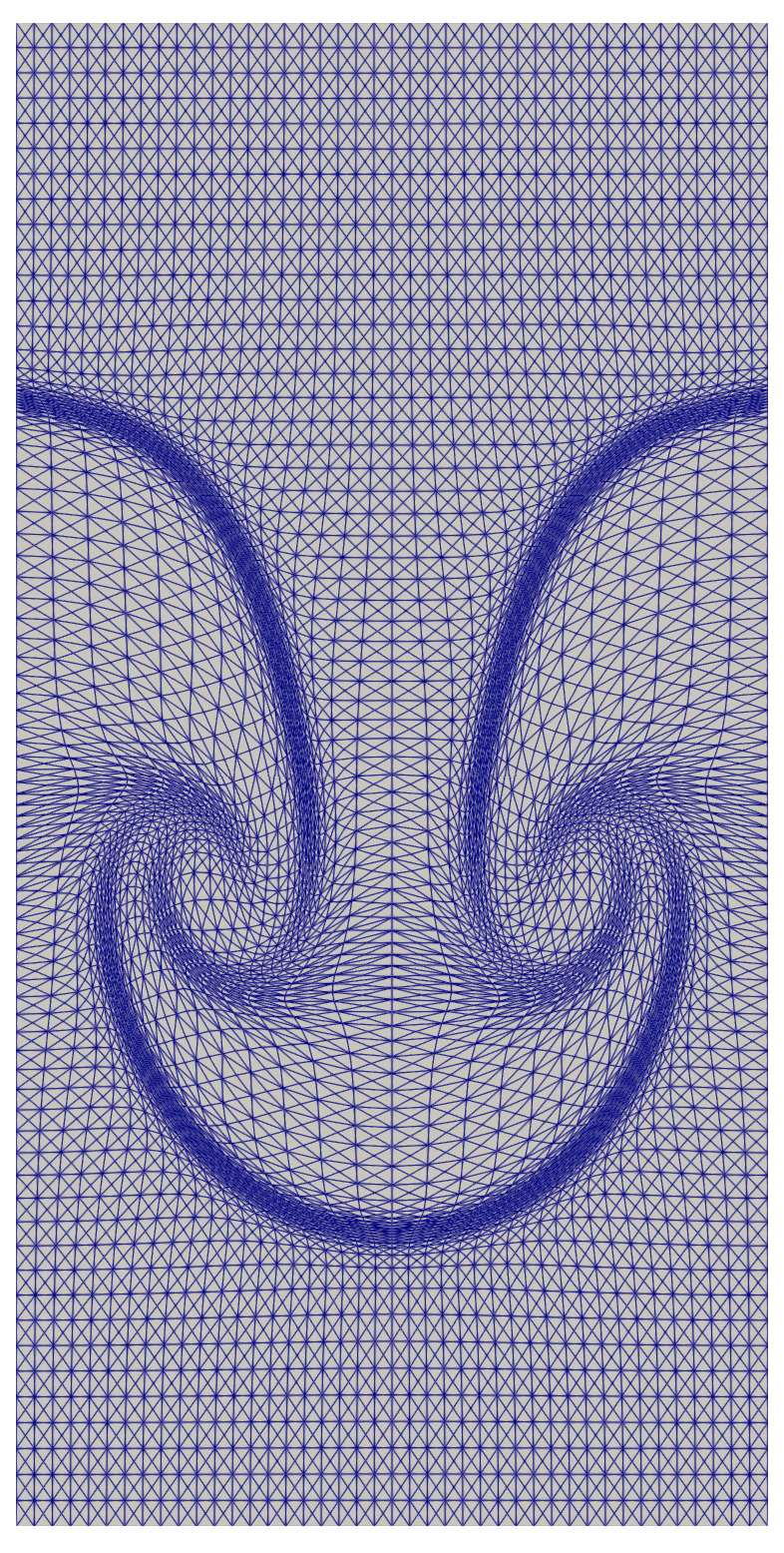}
    \includegraphics[width=0.155\textwidth]{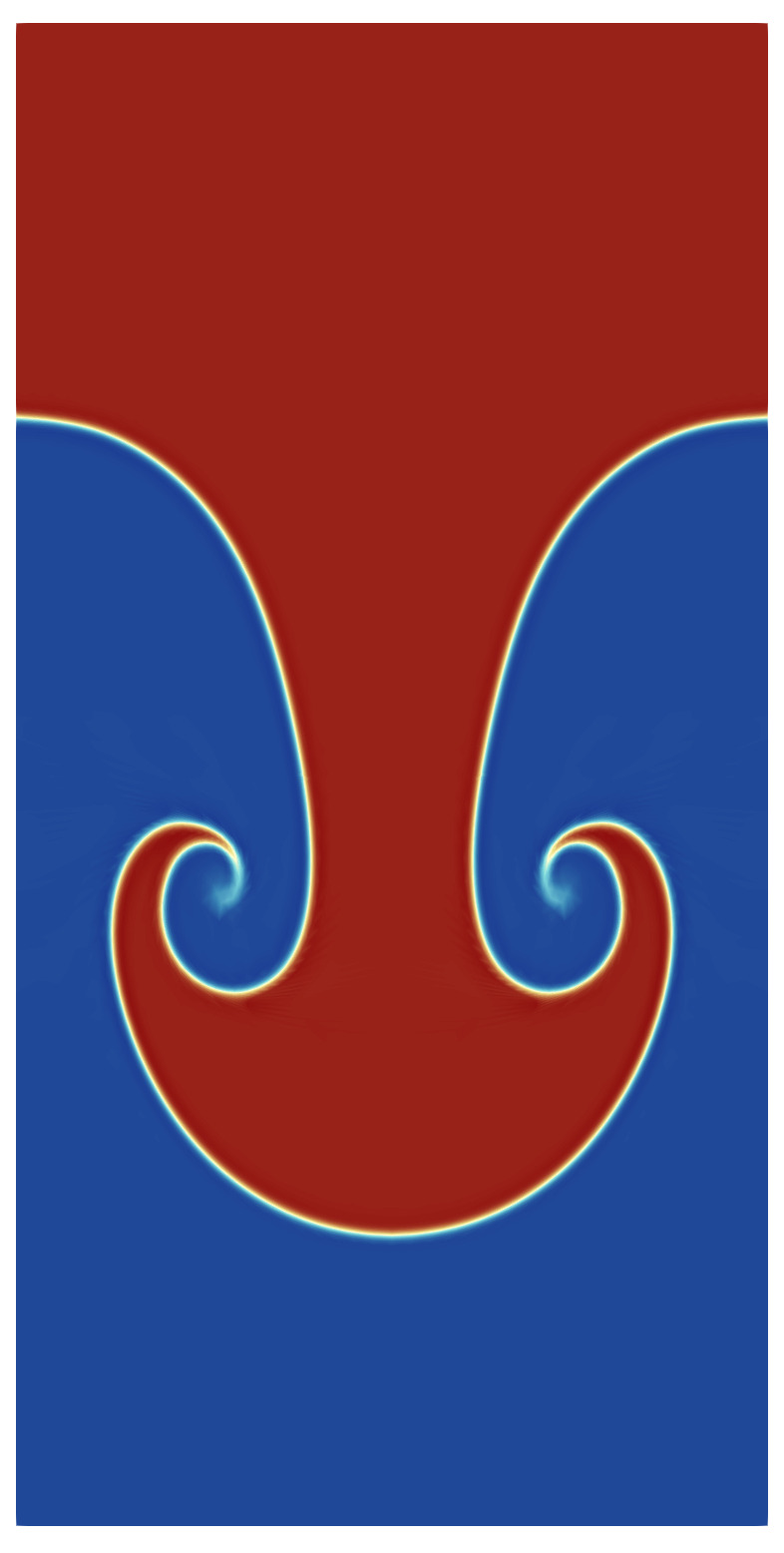}
    \includegraphics[width=0.155\textwidth]{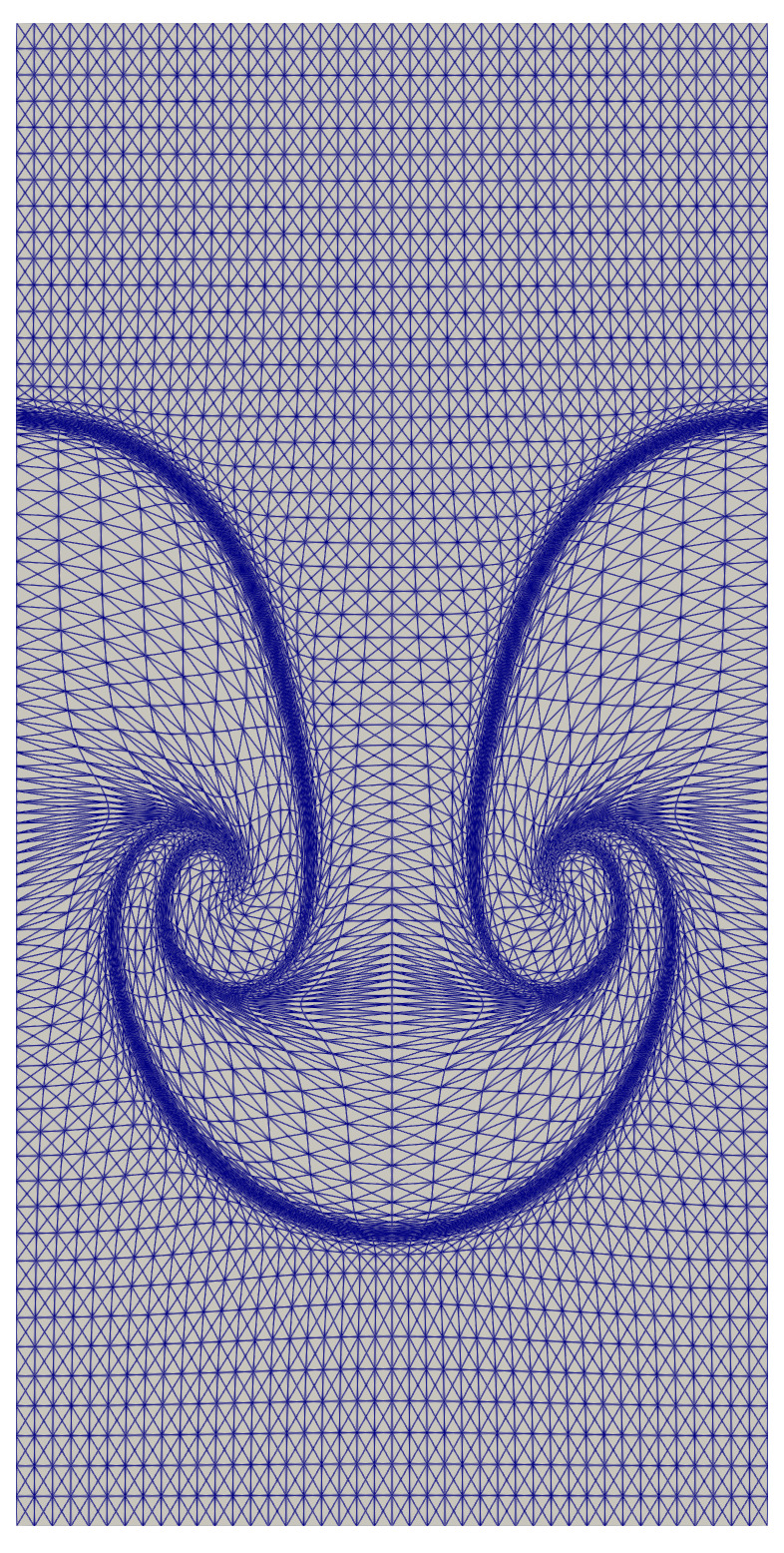}
    \caption{Interface and mesh distribution at $t=2.75$. From left to right: Hessian interface and mesh, followed by eigen-decomposition interface and mesh.}
    \label{Fig10.main}
\end{figure}

\begin{figure}[!htbp]
    \centering
    \includegraphics[width=0.27\textwidth]{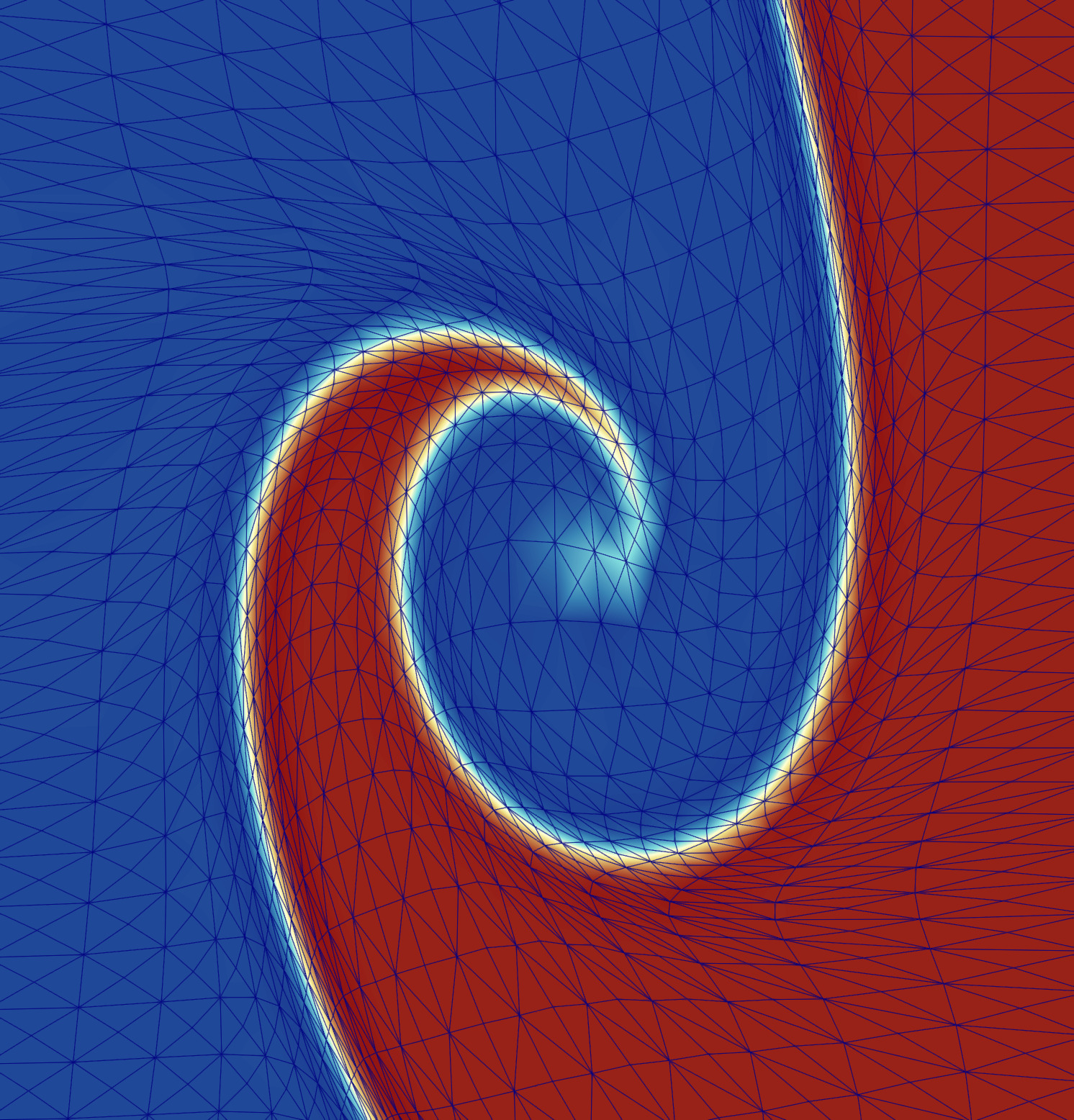}\qquad
    \includegraphics[width=0.27\textwidth]{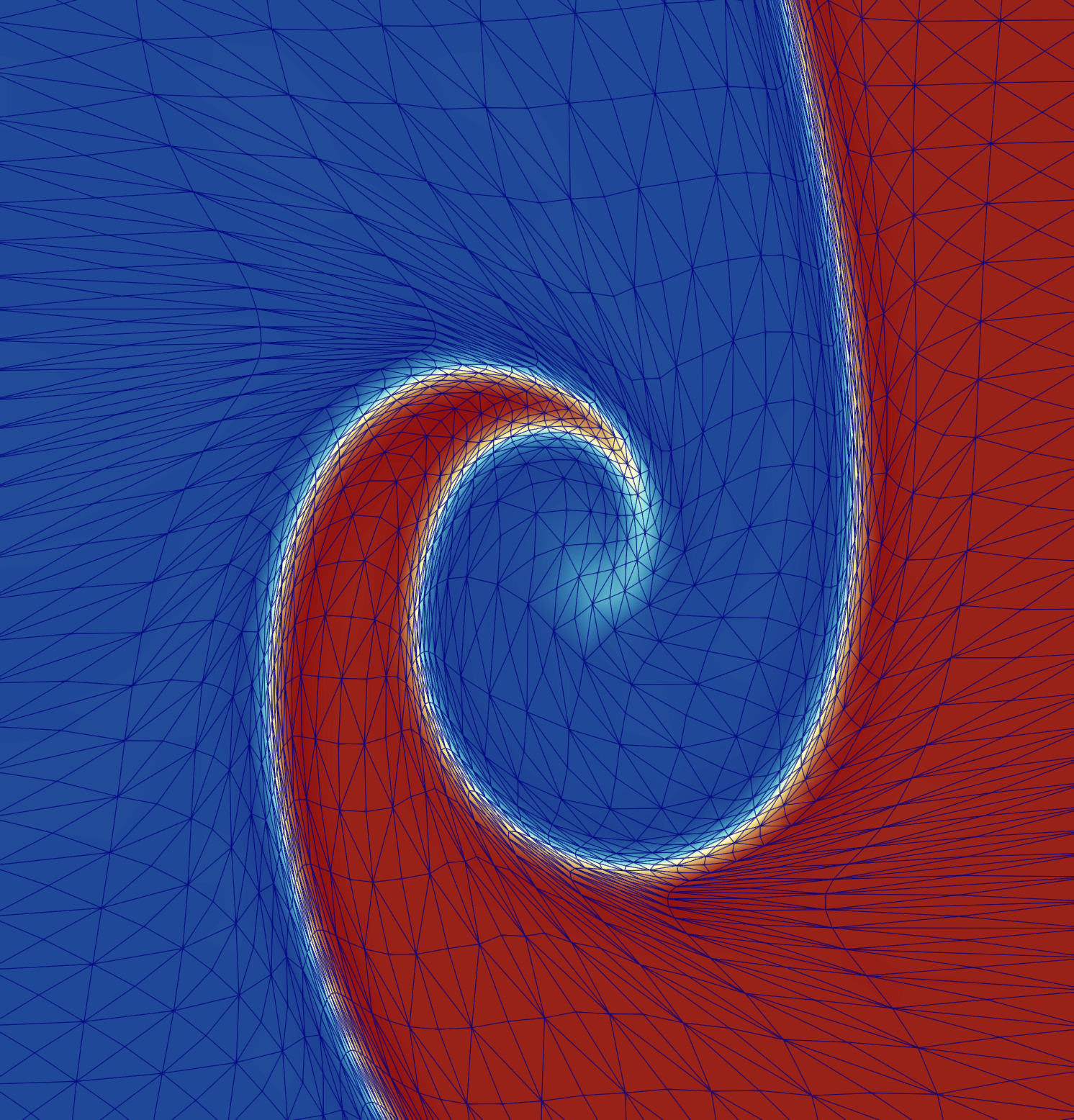}
    \caption{Local roll-up magnifications for the Hessian metric (left) and eigen-decomposition metric (right).}
    \label{Fig11.main}
\end{figure}

Figure~\ref{Fig10.main} shows the interfaces and meshes at $t=2.75$.  Both metrics resolve the nonlinear roll-up.  The eigen-decomposition metric concentrates more nodes near the interface and resolves finer structures in Figure~\ref{Fig11.main}; the Hessian metric gives a smoother distribution.  No interface breakup or visible spurious oscillation appears in the displayed simulations.

\begin{figure}[!htbp]
    \centering
    \includegraphics[width=0.3\textwidth]{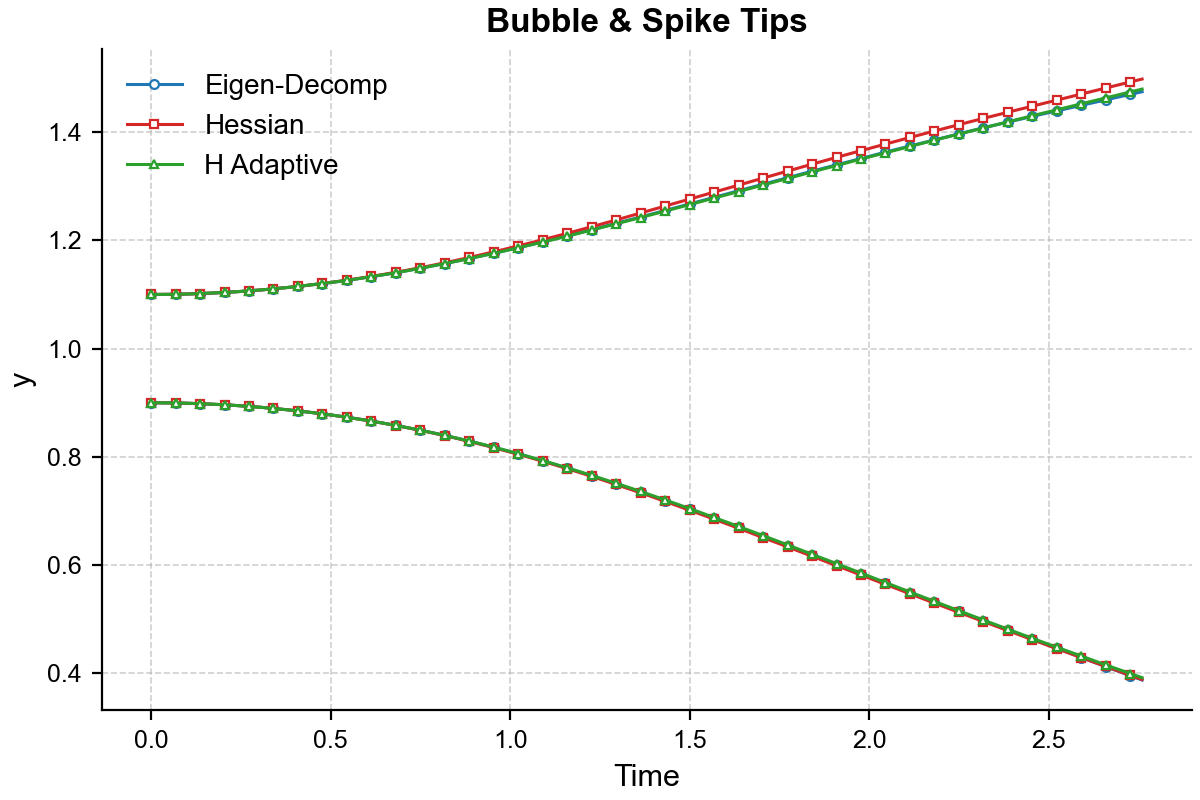}
    \includegraphics[width=0.3\textwidth]{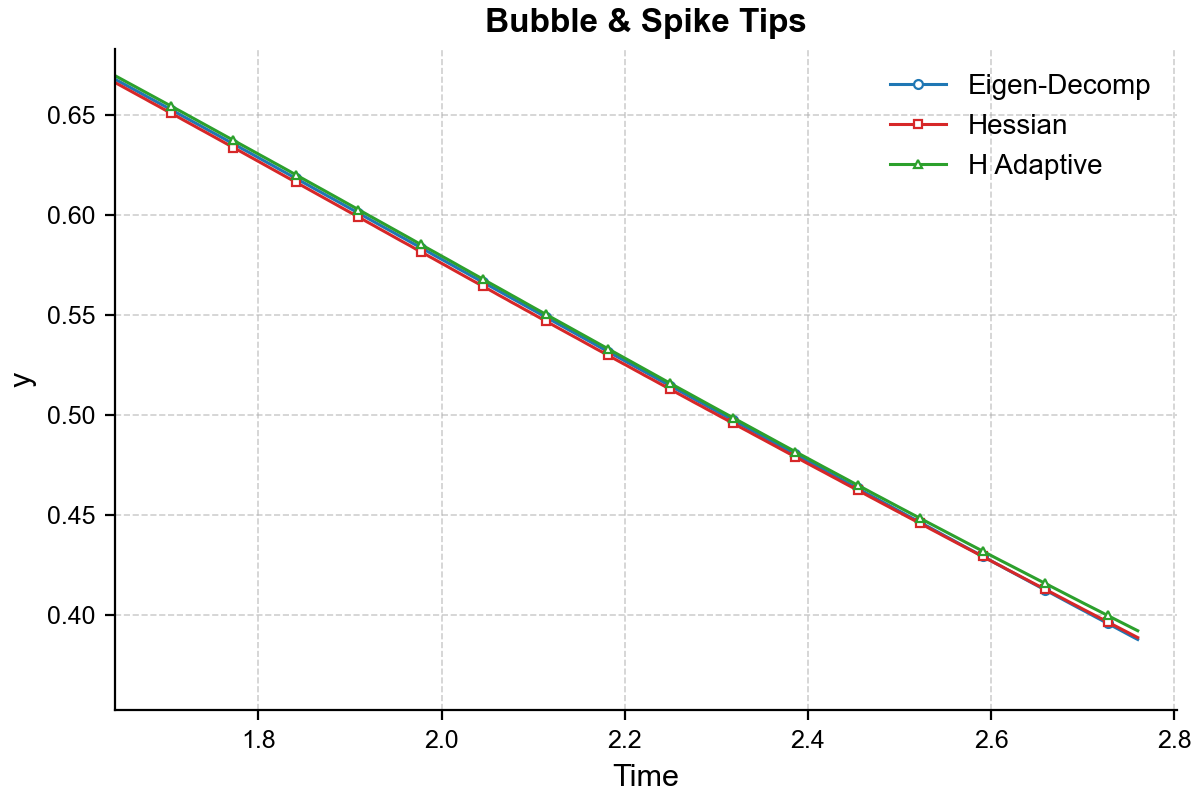}
    \includegraphics[width=0.3\textwidth]{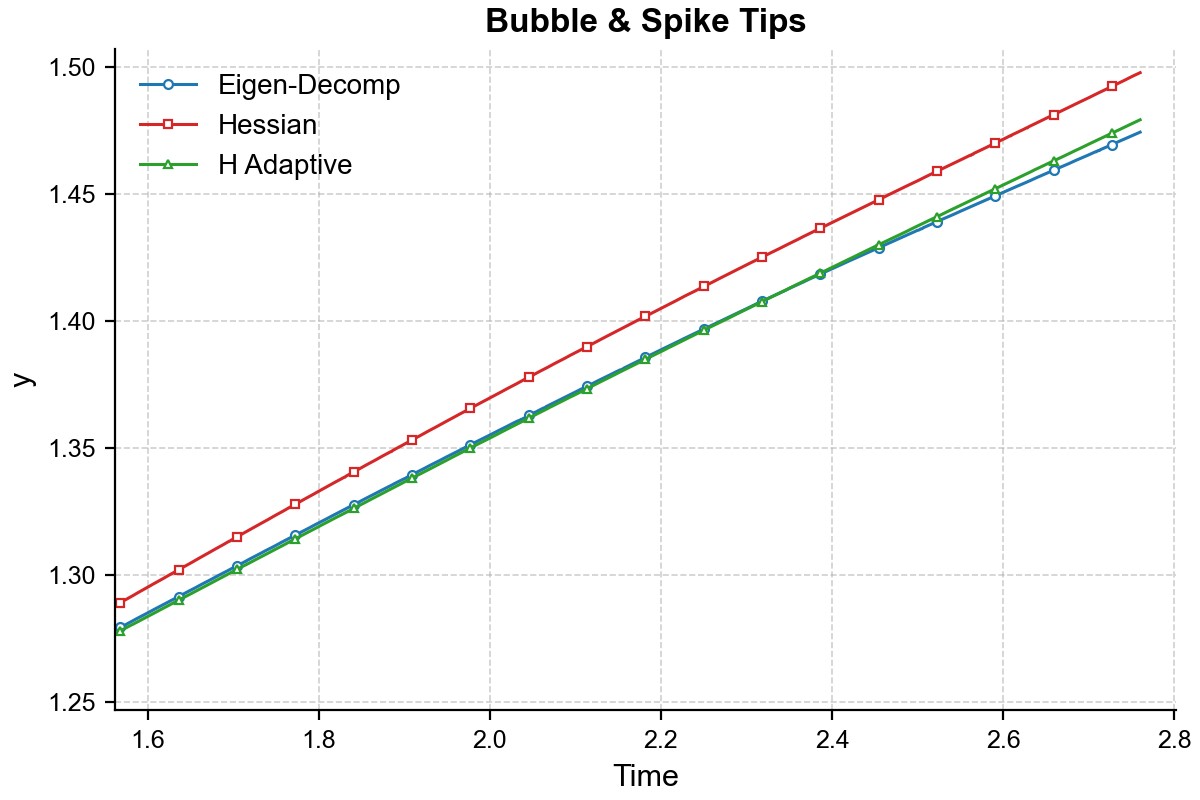}
    \caption{Tip positions versus dimensionless time: full histories (left), descending-fluid zoom (middle), and ascending-fluid zoom (right). Green: reference $h$-adaptive computation; red: Hessian metric; blue: eigen-decomposition metric.}
    \label{Fig12.main}
\end{figure}

Figure~\ref{Fig12.main} compares the vertical positions of the rising and
descending spikes with those obtained from a reference $h$-adaptive
simulation.
The spike trajectories are in close agreement. A quantitative cost comparison with the reference $h$-adaptive computation is deferred until matched tolerance, degree-of-freedom, and wall-time data are reported.


\section{Conclusions and Future Work}

We proposed the calibrated functional \cref{PF}, whose logarithmic coefficient is fixed by local EA stationarity.  The analysis establishes the SPD tensor target, inverse-Jacobian polyconvexity and coercivity, weak minimizer existence, and the Huang--Kamenski condition for semi-discrete physical-coordinate mesh nonsingularity.  The geometric discretization yields the compact computational-coordinate residual used by the BDF/direct-secant solver.

The function-induced and flow tests show stable metric-driven adaptation.  Under the matched finite-budget woven-junction protocol, the trace--log runs complete over a broader compression-parameter range than the scanned Huang benchmark; the largest completed levels differ by factors \(3.00,2.50,1.67\) for the three junction counts.  These comparisons are specific to the stated protocol and scanned parameters.  Extensions beyond simplicial, piecewise-linear geometry remain future work.

\section*{Acknowledgments}
\small
This work was supported by the National Key R\&D Program of China (2024YFA1012600), the National Natural Science Foundation of China (12371410, 12261131501), the 111 Project (D23017), the Program for Science and Technology Innovative Research Team in Higher Educational Institutions of Hunan Province of China, and the High Performance Computing Platform of Xiangtan University.  Generative AI tools assisted manuscript preparation; all mathematical content, computations, numerical results, references, and final text were checked by the authors.
\normalsize

\appendix
\section{Code Availability and Implementation Details}
\small
The implementation is available at \url{https://github.com/weihuayi/fealpy.git} (branch \texttt{develop}).  The moving-mesh routines are in \path{fealpy.mmesh.metrictensoradaptive}.  Representative scripts in \path{example/meshopt/mmesh} are \path{mmesher_with_solution.py}, \path{woven_junction_network.py}, \path{scalar_burgers_mm_solver.py}, and \path{rayleightaylor_mm_solver.py} (Sections~\ref{with_functions}--\ref{rayleightaylor}).

\section{Elementwise Differentiation for the Base Tangent}\label{app:base-tangent}
Let \(H_K=D\widehat E_K[\delta\boldsymbol\xi_f]\).  The local differential for \cref{eq:base-tangent} is
\[
D\mathcal G_K^\xi[H_K]=2\rho_K\boldsymbol R\widehat E_K^{-1}
\bigl(D\boldsymbol Q_K[H_K]-H_K\widehat E_K^{-1}\boldsymbol Q_K\bigr),
\]
with
\[
\begin{aligned}
D\boldsymbol P_K&=H_KE_K^{-1},\qquad
D\boldsymbol A_K=D\boldsymbol P_K\boldsymbol M_K^{-1}\boldsymbol P_K^T
+\boldsymbol P_K\boldsymbol M_K^{-1}D\boldsymbol P_K^T,\\
D\boldsymbol Q_K&=p s_K^{p-1}D\boldsymbol A_K
+p(p-1)s_K^{p-2}\operatorname{tr}(D\boldsymbol A_K)\boldsymbol A_K,
\qquad s_K=\operatorname{tr}(\boldsymbol A_K).
\end{aligned}
\]
Star assembly and the frozen \(\omega_i^{\xi,(n)}\) give \(D\mathcal V_h^{(n)}\).
\normalsize

\bibliographystyle{abbrv}
\bibliography{references}
\end{document}